\documentclass[a4paper]{amsart}

\usepackage{color}
\usepackage[latin1]{inputenc}
\usepackage[T1]{fontenc}
\usepackage{amsfonts}
\usepackage{amssymb}
\usepackage{amsmath}
\usepackage{amsthm}
\usepackage{stmaryrd}
\usepackage{enumerate}
\usepackage{multirow}
\usepackage{graphicx}
\usepackage{graphicx,type1cm,eso-pic,color}
\usepackage{pstricks} 
\usepackage{float}
\newtheorem{theorem}{Theorem}[section]
\newtheorem{lemma}[theorem]{Lemma}
\newtheorem{proposition}[theorem]{Proposition}
\newtheorem{corollary}[theorem]{Corollary}

\newtheorem{assumption}[theorem]{Assumption}
\newtheorem{remark}[theorem]{Remark}

\begin{document}
\setlength\arraycolsep{2pt}
\title[Nadaraya-Watson Estimator for I.I.D. Paths of Diffusion Processes]{Nadaraya-Watson Estimator for I.I.D. Paths of Diffusion Processes}
\author{Nicolas MARIE$^{\dag}$}
\address{$^{\dag,\diamond}$Laboratoire Modal'X, Universit\'e Paris Nanterre, Nanterre, France}
\email{nmarie@parisnanterre.fr}
\author{Am\'elie ROSIER$^{\diamond}$}
\address{$^{\diamond}$ESME Sudria, Paris, France}
\email{amelie.rosier@esme.fr}
\keywords{Diffusion processes ; Nonparametric drift estimation ; Nadaraya-Watson estimator ; PCO method ; Cross validation.}
\date{}
\maketitle
%


%
\begin{abstract}
This paper deals with a nonparametric Nadaraya-Watson estimator $\widehat b$ of the drift function computed from independent continuous observations of a diffusion process. Risk bounds on $\widehat b$ and its discrete-time approximation are established. The paper also deals with extensions of the PCO and leave-one-out cross validation bandwidth selection methods for $\widehat b$. Finally, some numerical experiments are provided.
\end{abstract}
\textbf{MSC2020 subject classifications.} 62G05 ; 62M05.
%


%
\section{Introduction}\label{section_introduction}\label{section_introduction}
Consider the stochastic differential equation
\begin{equation}\label{main_equation}
X_t = x_0 +\int_{0}^{t}b(X_s)ds +\int_{0}^{t}\sigma(X_s)dW_s,
\end{equation}
where $b,\sigma :\mathbb R\rightarrow\mathbb R$ are two continuous functions and $W = (W_t)_{t\in [0,T]}$ is a Brownian motion.
\\
\\
Since the 1980's, the statistical inference for stochastic differential equations (SDE) has been widely investigated by many authors in the parametric and in the nonparametric frameworks. Classically (see Hoffmann \cite{HOFFMANN99}, Kessler \cite{KESSLER00}, Kutoyants \cite{KUTOYANTS04}, Dalalyan \cite{DALALYAN05}, Comte et al. \cite{CGCR07}, etc.), the estimators of the drift function are computed from one path of the solution to Equation (\ref{main_equation}) and converge when $T$ goes to infinity. The existence and the uniqueness of the stationary solution to Equation (\ref{main_equation}) are then required, and obtained thanks to restrictive conditions on $b$.
\\
\\
Let $\mathcal I : (x,w)\mapsto\mathcal I(x,w)$ be the It\^o map for Equation (\ref{main_equation}) and, for $N\in\mathbb N^*$ copies $W^1,\dots,W^N$ of $W$, consider $X^i =\mathcal I(x_0,W^i)$ for every $i\in\{1,\dots,N\}$. The estimation of the drift function $b$ from continuous-time and discrete-time observations of $(X^1,\dots,X^N)$ is a functional data analysis problem already investigated in the parametric framework (see Ditlevsen and De Gaetano \cite{DDG05}, Overgaard et al. \cite{OJTM05}, Picchini, De Gaetano and Ditlevsen \cite{PDGD10}, Picchini and Ditlevsen \cite{PD11}, Comte, Genon-Catalot and Samson \cite{CGCS13}, Delattre and Lavielle \cite{DL13}, Delattre, Genon-Catalot and Samson \cite{DGCS13}, Dion and Genon-Catalot \cite{DGC16}, Delattre, Genon-Catalot and Lar\'edo \cite{DGCL13}, etc.) and more recently in the nonparametric framework (see Comte and Genon-Catalot \cite{CGC20b,CGC21}, Della Maestra and Hoffmann \cite{DMH21} and Denis et al. \cite{DDM21}). In \cite{DMH21}, the authors also study a Nadaraya-Watson type estimator, presented bellow, of the drift function in McKean-Vlasov models.
\\
\\
Under the appropriate conditions on $b$ and $\sigma$ recalled at Section \ref{section_preliminaries}, the distribution of $X_t$ has a density $p_t(x_0,.)$ for every $t\in (0,T]$, and then one can define
\begin{displaymath}
f(x) :=\frac{1}{T - t_0}\int_{t_0}^{T}p_t(x_0,x)dt
\textrm{ $;$ }x\in\mathbb R
\end{displaymath}
for any $t_0 > 0$. Clearly, $f$ is a density function:
\begin{displaymath}
\int_{-\infty}^{\infty}f(x)dx =\frac{1}{T - t_0}\int_{t_0}^{T}\int_{-\infty}^{\infty}p_t(x_0,x)dxdt = 1.
\end{displaymath}
Let $K :\mathbb R\rightarrow\mathbb R$ be a kernel (i.e. an integrable function such that $\int K = 1$) and consider $K_h(x) := h^{-1}K(h^{-1}x)$ with $h\in (0,1]$. In the spirit of Comte and Genon-Catalot \cite{CGC20b,CGC21}, our paper deals first with the continuous-time Nadaraya-Watson estimator
\begin{equation}\label{NW_estimator}
\widehat b_{N,h}(x) :=
\frac{\widehat{bf}_{N,h}(x)}{\widehat f_{N,h}(x)}
\end{equation}
of the drift function $b$, where
\begin{equation}\label{NW_denominator}
\widehat f_{N,h}(x) :=
\frac{1}{N(T - t_0)}\sum_{i = 1}^{N}\int_{t_0}^{T}K_h(X_{t}^{i} - x)dt
\end{equation}
is an estimator of $f$ and
\begin{equation}\label{NW_numerator}
\widehat{bf}_{N,h}(x) :=
\frac{1}{N(T - t_0)}\sum_{i = 1}^{N}\int_{t_0}^{T}K_h(X_{t}^{i} - x)dX_{t}^{i}
\end{equation}
is an estimator of $bf$. From independent copies of $X$ continuously observed on $[0,T]$, $\widehat b_{N,h}$ is a natural extension of the Nadaraya-Watson estimator already well-studied in the regression framework (see Comte \cite{COMTE17}, Chapter 4 or Gy\"orfi et al. \cite{GKKW02}, Chapter 5). The paper also deals with a discrete-time approximate of the previous Nadaraya-Watson estimator:
\begin{equation}\label{NW_estimator_discrete}
\widehat b_{n,N,h}(x) :=
\frac{\widehat{bf}_{n,N,h}(x)}{\widehat f_{n,N,h}(x)},
\end{equation}
where
\begin{equation}\label{NW_denominator_discrete}
\widehat f_{n,N,h}(x) :=
\frac{1}{nN}\sum_{i = 1}^{N}\sum_{j = 0}^{n - 1}K_h(X_{t_j}^{i} - x)
\end{equation}
is an estimator of $f$,
\begin{equation}\label{NW_numerator_discrete}
\widehat{bf}_{n,N,h}(x) :=
\frac{1}{N(T - t_0)}\sum_{i = 1}^{N}\sum_{j = 0}^{n - 1}K_h(X_{t_j}^{i} - x)(X_{t_{j + 1}}^{i} - X_{t_j}^{i})
\end{equation}
is an estimator of $bf$, and $(t_0,t_1,\dots,t_n)$ is the dissection of $[t_0,T]$ such that $t_j = t_0 + (T - t_0)j/n$ for every $j\in\{1,\dots,n\}$. Finally, our paper deals with a risk bound on the adaptive double bandwidths Nadaraya-Watson's estimator
\begin{displaymath}
\widehat b_{N,\widehat h,\widehat h'}(x) :=
\frac{\widehat{bf}_{N,\widehat h}(x)}{\widehat f_{N,\widehat h'}(x)},
\end{displaymath}
where $\widehat h$ (resp. $\widehat h'$) is selected via a penalized comparison to overfitting (PCO) type criterion for its numerator (resp. denominator). However, in the nonparametric regression framework, it is established in Comte and Marie \cite{CM21} that the leave-one-out cross-validation (looCV) bandwidth selection method for Nadaraya-Watson's estimator is numerically more satisfactory than two alternative procedures based on Goldenshluger-Lepski's method and on the PCO method. For this reason, an extension of the looCV method to $\widehat b_{n,N,h}$ is also provided, with numerical experiments, even if it seems difficult to establish a risk bound on the associated adaptive estimator.
\\
\\
Now, let us compare $\widehat b_{N,h}$ with the estimator of Della Maestra and Hoffmann \cite{DMH21} restricted to our framework:
\begin{displaymath}
\widehat b_{N,\mathbf h}(x) :=
\frac{\displaystyle{\sum_{i = 1}^{N}\int_{0}^{T}L_{h_1}(\tau - t)K_{h_2}(X_{t}^{i} - x)dX_{t}^{i}}}{\displaystyle{\sum_{i = 1}^{N}K_{h_3}(X_{\tau}^{i} - x)}},
\end{displaymath}
where $L :\mathbb R\rightarrow\mathbb R$ is another kernel, $h_1,h_2,h_3\in (0,1]$, $\mathbf h = (h_1,h_2,h_3)$ and $\tau\in (0,T)$. In \cite{DMH21}, the authors provide a nice risk bound on the adaptive estimator obtained by selecting $(h_1,h_2)$ (resp. $h_3$) via a Goldenshluger-Lepski type procedure on the numerator (resp. the denominator) of $\widehat b_{N,\mathbf h}$. As mentioned above, in Comte and Marie \cite{CM21}, it has been established that in the nonparametric regression framework, this approach is numerically less satisfactory than the looCV method. However, the looCV method provided in our paper doesn't extend to $\widehat b_{N,\mathbf h}$ because it cannot be written easily as a linear combination. Note also that even if it is numerically less satisfactory than the looCV method provided at Subsection \ref{section_looCV}, the PCO type method provided in our paper at Subsection \ref{section_PCO} is easier to implement and numerically faster than a Goldenshluger-Lepski type method because, as in the nonparametric regression framework, the criterion to minimize depends on one variable instead of two, and because there is no constant to calibrate. For technical reasons explained at Section \ref{section_conclusion}, the condition $(Nh^3)^{-1}\leqslant 1$ is required on the bandwidths collection to establish a risk bound on our PCO based adaptive estimator of $bf$, when Della Maestra and Hoffmann only need the condition $\log(N)^2(Nh)^{-1}\leqslant 1$ to establish a risk bound on their Goldenshluger-Lepski based adaptive estimator of $bp_{\tau}(x_0,.)$ in \cite{DMH21}. However, Remark \ref{uncommon_condition_bandwidths} explains why the condition $(Nh^3)^{-1}\leqslant 1$ on the bandwidths collection is not that uncomfortable. Finally, under similar conditions on $b$, $\sigma$ and $K$, the rate of convergence of our continuous-time Nadaraya-Watson estimator is of same order than the rate of convergence of the estimator of Della Maestra and Hoffmann \cite{DMH21} in the nonadaptive case. There is no discrete-time approximate of $\widehat b_{N,\mathbf h}$ studied in \cite{DMH21}.
\\
\\
Finally, even if they deal with a different type of nonparametric estimators, let us say few words on the recent papers of Comte and Genon-Catalot \cite{CGC20b} and Denis et al. \cite{DDM21}. On the one hand, in \cite{CGC20b}, the authors extend to the diffusion processes framework, for continuous-time observations, the least squares projection estimator already well studied in the regression framework (see Cohen et al. \cite{CDL13}, Comte and Genon-Catalot \cite{CGC20a}, etc.). In particular, they provide a model selection procedure and establish a risk bound on the associated adaptive estimator. As explained at Section \ref{section_continuous_time_NW}, in the nonadaptive case, the variance term of their estimator is comparable with the variance term of $\widehat b_{N,h}$, but as in the nonparametric regression framework, the rate of convergence of the least squares projection estimator depends on the regularity space associated to the projection basis. On the other hand, in \cite{DDM21}, the authors focus on a projection least squares estimator computed from discrete-time observations and on a $B$-spline space. They provide a model selection procedure and prove both upper and lower bounds on the associated adaptive estimator.
\\
\\
Section \ref{section_preliminaries} deals with the existence and the regularity of the density $p_t(x_0,.)$ of $X_t$ for every $t\in (0,T]$, and with a Nikol'skii type condition fulfilled by $f$. Section \ref{section_continuous_time_NW} deals with a risk bound on the continuous-time Nadaraya-Watson estimator and Section \ref{section_discrete_time_NW} with a risk bound on its discrete-time approximate. Finally, Section \ref{section_looCV_simulations} provides extensions of the PCO and looCV methods for the Nadaraya-Watson estimator studied in this paper. Some numerical experiments on the looCV based adaptive Nadaraya-Watson estimator are also provided. The proofs are postponed to Appendix \ref{section_proofs}.
\\
\\
\textbf{Notations and basic definitions:}
\begin{itemize}
 \item For every $A,B\in\mathbb R$ such that $A < B$, $C^0([A,B];\mathbb R)$ is equipped with the uniform norm $\|.\|_{\infty,A,B}$, and $C^0(\mathbb R)$ is equipped with the uniform (semi-)norm $\|.\|_{\infty}$.
 \item For every $p\in\overline{\mathbb N}$, $C_{\rm b}^{p}(\mathbb R) :=\cap_{j = 0}^{p}\{\varphi\in C^p(\mathbb R) :\varphi^{(j)}\textrm{ is bounded}\}$.
 \item For every $p\geqslant 1$, $\mathbb L^p(\mathbb R,dx)$ is equipped with its usual norm $\|.\|_p$ such that
 \begin{displaymath}
 \|\varphi\|_p :=\left(\int_{-\infty}^{\infty}\varphi(x)^pdx\right)^{1/p}
 \textrm{$;$ }
 \forall\varphi\in\mathbb L^p(\mathbb R,dx).
 \end{displaymath}
 \item $\mathbb H^2$ is the space of the processes $(Y_t)_{t\in [0,T]}$, adapted to the filtration generated by $W$, such that
 \begin{displaymath}
 \int_{0}^{T}\mathbb E(Y_{t}^{2})dt <\infty.
 \end{displaymath}
 \item For a given kernel $\delta$, the usual scalar product on $\mathbb L^2(\mathbb R,\delta(x)dx)$ is denoted by $\langle .,.\rangle_{2,\delta}$ and the associated norm by $\|.\|_{2,\delta}$.
\end{itemize}
%


%
\section{Preliminaries: regularity of the density and estimates}\label{section_preliminaries}
This section deals with the existence and the regularity of the density $p_t(x_0,.)$ of $X_t$ for every $t\in (0,T]$, with the Kusuoka-Stroock bounds on $(t,x)\mapsto p_t(x_0,x)$ and its derivatives, and then with a Nikol'skii type condition fulfilled by $f$.
\\
\\
In the sequel, in order to ensure the existence and the uniqueness of the (strong) solution to Equation (\ref{main_equation}), $b$ and $\sigma$ fulfill the following regularity assumption.
%


%
\begin{assumption}\label{assumption_b_sigma}
The functions $b$ and $\sigma$ are Lipschitz continuous.
\end{assumption}
\noindent
Now, assume that the solution $X$ to Equation (\ref{main_equation}) fulfills the following assumption.
%


%
\begin{assumption}\label{assumption_bound_derivative_density_x}
There exists $\beta\in\mathbb N^*$ such that, for any $t\in (0,T]$, the distribution of $X_t$ has a $\beta$ times continuously derivable density $p_t(x_0,.)$. Moreover, for every $x\in\mathbb R$,
\begin{displaymath}
0 < p_t(x_0,x)\leqslant
\frac{\mathfrak c_{\ref{assumption_bound_derivative_density_x},1}}{t^{1/2}}\exp\left[-\mathfrak 
m_{\ref{assumption_bound_derivative_density_x},1}\frac{(x - x_0)^2}{t}\right]
\end{displaymath}
and
\begin{displaymath}
|\partial_{x}^{\ell}p_t(x_0,x)|\leqslant
\frac{\mathfrak c_{\ref{assumption_bound_derivative_density_x},2}(\ell)}{t^{q_2(\ell)}}\exp\left[-\mathfrak m_{\ref{assumption_bound_derivative_density_x},2}(\ell)\frac{(x - x_0)^2}{t}\right]
\textrm{ $;$ }
\forall\ell\in\{1,\dots,\beta\},
\end{displaymath}
where all the constants are positive, depend on $T$, but not on $t$ and $x$.
\end{assumption}
\noindent
At Section \ref{section_discrete_time_NW}, the following assumption on $X$ is also required.
%


%
\begin{assumption}\label{assumption_bound_derivative_density_t}
For any $x\in\mathbb R$, the function $t\in (0,T]\mapsto p_t(x_0,x)$ is continuously derivable. Moreover,
\begin{displaymath}
|\partial_tp_t(x_0,x)|
\leqslant
\frac{\mathfrak c_{\ref{assumption_bound_derivative_density_t},3}}{t^{q_3}}\exp\left[-\mathfrak  m_{\ref{assumption_bound_derivative_density_t},3}\frac{(x - x_0)^2}{t}\right]
\textrm{ $;$ }
\forall t\in (0,T],
\end{displaymath}
where $\mathfrak c_{\ref{assumption_bound_derivative_density_t},3}$, $\mathfrak  m_{\ref{assumption_bound_derivative_density_t},3}$ and $q_3$ are three positive constants depending on $T$ but not on $t$ and $x$.
\end{assumption}
\noindent
Let us provide some examples of diffusion processes categories satisfying Assumptions \ref{assumption_bound_derivative_density_x} and/or \ref{assumption_bound_derivative_density_t}.
\\
\\
\textbf{Examples:}
\begin{enumerate}
 \item Assume that the functions $b$ and $\sigma$ belong to $C_{\rm b}^{\infty}(\mathbb R)$, and that there exists $\alpha > 0$ such that
 \begin{equation}\label{non_degeneracy}
 |\sigma(x)| >\alpha
 \textrm{ $;$ }
 \forall x\in\mathbb R.
 \end{equation}
 Then, by Kusuoka and Stroock \cite{KS85}, Corollary 3.25, $X$ fulfills Assumptions \ref{assumption_bound_derivative_density_x} and \ref{assumption_bound_derivative_density_t}.
 \item Assume that $b$ is Lipschitz continuous (but not bounded) and that $\sigma\in C_{\rm b}^{1}(\mathbb R)$. Assume also that $\sigma$ satisfies the non-degeneracy condition (\ref{non_degeneracy}) and that $\sigma'$ is H\"older continuous. Then, by Menozzi et al. \cite{MPZ21}, Theorem 1.2, $X$ fulfills Assumption \ref{assumption_bound_derivative_density_x} with $\beta = 1$ (but not necessarily Assumption \ref{assumption_bound_derivative_density_t}). Note that the conditions required to apply Menozzi et al. \cite{MPZ21}, Theorem 1.2 are fulfilled by the so-called Ornstein-Uhlenbeck process, that is the solution to the Langevin equation:
 \begin{equation}\label{Langevin}
 X_t = x_0 -\theta\int_{0}^{t}X_sds +\sigma W_t
 \textrm{ $;$ }
 t\in\mathbb R_+,
 \end{equation}
 where $\theta,\sigma > 0$ and $x_0\in\mathbb R_+$. In this special case, since it is well-known that the solution to Equation (\ref{Langevin}) is a Gaussian process such that
 \begin{displaymath}
\mathbb E(X_t) = x_0e^{-\theta t}
 \quad\textrm{and}\quad
{\rm var}(X_t) =
 \frac{\sigma^2}{2\theta}(1 - e^{-2\theta t})
 \textrm{ $;$ }\forall t\in [0,T],
 \end{displaymath}
 one can show that $X$ also fulfills Assumption \ref{assumption_bound_derivative_density_t}.
\end{enumerate}
%


%
\begin{remark}\label{integrability_functionals_X}
Under Assumptions \ref{assumption_b_sigma} and \ref{assumption_bound_derivative_density_x}, for any $p\geqslant 1$ and any continuous function $\varphi :\mathbb R\rightarrow\mathbb R$ having polynomial growth, $t\in [0,T]\mapsto\mathbb E(|\varphi(X_t)|^p)$ is bounded. Indeed, for any $t\in [0,T]$,
\begin{eqnarray*}
 \mathbb E(|\varphi(X_t)|^p) & \leqslant &
 \mathfrak c_1(1 +\mathbb E(|X_t|^{pq})) =
 \mathfrak c_1\int_{-\infty}^{\infty}(1 + |x|^{pq})p_t(x_0,x)dx\\
 & \leqslant &
 \mathfrak c_1\mathfrak c_{\ref{assumption_bound_derivative_density_x},1}
 \int_{-\infty}^{\infty}(1 + |t^{1/2}x + x_0|^{pq})e^{
 -\mathfrak m_{\ref{assumption_bound_derivative_density_x},1}x^2}dx
 \leqslant
 \mathfrak c_2(1\vee T^{pq/2})
\end{eqnarray*}
where
\begin{displaymath}
\mathfrak c_2 =
\mathfrak c_1\mathfrak c_{\ref{assumption_bound_derivative_density_x},1}
\int_{-\infty}^{\infty}[1 + (|x| + |x_0|)^{pq}]e^{
-\mathfrak m_{\ref{assumption_bound_derivative_density_x},1}x^2}dx
\end{displaymath}
and the constants $\mathfrak c_1,q > 0$ only depend on $\varphi$. Moreover,
\begin{eqnarray*}
 \|\varphi\|_{p,f} & := & \int_{-\infty}^{\infty}|\varphi(x)|^pf(x)dx\\
 & = &
 \frac{1}{T - t_0}\int_{t_0}^{T}
 \mathbb E(|\varphi(X_t)|^p)dt\leqslant\mathfrak c_2(1\vee T^{pq/2}).
\end{eqnarray*}
Then, $\varphi\in\mathbb L^p(\mathbb R,f(x)dx)$ and $\|\varphi\|_{p,f}$ is bounded by a constant which doesn't depend on $t_0$. In particular, the remark applies to $b$ and $\sigma$ with $q = 1$ by Assumption \ref{assumption_b_sigma}.
\end{remark}
\noindent
Finally, let us show that $f$ fulfills a Nikol'skii type condition.
%


%
\begin{corollary}\label{Nikolskii_f}
Under Assumption \ref{assumption_b_sigma}, $f(x) > 0$ for every $x\in\mathbb R$. Moreover, under Assumptions \ref{assumption_b_sigma} and \ref{assumption_bound_derivative_density_x}, there exists $\mathfrak c_{\ref{Nikolskii_f}} > 0$, depending on $T$ but not on $t_0$, such that for every $\ell\in\{0,\dots,\beta - 1\}$ and $\theta\in\mathbb R$,
\begin{displaymath}
\int_{-\infty}^{\infty}[f^{(\ell)}(x +\theta) - f^{(\ell)}(x)]^2dx
\leqslant\frac{\mathfrak c_{\ref{Nikolskii_f}}}{t_{0}^{2q_2(\ell + 1)}}(\theta^2 + |\theta|^3).
\end{displaymath}
\end{corollary}
%


%
\begin{remark}\label{how_to_choose_t_0}
Assumption \ref{assumption_bound_derivative_density_x}, Assumption \ref{assumption_bound_derivative_density_t} and Corollary \ref{Nikolskii_f} are crucial in the sequel, but $t_0$ has to be chosen carefully to get reasonable risk bounds on the estimators $\widehat b_{N,h}$ and $\widehat b_{n,N,h}$. Indeed, the behavior of the Kusuoka-Stroock bounds on $(t,x)\mapsto p_t(x_0,x)$ and its derivatives is singular at point $(0,x_0)$. This is due to the fact that the distribution of $X$ at time $0$ is a Dirac measure while that it has a smooth density with respect to Lebesgue's measure for every $t\in (0,T]$. Moreover, since $X$ is not a stationary process in general, the Kusuoka-Stroock bounds on $(t,x)\mapsto p_t(x_0,x)$ and its derivatives explode when $T\rightarrow\infty$. The same difficulty appears with the estimators studied in Comte and Genon-Catalot \cite{CGC20b} and in Della Maestra and Hoffmann \cite{DMH21}. So, it is recommended to take $T$ as small as possible in practice. In the sequel, only the dependence in $t_0$ is tracked in the risk bounds derived from Assumption \ref{assumption_bound_derivative_density_x}, Assumption \ref{assumption_bound_derivative_density_t} and Corollary \ref{Nikolskii_f} because it is specific to our approach. Finally, these risk bounds only depend on $t_0$ through a multiplicative constant of order $1/\min\{t_{0}^{\alpha},T - t_0\}$ ($\alpha > 0$). So, to take $t_0\in [1,T - 1]$ when $T > 1$ gives constants not depending on $t_0$.
\end{remark}
%


%
\section{Risk bound on the continuous-time Nadaraya-Watson estimator}\label{section_continuous_time_NW}
This section deals with risk bounds on $\widehat f_{N,h}$, $\widehat{bf}_{N,h}$, and then on the Nadaraya-Watson estimator $\widehat b_{N,h}$.
\\
\\
In the sequel, the kernel $K$ fulfills the following usual assumptions.
%


%
\begin{assumption}\label{assumption_K_1}
The kernel $K$ is symmetric, continuous and belongs to $\mathbb L^2(\mathbb R,dx)$.
\end{assumption}
%


%
\begin{assumption}\label{assumption_K_2}
There exists $\upsilon\in\mathbb N^*$ such that
\begin{displaymath}
\int_{-\infty}^{\infty}|z^{\upsilon + 1}K(z)|dz <\infty
\quad\textrm{and}\quad
\int_{-\infty}^{\infty}z^{\ell}K(z)dz = 0
\textrm{ $;$ }
\forall\ell\in\{1,\dots,\upsilon\}.
\end{displaymath}
\end{assumption}
\noindent
About the construction of kernels fulfilling both Assumptions \ref{assumption_K_1} and \ref{assumption_K_2}, the reader can refer to Comte \cite{COMTE17}, Proposition 2.10. The following proposition provides a risk bound on $\widehat f_{N,h}$ (see (\ref{NW_denominator})).
%


%
\begin{proposition}\label{risk_bound_denominator}
Under Assumptions \ref{assumption_b_sigma}, \ref{assumption_bound_derivative_density_x}, \ref{assumption_K_1} and \ref{assumption_K_2} with $\upsilon =\beta$,
\begin{displaymath}
\mathbb E(\|\widehat f_{N,h} - f\|_{2}^{2})
\leqslant
\mathfrak c_{\ref{risk_bound_denominator}}(t_0)h^{2\beta} +\frac{\|K\|_{2}^{2}}{Nh}
\end{displaymath}
with
\begin{displaymath}
\mathfrak c_{\ref{risk_bound_denominator}}(t_0) =
\frac{\mathfrak c_{\ref{Nikolskii_f}}}{|(\beta - 2)!|^2t_{0}^{2q_2(\beta)}}
\left(\int_{-\infty}^{\infty}|z|^{\beta}(1 + |z|^{1/2})|K(z)|dz\right)^2.
\end{displaymath}
\end{proposition}
\noindent
Note that thanks to Proposition \ref{risk_bound_denominator}, the bias-variance tradeoff is reached by (the risk bound on) $\widehat f_{N,h}$ when $h$ is of order $N^{-1/(2\beta + 1)}$, leading to a rate of order $N^{-2\beta/(2\beta + 1)}$. Moreover, by Remark \ref{how_to_choose_t_0}, to take $t_0\geqslant 1$ when $T > 1$ gives
\begin{displaymath}
\mathbb E(\|\widehat f_{N,h} - f\|_{2}^{2})
\leqslant
\mathfrak c_{\ref{risk_bound_denominator}}h^{2\beta}  +\frac{\|K\|_{2}^{2}}{Nh}
\quad {\rm with}\quad
\mathfrak c_{\ref{risk_bound_denominator}} =
\frac{\mathfrak c_{\ref{Nikolskii_f}}}{|(\beta - 2)!|^2}
\left(\int_{-\infty}^{\infty}|z|^{\beta}(1 + |z|^{1/2})|K(z)|dz\right)^2.
\end{displaymath}
Note also that in the risk bound on $\widehat f_{N,h}$ of Proposition \ref{risk_bound_denominator}, only the control of the bias term depends on $T$, trough the constant $\mathfrak c_{\ref{Nikolskii_f}}$, depending itself on the constants $\mathfrak c_{\ref{assumption_bound_derivative_density_x},2}(\ell)$, $\ell\in\{1,\dots,\beta\}$, involved in the Kusuoka-Strook bounds (see Assumption \ref{assumption_bound_derivative_density_x}). Indeed, except in the special case of the Ornstein-Uhlenbeck process which is stationary, for all the examples of diffusion processes fulfilling Assumption \ref{assumption_bound_derivative_density_x} (see Kusuoka and Stroock \cite{KS85} and Menozzi et al. \cite{MPZ21}), the constants  $\mathfrak c_{\ref{assumption_bound_derivative_density_x},2}(\ell)$, $\ell\in\{1,\dots,\beta\}$, depend on $T$. The variance term doesn't depend on time at all.
\\
\\
In the sequel, $\mathbb L^2(\mathbb R,f(x)dx)$ is equipped with the $f$-weighted norm $\|.\|_{2,f}$ defined at the end of the introduction section. Let us recall that by Remark \ref{integrability_functionals_X}, for every $\varphi\in\mathbb L^2(\mathbb R,f(x)dx)$, $\|\varphi\|_{2,f}$ is bounded by a constant which doesn't depend on $t_0$.
\\
\\
The following proposition provides a risk bound on $\widehat{bf}_{N,h}$ (see (\ref{NW_numerator})).
%


%
\begin{proposition}\label{risk_bound_numerator}
Under Assumptions \ref{assumption_b_sigma} and \ref{assumption_K_1},
\begin{displaymath}
\mathbb E(\|\widehat{bf}_{N,h} - bf\|_{2}^{2})
\leqslant
\|(bf)_h - bf\|_{2}^{2} +
\frac{\mathfrak c_{\ref{risk_bound_numerator}}(t_0)}{Nh}
\end{displaymath}
with $(bf)_h := K_h\ast (bf)$ and
\begin{displaymath}
\mathfrak c_{\ref{risk_bound_numerator}}(t_0) =
2\|K\|_{2}^2
\left(\|b\|_{2,f}^{2} + \frac{1}{T - t_0}\|\sigma\|_{2,f}^{2}\right).
\end{displaymath}
\end{proposition}
\noindent
Assume that $bf$ is $\gamma\in\mathbb N^*$ times continuously derivable and that there exists $\varphi\in\mathbb L^1(\mathbb R,|z|^{\gamma - 1}K(z)dz)$ such that, for every $\theta\in\mathbb R$ and $h\in (0,1]$,
\begin{equation}\label{Nikolskii_type_condition}
\int_{-\infty}^{\infty}[(bf)^{(\gamma - 1)}(x + h\theta) - (bf)^{(\gamma - 1)}(x)]^2dx
\leqslant\varphi(\theta)h^2.
\end{equation}
If in addition $K$ fulfills Assumption \ref{assumption_K_2} with $\upsilon =\gamma$, then $\|(bf)_h - bf\|_{2}^{2}$ is of order $h^{2\gamma}$, and by Proposition \ref{risk_bound_numerator}, the bias-variance tradeoff is reached by $\widehat{bf}_{N,h}$ when $h$ is of order $N^{-1/(2\gamma + 1)}$, leading to the rate $N^{-2\gamma/(2\gamma + 1)}$. Moreover, by Remark \ref{how_to_choose_t_0}, to take $t_0\leqslant T - 1$ when $T > 1$ gives
\begin{displaymath}
\mathbb E(\|\widehat{bf}_{N,h} - bf\|_{2}^{2})
\leqslant
\|(bf)_h - bf\|_{2}^{2} +
\frac{\mathfrak c_{\ref{risk_bound_numerator}}}{Nh}
\quad {\rm with}\quad
\mathfrak c_{\ref{risk_bound_numerator}} =
2\|K\|_{2}^2
(\|b\|_{2,f}^{2} +\|\sigma\|_{2,f}^{2}).
\end{displaymath}
Note also that the variance term in this risk bound doesn't depend on $T$.
\\
\\
Finally, Propositions \ref{risk_bound_denominator} and \ref{risk_bound_numerator} allow to provide a risk bound on a truncated version of the Nadaraya-Watson estimator $\widehat b_{N,h}$ (see (\ref{NW_estimator})).
%


%
\begin{proposition}\label{risk_bound_NW}
Consider the 2 bandwidths (truncated) Nadaraya-Watson (2bNW) estimator
\begin{displaymath}
\widehat b_{N,h,h'}(x) :=
\frac{\widehat{bf}_{N,h}(x)}{\widehat f_{N,h'}(x)}
\mathbf 1_{\widehat f_{N,h'}(x) > m/2}
\quad
\textrm{with}
\quad
h,h' > 0,
\end{displaymath}
and assume that $f(x) > m > 0$ for every $x\in [A,B]$ ($m\in (0,1]$ and $A,B\in\mathbb R$ such that $A < B$). Under Assumptions \ref{assumption_b_sigma}, \ref{assumption_bound_derivative_density_x}, \ref{assumption_K_1} and \ref{assumption_K_2} with $\upsilon =\beta$,
\begin{displaymath}
\mathbb E(\|\widehat b_{N,h,h'} - b\|_{f,A,B}^{2})
\leqslant
\frac{\mathfrak c_{\ref{risk_bound_NW}}}{m^2}\left[\|(bf)_h - bf\|_{2}^{2} +
\frac{\mathfrak c_{\ref{risk_bound_numerator}}(t_0)}{Nh} +
2\|b\|_{2,f}^{2}\left(\mathfrak c_{\ref{risk_bound_denominator}}(t_0)(h')^{2\beta} +
\frac{\|K\|_{2}^{2}}{Nh'}\right)\right]
\end{displaymath}
with $\mathfrak c_{\ref{risk_bound_NW}} := 8(\|f\|_{\infty}\vee\|b^2f\|_{\infty})$ and $\|\varphi\|_{f,A,B} :=\|\varphi\mathbf 1_{[A,B]}\|_{2,f}$ for every $\varphi\in\mathbb L^2(\mathbb R,f(x)dx)$.
\end{proposition}
\noindent
Proposition \ref{risk_bound_NW} says that the risk of $\widehat b_{N,h,h'}$ can be controlled by the sum of those of $\widehat{bf}_{N,h}$ and $\widehat f_{N,h'}$ up to a multiplicative constant. Now, if $K$ fulfills Assumption \ref{assumption_K_2} with $\upsilon =\beta\vee\gamma$, and if $bf$ satisfies Condition (\ref{Nikolskii_type_condition}), then the risk bound on $\widehat b_{N,h,h'}$ is of order $h^{2\gamma} + (h')^{2\beta} + 1/(Nh) + 1/(Nh')$, and the bias-variance tradeoff is reached when $h$ (resp. $h'$) is of order $N^{-1/(2\gamma + 1)}$ (resp. $N^{-1/(2\beta + 1)}$), leading to the rate
\begin{displaymath}
N^{-2\left[\left(\frac{\gamma}{2\gamma + 1}\right)\wedge
\left(\frac{\beta}{2\beta + 1}\right)\right]} =
N^{-\frac{2(\beta\wedge\gamma)}{2(\beta\wedge\gamma) + 1}},
\end{displaymath}
which is of same order than the rate of the nonadaptive version of the estimator of Della Maestra and Hoffmann \cite{DMH21} (see their Theorem 15). Note also that to consider the 2bNW estimator is crucial to extend the PCO method to our framework in the spirit of Comte and Marie \cite{CM21} (see Subsection \ref{section_PCO}). However, by taking $h = h'$ of order $N^{-1/(2(\beta\wedge\gamma) + 1)}$, the bias-variance tradeoff is reached by the $1$ bandwidth (truncated) Nadaraya-Watson estimator with the same rate. Finally, if $h = h'$ and $\sigma$ is bounded, then the variance term in the risk bound of Proposition \ref{risk_bound_NW} is comparable to the variance term in the risk bound obtained by Comte and Genon-Catalot in \cite{CGC20b} for their least squares projection estimator (see \cite{CGC20b}, Propositions 2.1 and 2.2). Indeed, for a $d$-dimensional projection space, the variance term in the risk bound of Comte and Genon-Catalot \cite{CGC20b}, Proposition 2.1 is of order $d/N$ which is comparable to $1/(Nh)$. The rate of convergence of their least squares projection estimator depends on the regularity space associated to the projection basis but, as in the nonparametric regression framework, not on the regularity of $f$.\\
The limitation of our Proposition \ref{risk_bound_NW} is that $m$ is unknown in general and must be replaced by an estimator as well. Most of the time, as stated in Comte \cite{COMTE17}, Chapter 4, the minimum of an estimator of $f$ is taken to choose $m$ in practice:
\begin{displaymath}
\widehat m_{N,h'} =\min\{\widehat f_{N,h'}(x)\textrm{ $;$ }x\in [A,B]\}
\end{displaymath}
for instance. A more naive way to solve this difficulty in practice is to take
\begin{displaymath}
m = m_N =
\mathfrak cN^{-\frac{\varepsilon}{2}\cdot\frac{2(\beta\wedge\gamma)}{2(\beta\wedge\gamma) + 1}}
\xrightarrow[N\rightarrow\infty]{} 0,
\end{displaymath}
where $\mathfrak c > 0$ is a fixed constant and $\varepsilon\in (0,1)$ is chosen as close as possible to $0$. Under Assumption \ref{assumption_bound_derivative_density_x}, by Corollary \ref{Nikolskii_f},
\begin{displaymath}
\exists N_0\in\mathbb N :
\forall N > N_0\textrm{, }
\forall x\in [A,B]\textrm{, }
f(x) > m_N.
\end{displaymath}
So, by Proposition \ref{risk_bound_NW}, when $h$ (resp. $h'$) is of order $N^{-1/(2\gamma + 1)}$ (resp. $N^{-1/(2\beta + 1)}$), $\widehat b_{N,h,h'}$ converges with the slightly degraded rate
\begin{displaymath}
N^{-(1 -\varepsilon)\frac{2(\beta\wedge\gamma)}{2(\beta\wedge\gamma) + 1}}.
\end{displaymath}
This last comment remains true for Proposition \ref{risk_bound_approximate_NW} and Corollary \ref{risk_bound_NW_PCO}.
%


%
\section{Risk bound on the discrete-time approximate Nadaraya-Watson estimator}\label{section_discrete_time_NW}
This section deals with risk bounds on $\widehat f_{n,N,h}$, $\widehat{bf}_{n,N,h}$, and then on the approximate Nadaraya-Watson estimator $\widehat b_{n,N,h}$.
\\
\\
In the sequel, in addition to Assumptions \ref{assumption_K_1} and \ref{assumption_K_2}, $K$ fulfills the following one.
%


%
\begin{assumption}\label{assumption_K_+}
The kernel $K$ is two times continuously derivable on $\mathbb R$ and $K',K''\in\mathbb L^2(\mathbb R,dx)$.
\end{assumption}
\noindent
Compactly supported kernels belonging to $C^2(\mathbb R)$ or Gaussian kernels fulfill Assumption \ref{assumption_K_+}. The following proposition provides a risk bound on $\widehat f_{n,N,h}$ (see (\ref{NW_denominator_discrete})).
%


%
\begin{proposition}\label{risk_bound_approximate_denominator}
Under Assumptions \ref{assumption_b_sigma}, \ref{assumption_bound_derivative_density_x}, \ref{assumption_bound_derivative_density_t}, \ref{assumption_K_1}, \ref{assumption_K_2} with $\upsilon =\beta$, and \ref{assumption_K_+}, if
\begin{displaymath}
\frac{1}{nh^2}\leqslant 1,
\end{displaymath}
then there exists a constant $\mathfrak c_{\ref{risk_bound_approximate_denominator}} > 0$, not depending on $h$, $N$, $n$ and $t_0$, such that
\begin{displaymath}
\mathbb E(\|\widehat f_{n,N,h} - f\|_{2}^{2})
\leqslant
\frac{\mathfrak c_{\ref{risk_bound_approximate_denominator}}}{\min\{t_{0}^{2q_2(\beta)},t_{0}^{2q_3}\}}
\left(h^{2\beta} +\frac{1}{Nh} +\frac{1}{n^2}\right) +\frac{1}{Nnh^3}.
\end{displaymath}
\end{proposition}
\noindent
Assume that $\beta = 1$ (extreme case) and $h$ is of order $N^{-1/3}$. As mentioned at Section \ref{section_continuous_time_NW}, under this condition, the bias-variance tradeoff is reached by the continuous-time estimator of $f$. Then, the approximation error of $\widehat f_{n,N,h}$ is of order $1/n$, which is the order of the variance of the Brownian motion increments along the dissection $(t_0,t_1,\dots,t_n)$ of $[t_0,T]$. For this reason, the risk bound established in Proposition \ref{risk_bound_approximate_denominator} is satisfactory. Moreover, by Remark \ref{how_to_choose_t_0}, to take $t_0\geqslant 1$ when $T > 1$ gives
\begin{displaymath}
\mathbb E(\|\widehat f_{n,N,h} - f\|_{2}^{2})
\leqslant
\mathfrak c_{\ref{risk_bound_approximate_denominator}}
\left(h^{2\beta} +\frac{1}{Nh} +\frac{1}{n^2}\right) +\frac{1}{Nnh^3}.
\end{displaymath}
The following proposition provides a risk bound on $\widehat{bf}_{n,N,h}$ (see (\ref{NW_numerator_discrete})).
%


%
\begin{proposition}\label{risk_bound_approximate_numerator}
Consider $\varepsilon\in (0,1)$. Under Assumptions \ref{assumption_b_sigma}, \ref{assumption_bound_derivative_density_x}, \ref{assumption_bound_derivative_density_t}, \ref{assumption_K_1} and \ref{assumption_K_+}, if
\begin{displaymath}
\frac{1}{nh^{2 -\varepsilon}}\leqslant 1,
\end{displaymath}
the kernel $K$ belongs to $\mathbb L^4(\mathbb R,dx)$ and $z\mapsto zK'(z)$ belongs to $\mathbb L^2(\mathbb R,dx)$, then there exist a constant $\mathfrak c_{\ref{risk_bound_approximate_numerator}} > 0$, not depending on $\varepsilon$, $h$, $N$, $n$ and $t_0$, and a constant $\mathfrak c_{\ref{risk_bound_approximate_numerator}}(\varepsilon) > 0$, depending on $\varepsilon$ but not on $h$, $N$, $n$ and $t_0$, such that
\begin{eqnarray*}
 \mathbb E(\|\widehat{bf}_{n,N,h} - bf\|_{2}^{2})
 & \leqslant &
 \frac{\mathfrak c_{\ref{risk_bound_approximate_numerator}}}{\min\{t_{0}^{1/2},t_{0}^{2q_3},T - t_0\}}
 \left(\|(bf)_h - bf\|_{2}^{2} +\frac{1}{Nh} +\frac{1}{n}\right)\\
 & &
 \hspace{6cm}
 +\frac{\mathfrak c_{\ref{risk_bound_approximate_numerator}}(\varepsilon)}{\min\{1,t_{0}^{(1 -\varepsilon)/2}\}}
 \cdot\frac{1}{Nnh^{3 +\varepsilon}}.
\end{eqnarray*}
\end{proposition}
%


%
\begin{remark}\label{special_case_approximate_numerator}
Note that if $b$ and $\sigma$ are bounded, Proposition \ref{risk_bound_approximate_numerator} can be improved. Precisely, with $\varepsilon = 0$ and without the additional conditions $K\in\mathbb L^4(\mathbb R,dx)$ and $z\mapsto zK'(z)$ belongs to $\mathbb L^2(\mathbb R,dx)$, the risk bound on $\widehat{bf}_{n,N,h}$ is of same order than in Proposition \ref{risk_bound_approximate_denominator} (see Remark \ref{special_case_approximate_numerator_details} for details).
\end{remark}
\noindent
Assume that $bf$ fulfills Condition (\ref{Nikolskii_type_condition}) with $\gamma = 1$ (extreme case), and that $h$ is of order $N^{-1/3}$. Then, for $\varepsilon > 0$ as close as possible to $0$, the approximation error of $\widehat{bf}_{n,N,h}$ is of order $N^{\varepsilon/3}/n$. If in addition $b$ and $\sigma$ are bounded, thanks to Remark \ref{special_case_approximate_numerator}, with $\varepsilon = 0$ and without the additional conditions $K\in\mathbb L^4(\mathbb R,dx)$ and $z\mapsto zK'(z)$ belongs to $\mathbb L^2(\mathbb R,dx)$, then the approximation error of $\widehat{bf}_{n,N,h}$ is of order $1/n$ as the error of $\widehat f_{n,N,h}$. Moreover, by Remark \ref{how_to_choose_t_0}, to take $t_0\in [1,T - 1]$ when $T > 1$ gives
\begin{displaymath}
\mathbb E(\|\widehat{bf}_{n,N,h} - bf\|_{2}^{2})
\leqslant
\mathfrak c_{\ref{risk_bound_approximate_numerator}}
\left(\|(bf)_h - bf\|_{2}^{2} +\frac{1}{Nh} +\frac{1}{n}\right) +
\frac{\mathfrak c_{\ref{risk_bound_approximate_numerator}}(\varepsilon)}{Nnh^{3 +\varepsilon}}.
\end{displaymath}
Finally, Propositions \ref{risk_bound_approximate_denominator} and \ref{risk_bound_approximate_numerator} allow to provide a risk bound on a truncated version the approximate Nadaraya-Watson estimator $\widehat b_{n,N,h}$ (see (\ref{NW_estimator_discrete})).
%


%
\begin{proposition}\label{risk_bound_approximate_NW}
Consider $\varepsilon > 0$, $m\in (0,1]$, and assume that $f(x) > m > 0$ for every $x\in [A,B]$ ($A,B\in\mathbb R$ such that $A < B$). Under the assumptions of Proposition \ref{risk_bound_approximate_numerator} and, in addition, Assumptions \ref{assumption_bound_derivative_density_t} and \ref{assumption_K_+}, there exist a constant $\mathfrak c_{\ref{risk_bound_approximate_NW}} > 0$, not depending on $\varepsilon$, $A$, $B$, $h$, $N$, $n$ and $t_0$, and a constant $\mathfrak c_{\ref{risk_bound_approximate_NW}}(\varepsilon) > 0$, depending on $\varepsilon$ but not on $A$, $B$, $h$, $N$, $n$ and $t_0$, such that
\begin{eqnarray*}
 \mathbb E(\|\widetilde b_{n,N,h} - b\|_{f,A,B}^{2})
 & \leqslant &
 \frac{\mathfrak c_{\ref{risk_bound_NW}}}{m^2
 \min\{1,t_{0}^{(1 -\varepsilon)/2},t_{0}^{1/2},t_{0}^{2q_2(\beta)},t_{0}^{2q_3},T - t_0\}}\\
 & &
 \hspace{1.5cm}\times\left[
\mathfrak c_{\ref{risk_bound_approximate_NW}}\left(
\|(bf)_h - bf\|_{2}^{2} + h^{2\beta} +\frac{1}{Nh} +\frac{1}{n}\right) +
\frac{\mathfrak c_{\ref{risk_bound_approximate_NW}}(\varepsilon)}{Nnh^{3 +\varepsilon}}\right]
\end{eqnarray*}
with $\widetilde b_{n,N,h}(.) :=\widehat b_{n,N,h}(.)\mathbf 1_{\widehat f_{n,N,h}(.) > m/2}$.
\end{proposition}
\noindent
The proof of Proposition \ref{risk_bound_approximate_NW} given Propositions \ref{risk_bound_approximate_denominator} and \ref{risk_bound_approximate_numerator} is almost the same than the proof of Proposition \ref{risk_bound_NW} given Propositions \ref{risk_bound_denominator} and \ref{risk_bound_numerator}. Of course one can establish a risk bound on the discrete-time approximate 2bNW estimator, but to focus on the 1 bandwidth estimator is clearer and sufficient to introduce the looCV selection method based on discrete-time observations of $X^1,\dots,X^N$ at Subsection \ref{section_looCV}. Now, assume that $bf$ satisfies Condition (\ref{Nikolskii_type_condition}) with $\gamma =\beta$, and that $K$ fulfills Assumption \ref{assumption_K_2} with $\upsilon =\beta$. Then, $\|(bf)_h - bf\|_{2}^{2}$ is of order $h^{2\beta}$. For the sake of simplicity, assume also that $b$ is bounded, and then let's take $\varepsilon = 0$ in Proposition \ref{risk_bound_approximate_NW}. First, note that the minimization problem
\begin{displaymath}
\min_{h\in (0,\infty)}\left\{
h^{2\beta} +\frac{1}{Nh} +\frac{1}{n} +\frac{1}{Nnh^3}
\right\}
\end{displaymath}
has unfortunately no explicit solutions. However, let us provide an upper-bound on the rate of our discrete-time estimator. Since $(nh^2)^{-1}\leqslant 1$, Proposition \ref{risk_bound_approximate_NW} says that the risk of $\widetilde b_{n,N,h}$ is at most of order $h^{2\beta} + 1/(Nh) + 1/n$. So, the optimal bandwidth for this bound is of order $N^{-1/(2\beta + 1)}$, leading to the rate
\begin{displaymath}
N^{-\frac{2\beta}{2\beta + 1}} +\frac{1}{n}.
\end{displaymath}
Moreover, by taking a bandwidth of order $N^{-1/(2\beta + 1)}$ such that $(nh^2)^{-1}\leqslant 1$, $N$ is at most of order $n^{(2\beta + 1)/2}$. So, clearly, the more $f$ and $bf$ are regular, the more $N$ can be chosen freely with respect to $n$, and if $N$ is of order $n^{(2\beta + 1)/2}$, then the risk of $\widetilde b_{n,N,h}$ is at most of order $1/n$. Finally, note that if $\beta = 1$, for a bandwidth of order $N^{-1/3}$ such that $(nh^2)^{-1}\leqslant 1$, then
\begin{displaymath}
1/n\leqslant h^2\propto N^{-2/3},
\end{displaymath}
and the rate of $\widetilde b_{n,N,h}$ is of order $N^{-2/3}$ (the optimal rate).
%


%
\section{Bandwidth selection and numerical experiments}\label{section_looCV_simulations}
This section deals with extensions of the PCO (see Lacour et al. \cite{LMR17}) and looCV methods to the Nadaraya-Watson estimator studied in this paper (see Subsections \ref{section_looCV} and \ref{section_PCO}). Subsection \ref{section_simulations} deals with some numerical experiments on the looCV based adaptive Nadaraya-Watson estimator which is, as explained in Comte and Marie \cite{CM21} in the nonparametric regression framework, numerically more satisfactory than the PCO based one. However, and this is its main advantage, the PCO based adaptive Nadaraya-Watson estimator offers theoretical guarantees: an oracle inequality is established in Subsection \ref{section_PCO}. Note also that the PCO method is easier to implement and numerically faster than the Goldenshluger-Lepski method which has been extended by Della Maestra and Hoffmann in \cite{DMH21} for their estimator of the drift function in McKean-Vlasov models.
%


%
\subsection{An extension of the Penalized Comparison to Overfitting method}\label{section_PCO}
Let $\mathcal H_N$ (resp. $\mathcal H_N'$) be a finite subset of $[h_0,1]$ (resp. $[h_0',1]$), where $h_0 > 0$ and $(Nh_{0}^{3})^{-1}\leqslant 1$ (resp. $h_0' > 0$ and $(Nh_0')^{-1}\leqslant 1$). Consider an additional kernel $\delta$,
\begin{equation}\label{penalty_proposal_bf}
\widehat h\in\arg\min_{h\in\mathcal H_N}
\{\|\widehat{bf}_{N,h} -\widehat{bf}_{N,h_0}\|_{2,\delta}^{2} +\textrm{pen}(h)\}
\end{equation}
with
\begin{equation}\label{penalty_PCO}
\textrm{pen}(h) :=
\frac{2}{(T - t_0)^{2}N^2}
\sum_{i = 1}^{N}\left\langle\int_{t_0}^{T}K_h(X_{s}^{i} -\cdot)dX_{s}^{i},
\int_{t_0}^{T}K_{h_0}(X_{s}^{i} -\cdot)dX_{s}^{i}\right\rangle_{2,\delta}
\textrm{ $;$ }
\forall h\in\mathcal H_N,
\end{equation}
and
\begin{equation}\label{penalty_proposal_f}
\widehat h'\in\arg\min_{h\in\mathcal H_N'}
\{\|\widehat f_{N,h} -\widehat f_{N,h_0'}\|_{2}^{2} +\textrm{pen}'(h)\}
\end{equation}
with
\begin{displaymath}
\textrm{pen}'(h) :=
\frac{2}{(T - t_0)^{2}N^2}
\sum_{i = 1}^{N}\left\langle\int_{t_0}^{T}K_h(X_{s}^{i} -\cdot)ds,
\int_{t_0}^{T}K_{h_0'}(X_{s}^{i} -\cdot)ds\right\rangle_2
\textrm{ $;$ }
\forall h\in\mathcal H_N'.
\end{displaymath}
This subsection deals with risk bounds on the adaptive estimators $\widehat{bf}_{N,\widehat h}(.)$ (see (\ref{NW_denominator})), $\widehat f_{N,\widehat h'}(.)$ (see (\ref{NW_numerator})) and
\begin{displaymath}
\widehat b_{N,\widehat h,\widehat h'}(x) =
\frac{\widehat{bf}_{N,\widehat h}(x)}{\widehat f_{N,\widehat h'}(x)}
\mathbf 1_{\widehat f_{N,\widehat h'}(x) > m/2}
\textrm{ $;$ }
x\in [A,B]
\end{displaymath}
with the notations of Proposition \ref{risk_bound_NW}. In the sequel, $K$, $\delta$ and $\sigma$ fulfill the following technical assumption.
%


%
\begin{assumption}\label{assumption_K_sigma_PCO}
The kernels $K$ and $\delta$ are continuously derivable on $\mathbb R$, the derivative of $K$ belongs to $\mathbb L^2(\mathbb R,dx)$, $\delta$ is positive and its derivative is bounded, and $\sigma$ is bounded.
\end{assumption}
\noindent
Moreover, recall that under Assumptions \ref{assumption_b_sigma} and \ref{assumption_bound_derivative_density_x}, $b^2$ and $\sigma^2$ belong to $\mathbb L^1(\mathbb R,f(x)dx)$ (see Remark \ref{integrability_functionals_X}).
%


%
\begin{theorem}\label{risk_bound_PCO_estimator_bf}
Under Assumptions \ref{assumption_b_sigma}, \ref{assumption_bound_derivative_density_x}, \ref{assumption_K_1} and \ref{assumption_K_sigma_PCO},
\begin{enumerate}
 \item There exist two deterministic constants $\mathfrak c_{\ref{risk_bound_PCO_estimator_bf},1},\mathfrak c_{\ref{risk_bound_PCO_estimator_bf},2} > 0$, not depending on $N$, such that for every $\vartheta\in (0,1)$ and $\lambda > 0$, with probability larger than $1 -\mathfrak c_{\ref{risk_bound_PCO_estimator_bf},1}|\mathcal H_N|e^{-\lambda}$,
 \begin{displaymath}
 \|\widehat{bf}_{N,\widehat h} - bf\|_{2,\delta}^{2}
 \leqslant
 (1 +\vartheta)\min_{h\in\mathcal H_N}
 \|\widehat{bf}_{N,h} - bf\|_{2,\delta}^{2} +
 \frac{\mathfrak c_{\ref{risk_bound_PCO_estimator_bf},2}}{\vartheta}
 \left[\|(bf)_{h_0} - bf\|_{2,\delta}^{2} +\frac{(1 +\lambda)^3}{N}\right].
 \end{displaymath}
 \item There exist two deterministic constants $\overline{\mathfrak c}_{\ref{risk_bound_PCO_estimator_bf},1},\overline{\mathfrak c}_{\ref{risk_bound_PCO_estimator_bf},2} > 0$, not depending on $N$, such that for every $\vartheta\in (0,1)$ and $\lambda > 0$, with probability larger than $1 -\overline{\mathfrak c}_{\ref{risk_bound_PCO_estimator_bf},1}|\mathcal H_N|e^{-\lambda}$,
 \begin{displaymath}
 \|\widehat f_{N,\widehat h'} - f\|_{2}^{2}
 \leqslant
 (1 +\vartheta)\min_{h'\in\mathcal H_N}
 \|\widehat f_{N,h'} - f\|_{2}^{2} +
 \frac{\overline{\mathfrak c}_{\ref{risk_bound_PCO_estimator_bf},2}}{\vartheta}
 \left[\|f_{h_0'} - f\|_{2}^{2} +\frac{(1 +\lambda)^3}{N}\right].
 \end{displaymath}
\end{enumerate}
\end{theorem}
%


%
\begin{corollary}\label{risk_bound_NW_PCO}
Under Assumptions \ref{assumption_b_sigma}, \ref{assumption_bound_derivative_density_x}, \ref{assumption_K_1} and \ref{assumption_K_sigma_PCO}, if $f(x),\delta(x) > m > 0$ for every $x\in [A,B]$ ($m\in (0,1]$ and $A,B\in\mathbb R$ such that $A < B$), then there exists a deterministic constant $\mathfrak c_{\ref{risk_bound_NW_PCO}} > 0$, not depending on $N$, $A$ and $B$, such that for every $\vartheta\in (0,1)$,
\begin{eqnarray*}
 \mathbb E(\|\widehat b_{N,\widehat h,\widehat h'} - b\|_{f,A,B}^{2})
 & \leqslant &
 \frac{2\mathfrak c_{\ref{risk_bound_NW}}(1\vee\|\delta\|_{\infty})}{m^3}\left[
 (1 +\vartheta)\min_{(h,h')\in\mathcal H_N\times\mathcal H_N'}
 \{\mathbb E(\|\widehat{bf}_{N,h} - bf\|_{2}^{2}) +
 \mathbb E(\|\widehat f_{N,h'} - f\|_{2}^{2})\}\right.\\
 & &
 \hspace{4cm}\left.
 +\frac{\mathfrak c_{\ref{risk_bound_NW_PCO}}}{\vartheta}\left(
 \|(bf)_{h_0} - bf\|_{2}^{2} +\|f_{h_0'} - f\|_{2}^{2} +\frac{1}{N}
 \right)\right].
\end{eqnarray*}
\end{corollary}
\noindent
Corollary \ref{risk_bound_NW_PCO} says that the risk of the adaptive estimator $\widehat b_{N,\widehat h,\widehat h'}$ is controlled by the sum of the minimal risks of
\begin{displaymath}
\widehat{bf}_{N,h}
\quad {\rm and}\quad
\widehat f_{N,h'}
\textrm{ $;$ }
(h,h')\in\mathcal H_N,
\end{displaymath}
up to a multiplicative constant and a negligible additive term.
%


%
\begin{remark}\label{uncommon_condition_bandwidths}
The condition $(Nh_{0}^{3})^{-1}\leqslant 1$ on the bandwidths collection $\mathcal H_N$ is quite uncomfortable but not that much because if $bf$ satisfies Condition (\ref{Nikolskii_type_condition}) with $\gamma =\beta\geqslant 2$, then the (unknown) bandwidth $h^*$ of order $N^{-1/(2\beta + 1)}$ such that our estimator of $bf$ reaches the bias-variance tradeoff (see Section \ref{section_continuous_time_NW}) possibly belongs to $\mathcal H_N$. Indeed, there exists an unknown constant $\mathfrak c^* > 0$ such that $h^* =\mathfrak c^*N^{-1/(2\beta + 1)}$, and then
\begin{displaymath}
\frac{1}{N(h^*)^3} =
(\mathfrak c^*)^{-3}N^{\frac{2}{2\beta + 1}(1 -\beta)}\leqslant 1
\end{displaymath}
for $N$ large enough. Moreover, the proof of Proposition \ref{risk_bound_PCO_estimator_bf} remains true by replacing the condition $(Nh_{0}^{3})^{-1}\leqslant 1$ by $(Nh_{0}^{3})^{-1}\leqslant\mathfrak m$ with $\mathfrak m > 0$. So, even for $\beta = 1$, $(N(h^*)^3)^{-1}\leqslant (\mathfrak c^*)^{-3}$ and then $h^*$ possibly belongs to $\mathcal H_N$ when $\mathfrak m$ is large enough.
\end{remark}
%


%
\begin{remark}\label{choice_additional_kernel}
A nice choice for $\delta$ is the standard normal density:
\begin{displaymath}
\delta(x) :=
\frac{e^{-x^2/2}}{\sqrt{2\pi}}
\textrm{ $;$ }
\forall x\in\mathbb R.
\end{displaymath}
First, $\delta$ obviously fulfills Assumption \ref{assumption_K_sigma_PCO}. Moreover, $\|\delta\|_{\infty}\leqslant 1$. Finally, by assuming that $f(x) > m_1$ for every $x\in [A,B]$ ($m_1\in (0,1]$ and $A,B\in\mathbb R$ such that $A < B$), since $\delta$ is continuous and positive on $\mathbb R$ (${\rm supp}(\delta) =\mathbb R$), necessarily there exists $m_2 > 0$ such that $\delta(x) > m_2$ for every $x\in [A,B]$. So, $f(x),\delta(x) > m = m_1\wedge m_2 > 0$ for every $x\in [A,B]$. Therefore, under Assumptions \ref{assumption_b_sigma}, \ref{assumption_bound_derivative_density_x}, \ref{assumption_K_1} and \ref{assumption_K_sigma_PCO}, by Corollary \ref{risk_bound_NW_PCO},
\begin{eqnarray*}
 \mathbb E(\|\widehat b_{N,\widehat h,\widehat h'} - b\|_{f,A,B}^{2})
 & \leqslant &
 \frac{2\mathfrak c_{\ref{risk_bound_NW}}}{m^3}\left[
 (1 +\vartheta)\min_{(h,h')\in\mathcal H_N\times\mathcal H_N'}
 \{\mathbb E(\|\widehat{bf}_{N,h} - bf\|_{2}^{2}) +
 \mathbb E(\|\widehat f_{N,h'} - f\|_{2}^{2})\}\right.\\
 & &
 \hspace{4cm}\left.
 +\frac{\mathfrak c_{\ref{risk_bound_NW_PCO}}}{\vartheta}\left(
 \|(bf)_{h_0} - bf\|_{2}^{2} +\|f_{h_0'} - f\|_{2}^{2} +\frac{1}{N}
 \right)\right]
\end{eqnarray*}
for every $\vartheta\in (0,1)$.
\end{remark}
%


%
\subsection{An extension of the leave-one-out cross-validation method}\label{section_looCV}
First of all, note that the estimator $\widehat b_{n,N,h}$ (see (\ref{NW_estimator_discrete})) can be written the following way:
\begin{displaymath}
\widehat b_{n,N,h}(x) =
\sum_{i = 1}^{N}\sum_{j = 0}^{n - 1}
\omega_{j}^{i}(x)(X_{t_{j + 1}}^{i} - X_{t_{j}^{i}})
\end{displaymath}
with
\begin{displaymath}
\omega_{j}^{i}(x) :=
\frac{K_h(X_{t_j}^{i} - x)}{\displaystyle{
\sum_{k = 1}^{N}\sum_{\ell = 0}^{n - 1}
K_h(X_{t_{\ell}}^{k} - x)(t_{\ell + 1} - t_{\ell})}}
\textrm{ $;$ }
\forall (j,i)\in\{0,\dots,n - 1\}\times\{1,\dots,N\},
\end{displaymath}
satisfying
\begin{displaymath}
\sum_{i = 1}^{N}\sum_{j = 0}^{n - 1}
\omega_{j}^{i}(x)(t_{j + 1} - t_j) = 1.
\end{displaymath}
This nice (weighted) representation of $\widehat b_{n,N,h}(x)$ allows us to consider the following extension of the well-known looCV criterion in our framework:
\begin{displaymath}
{\rm CV}(h) :=
\sum_{i = 1}^{N}\left[
\sum_{j = 0}^{n - 1}\widehat b_{n,N,h}^{-i}(X_{t_j}^{i})^2(t_{j + 1} - t_j)
- 2\sum_{j = 0}^{n - 1}\widehat b_{n,N,h}^{-i}(X_{t_j}^{i})(X_{t_{j + 1}}^{i} - X_{t_j}^{i})\right]
\end{displaymath}
with
\begin{displaymath}
\widehat b_{n,N,h}^{-i}(x) :=
\sum_{k\in\{1,\dots,N\}\backslash\{i\}}
\sum_{j = 0}^{n - 1}\omega_{j}^{k}(x)(X_{t_{j + 1}}^{k} - X_{t_j}^{k})
\textrm{ $;$ }
\forall i\in\{1,\dots,N\}.
\end{displaymath}
Let us explain heuristically this extension of the looCV criterion. By assuming that $dX_t = Y_tdt$, Equation (\ref{main_equation}) leads to the regression model
\begin{displaymath}
Y_{t_j} = b(X_{t_j}) +\varepsilon_{t_j}
\quad {\rm with}\quad
\int_{0}^{t_j}\varepsilon_sds =
\int_{0}^{t_j}\sigma(X_s)dW_s.
\end{displaymath}
Then, a natural extension of the looCV criterion is
\begin{eqnarray*}
 {\rm CV}^*(h) & := &
 \sum_{i = 1}^{N}\sum_{j = 0}^{n - 1}
 (Y_{t_j}^{i} -\widehat b_{n,N,h}^{-i}(X_{t_j}^{i}))^2(t_{j + 1} - t_j)\\
 & \approx &
 {\rm CV}(h) +\sum_{i = 1}^{N}\sum_{j = 0}^{n - 1}(Y_{t_j}^{i})^2(t_{j + 1} - t_j)
\end{eqnarray*}
because $Y_{t_j}(t_{j + 1} - t_j)\approx X_{t_{j + 1}} - X_{t_j}$ thanks to the assumption $dX_t = Y_tdt$. Of course ${\rm CV}^*(h)$ is not satisfactory because the last term of its previous decomposition doesn't exist, but since this term doesn't depend on $h$, to minimize ${\rm CV}^*(.)$ is almost equivalent to minimize ${\rm CV}(.)$ which only involves quantities existing without the condition $dX_t = Y_tdt$.
%


%
\subsection{Numerical experiments}\label{section_simulations}
\begin{figure}[!ht]
\centering
\includegraphics[scale=0.4]{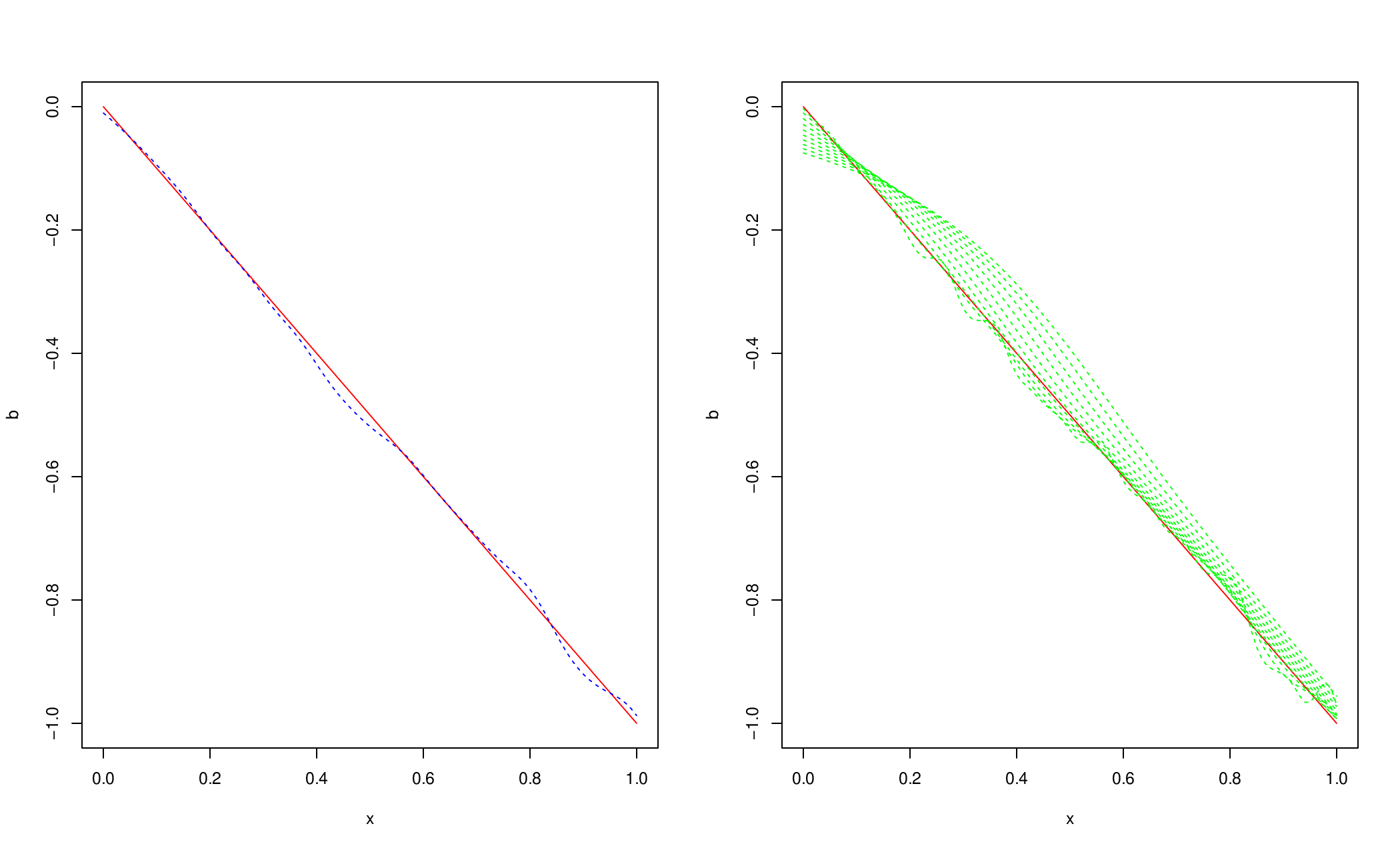} 
\caption{LooCV NW estimation for Model 1 (Langevin equation), $\widehat h = 0.04$.}
\label{graph_EDS_Lang}
\end{figure}

\begin{figure}[!ht]
\centering
\includegraphics[scale=0.4]{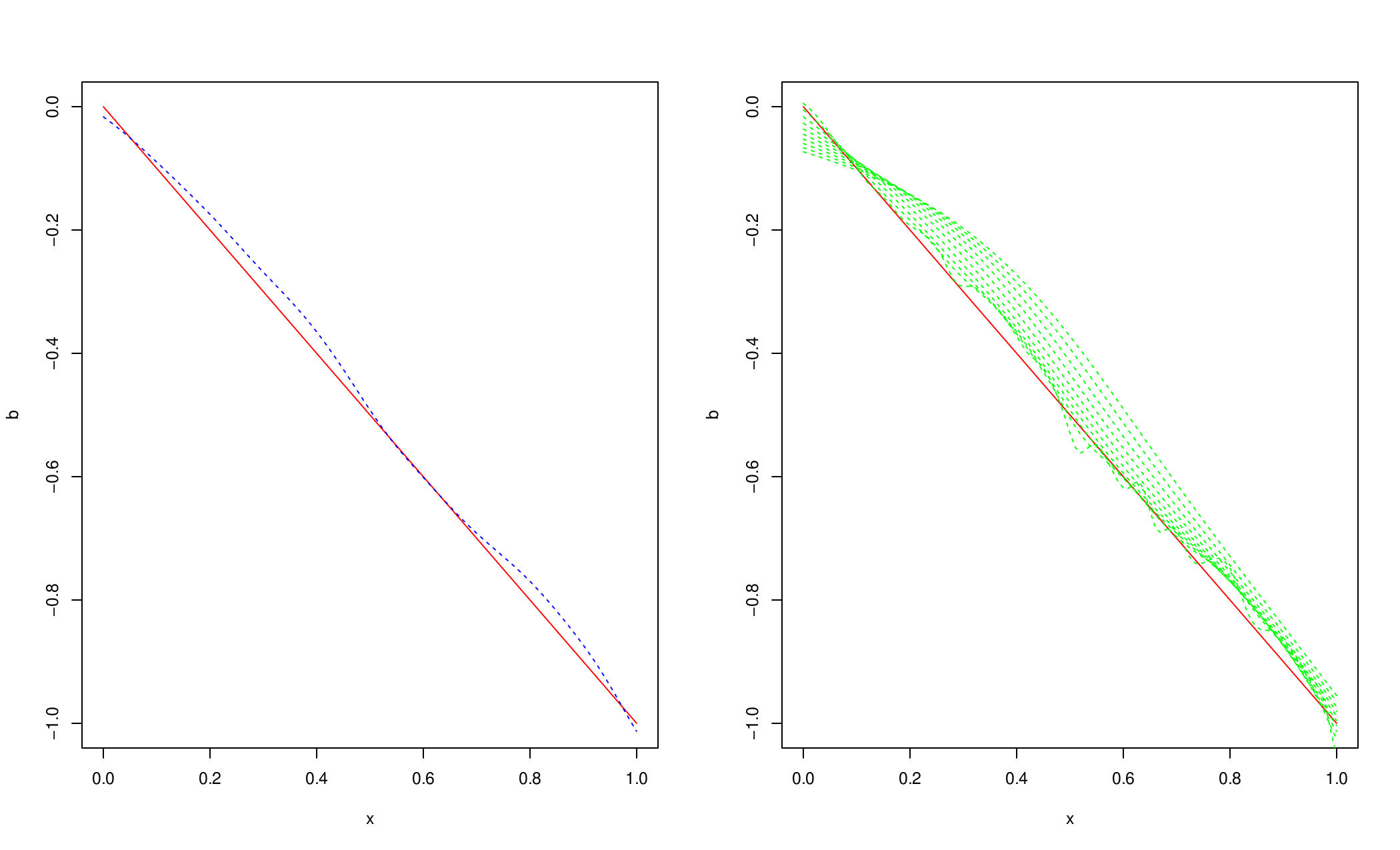} 
\caption{LooCV NW estimation for Model 2 (hyperbolic diffusion process), $\widehat h = 0.04$.}
\label{graph_EDS_diff_hyper}
\end{figure}

\begin{figure}[!h]
\centering
\includegraphics[scale=0.4]{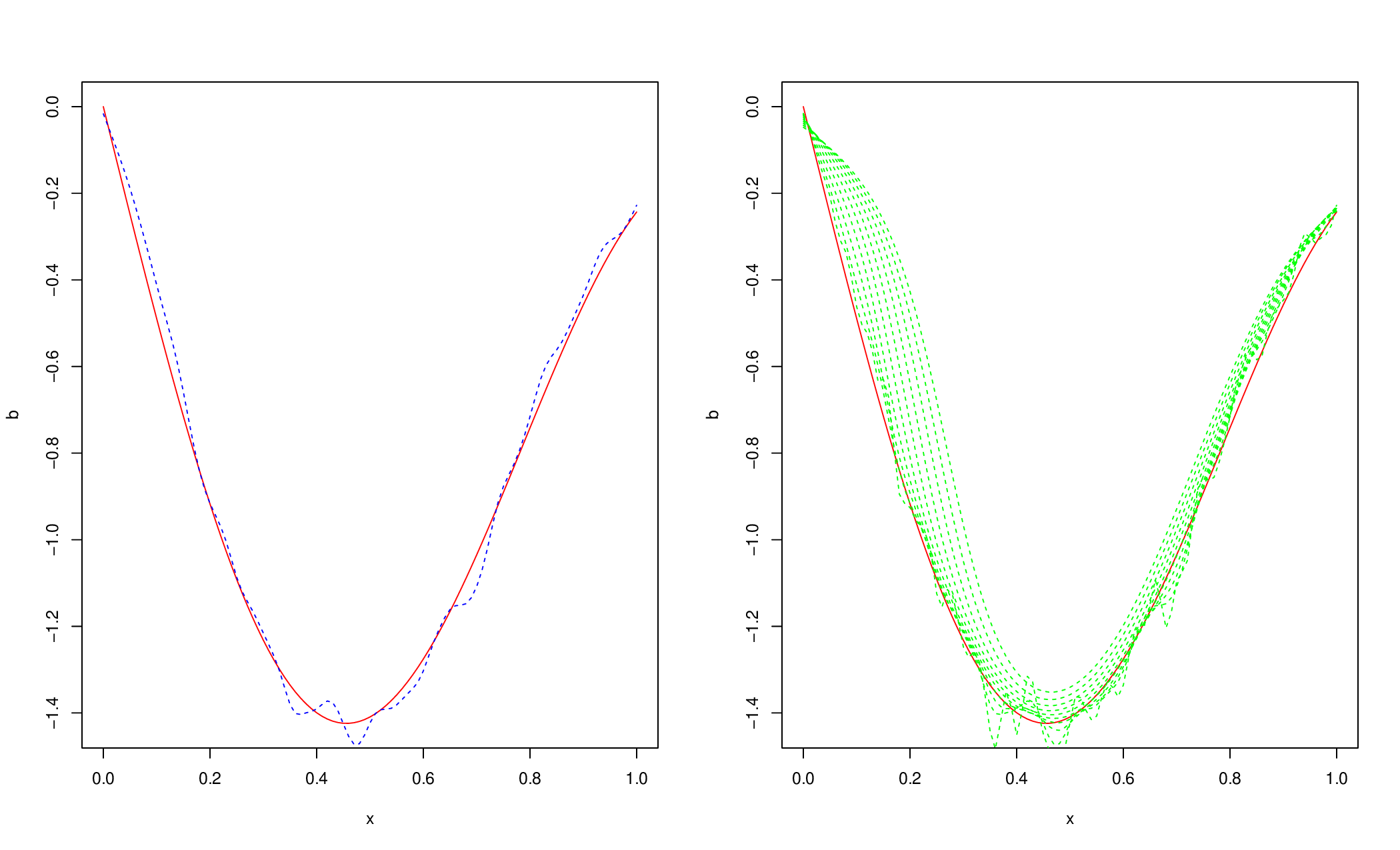} 
\caption{LooCV NW estimation for Model 3, $\widehat h = 0.02$.}
\label{graph_EDS_sin_add}
\end{figure}

\begin{figure}[!h]
\centering
\includegraphics[scale=0.4]{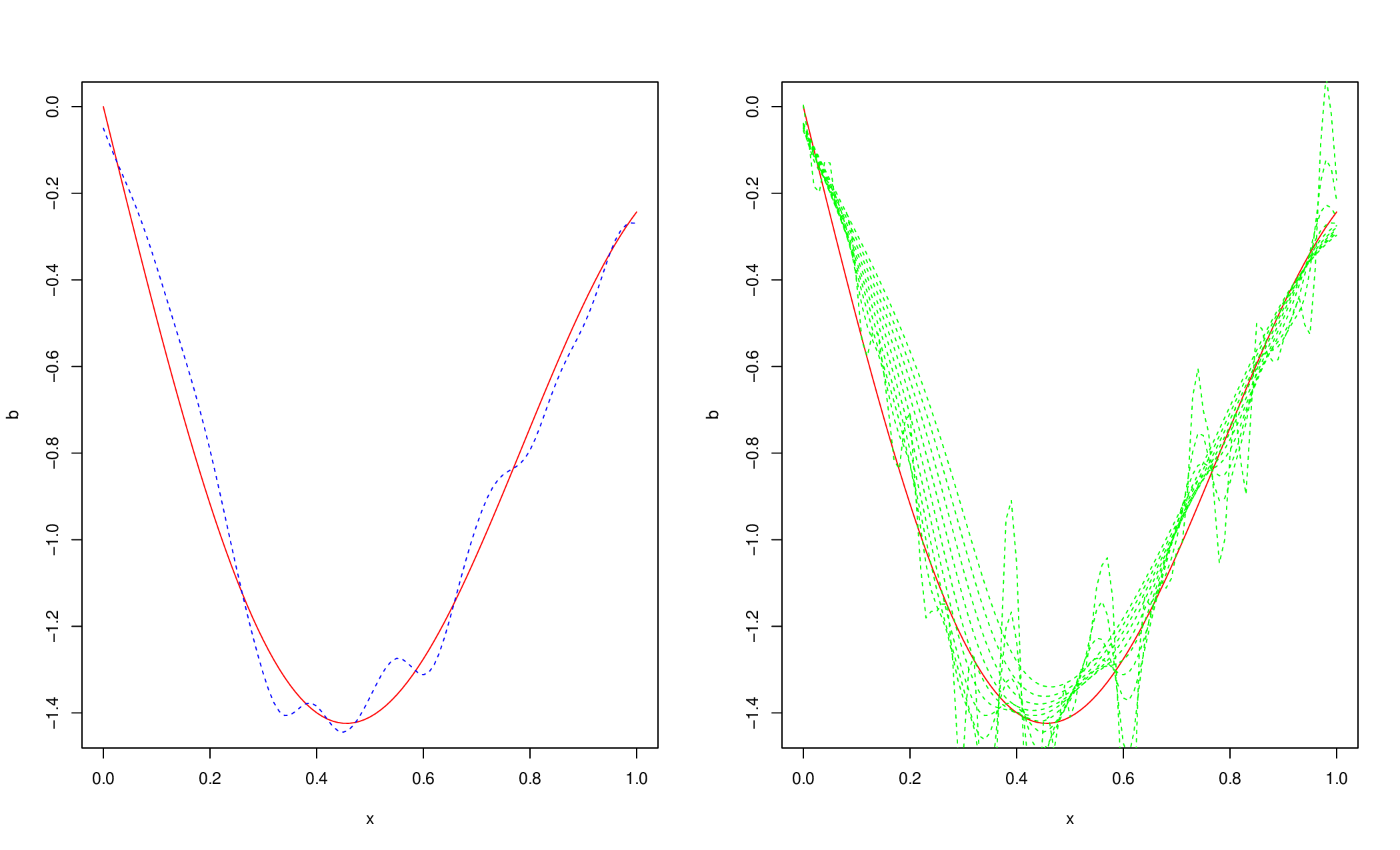} 
\caption{LooCV NW estimation for Model 4, $\widehat h = 0.06$.}
\label{graph_EDS_sin_multi}
\end{figure}

Some numerical experiments on our estimation method are presented in this subsection. The discrete-time approximate Nadaraya-Watson (NW) (see (\ref{NW_estimator_discrete})) estimator is computed on 4 datasets generated by SDEs with various types of vector fields. In each case, the bandwidth of the NW estimator is selected via the looCV method introduced at Subsection \ref{section_looCV}. On the one hand, two models with the same linear drift function are considered, but with an additive noise for the first one and a multiplicative noise for the second one:
\begin{itemize}
	\item[1.] The so-called Langevin equation, that is
		\begin{displaymath}
		X_t = x_0 -\int_{0}^{t}X_sds + 0.1\cdot W_t.
		\end{displaymath}
	\item[2.] The hyperbolic diffusion process, that is
		\begin{displaymath}
		X_t = x_0 -\int_{0}^{t}X_sds + 0.1\int_{0}^{t}\sqrt{1+X_{s}^{2}}dW_s.
		\end{displaymath}
\end{itemize}
On the other hand, two models having the same non-linear drift function involving $\sin(.)$ are considered, but here again with an additive noise for the first one and a multiplicative noise for the second one:
\begin{itemize}
	\item[3.] The third model is defined by
		\begin{displaymath}
		X_t = x_0 -\int_{0}^{t}(X_s + \sin(4X_s))ds + 0.1\cdot W_t.
		\end{displaymath}
	\item[4.] The fourth model is defined by
		\begin{displaymath}
		X_t = x_0 -\int_{0}^{t}(X_s + \sin(4X_s))ds +
		0.1\int_{0}^{t}(2 +\cos(X_s))dW_s.
		\end{displaymath}
\end{itemize}
The models and the estimator are implemented by taking $N = 200$, $n = 50$, $T = 5$, $x_0 = 2$, $t_0 = 1$ and $K$ the Gaussian kernel $z\mapsto (2\pi)^{-1/2}e^{-z^2/2}$. For Models 1 and 2, the estimator of the drift function is computed for the bandwidths set
\begin{displaymath}
\mathcal H_1 :=\{0.02k\textrm{ $;$ }k = 1,\dots,10\},
\end{displaymath}
and for Models 3 and 4, it is computed for the bandwidths set
\begin{displaymath}
\mathcal H_2 :=\{0.01k\textrm{ $;$ }k = 1,\dots,10\}.
\end{displaymath}
Each set of bandwidths has been chosen after testing different values of $h$, to see with which ones the estimation performs better. To choose smaller values in the second set of bandwidths allows to check that the looCV method does not systematically select the smallest bandwidth for Models 3 and 4.
\\
\\
For each of the previous models, on Figures \ref{graph_EDS_Lang}, \ref{graph_EDS_diff_hyper}, \ref{graph_EDS_sin_add} and \ref{graph_EDS_sin_multi} respectively, the true drift function (in red) and the looCV adaptive NW estimator (in blue) are plotted on the left-hand side, and the beam of proposals is plotted in green on the right-hand side. On Figures \ref{graph_EDS_Lang} and \ref{graph_EDS_diff_hyper}, one can see that the drift function is well estimated by the looCV adaptive NW estimator, with a MSE equal to $2.95\cdot 10^{-4}$ for the Langevin equation and to $8.31\cdot 10^{-4}$ for the hyperbolic diffusion process. As presumed, the multiplicative noise in Model 2 slightly degrades the MSE. Note that when the bandwidth is too small, the estimation degrades, but the looCV method selects a higher value of $h$ which performs better on the estimation. This means that, as expected, the looCV method selects a reasonable approximation of the bandwidth for which the NW estimator reaches the bias-variance tradeoff. On Figures \ref{graph_EDS_sin_add} and \ref{graph_EDS_sin_multi}, one can see that the drift function of Models 3 and 4 is still well estimated by our looCV adaptive NW estimator. However, note that there is a significant degradation of the MSE, which is equal to $2.89\cdot 10^{-3}$ for Model 3 and to $9.26\cdot 10^{-3}$ for Model 4. This is probably related to the {\it nonlinearity} of the drift function to estimate. Once again, to consider a multiplicative noise in Model 4 degrades the estimation quality with respect to Model 3. As for Models 1 and 2, note that the looCV does not systematically select the smallest bandwidth.
\newline
\newline
For Model 1, at levels $n = 10,20,\dots,100$, Figure \ref{MSE_f(N,n)_Lang} shows the evolution of the MSE of the looCV adaptive NW estimator as a function of $N$. For this study, the value of $N$ ranges from $20$ to $200$. Figure \ref{MSE_f(N,n)_Lang} shows that the MSE of our adaptive estimator remains low regardless to the value of $(n,N)$ (from $4.50\cdot 10^{-5}$ to $9.01\cdot 10^{-3}$), decreases when $N$ increases (for each $n$), and decreases when $n$ increases (for a fixed $N$). This is consistent with the risk bounds of Section 4. Note also that for $N\geqslant 70$, there is no significant gain to take $n$ larger than $30$. For Model 3, Figure \ref{MSE_f(N,n)_model_3} shows the evolution of the MSE of the looCV adaptive NW estimator as a function of $N$ and leads to the same conclusions than for Model 1. Note anyway that due to the {\it nonlinearity} of $b$, the MSE of our adaptive estimator reaches higher values (from $2.67\cdot 10^{-4}$ to $4.49\cdot 10^{-2}$) than for Model 1. Again, there is no significant gain to take $n$ larger than $30$, and above all larger than $70$.
\begin{figure}[!h]
\centering
\includegraphics[scale=0.55]{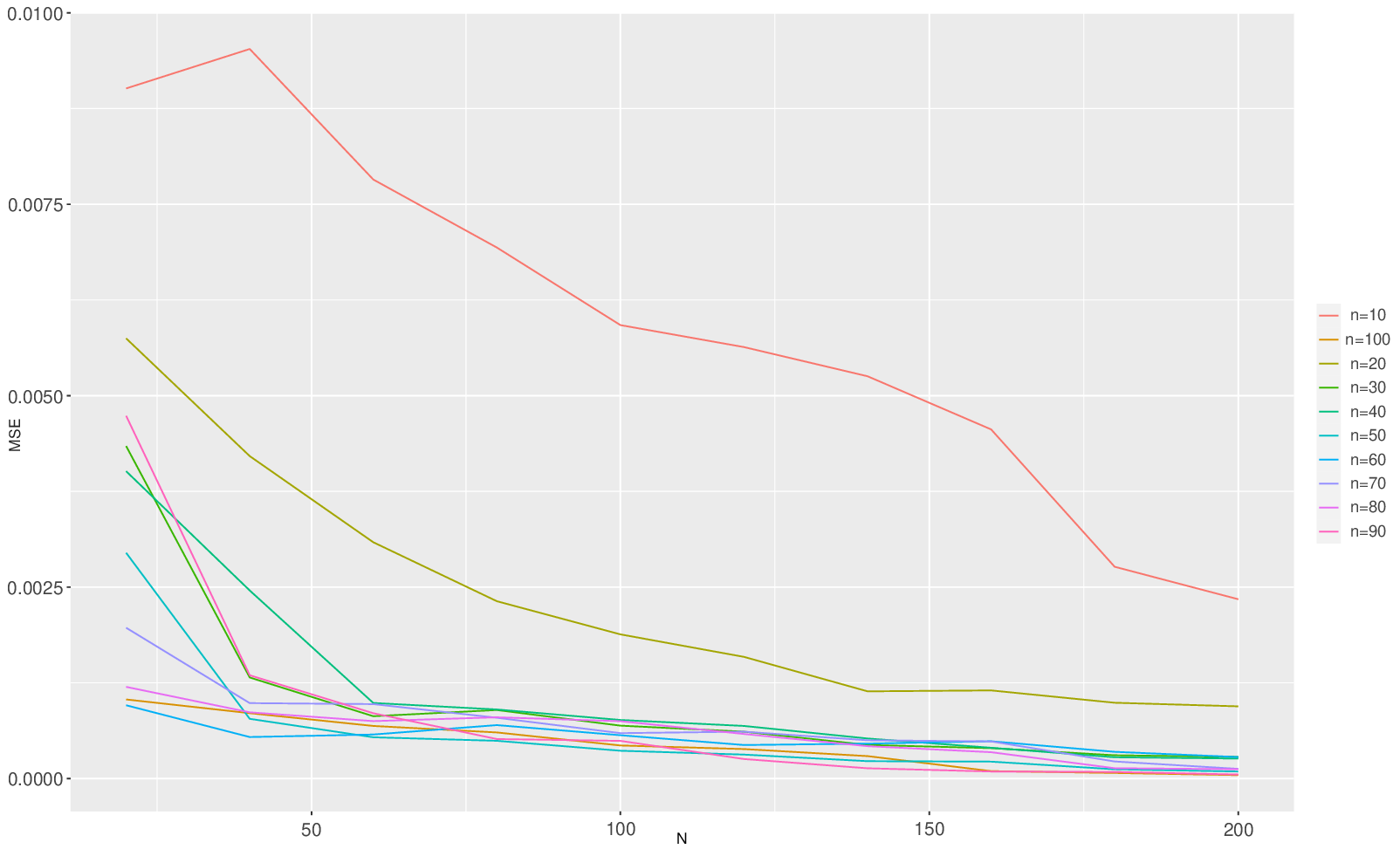} 
\caption{MSE of the looCV estimator with respect to $N$ and $n$ for Model 1.}
\label{MSE_f(N,n)_Lang}
\end{figure}
\begin{figure}[!h]
\centering
\includegraphics[scale=0.55]{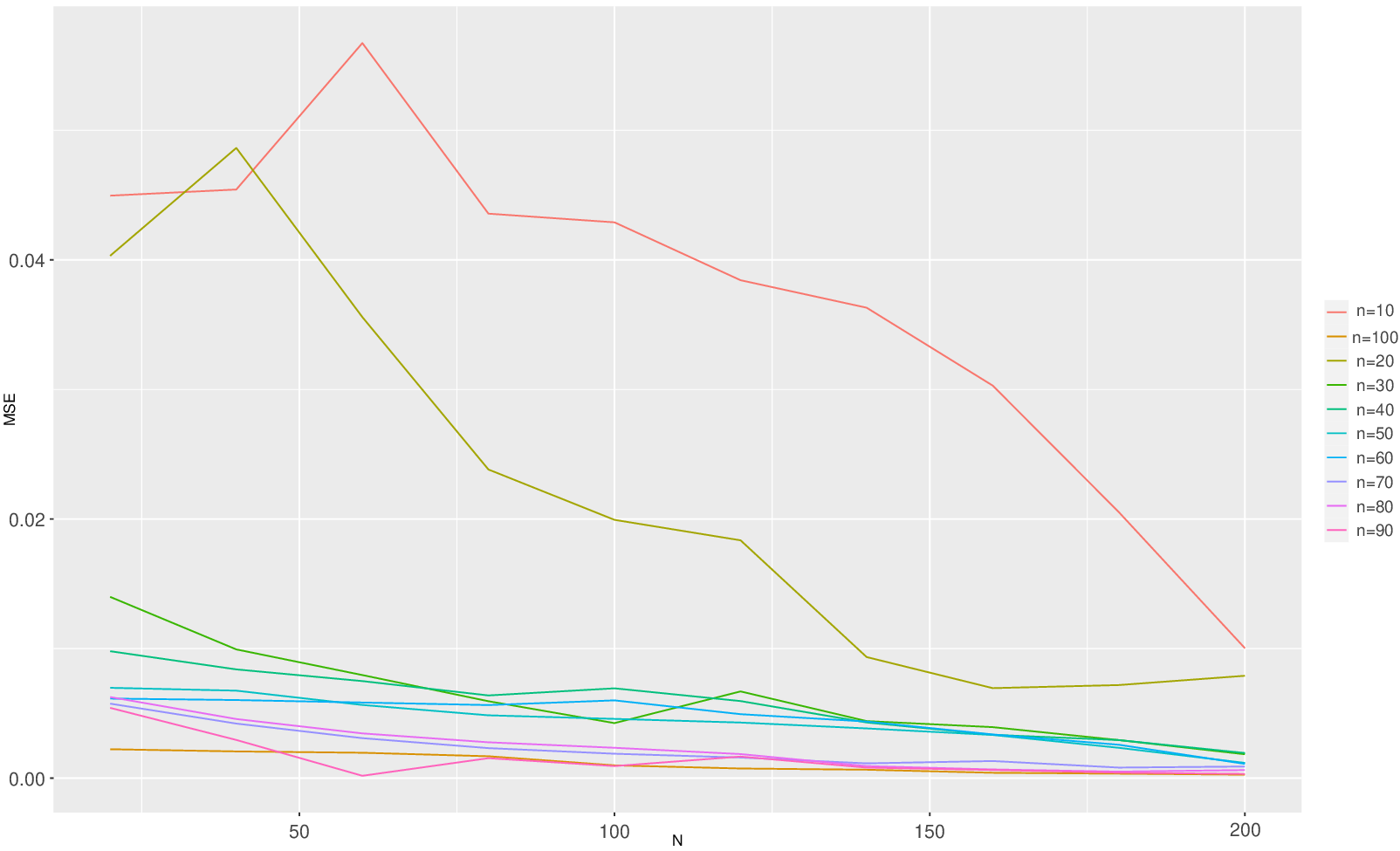} 
\caption{MSE of the looCV estimator with respect to $N$ and $n$ for Model 3.}
\label{MSE_f(N,n)_model_3}
\end{figure}
\newline
\newline
Finally, for each model, Table \ref{table_MSE} gathers the mean MSE of 100 looCV NW estimations of the drift function as well as the mean MSE of the corresponding 100 oracle estimations. The mean MSEs are globally low, but significantly higher for the models with a nonlinear drift function (Models 3 and 4) than for the models with a linear one (Models 1 and 2). Moreover, for each drift function, the mean MSE is slightly degraded for the models with a multiplicative noise (Models 2 and 4) with respect to the models with an additive one (Models 1 and 3). Note also that for each model, the mean MSE of the looCV estimations is close to the mean MSE of the corresponding oracle estimations. This means that our looCV method performs well in practice.
\begin{table}[h!]
\begin{center}
\begin{tabular}{|c|c|c|}
\hline
& looCV & Oracle \\
	\hline
	\hline
Model 1 & $3.03\cdot 10^{-4}$ & $2.67\cdot 10^{-4}$\\
	\hline
Model 2 & $6.52\cdot 10^{-4}$ & $4.96\cdot 10^{-4}$\\
	\hline
Model 3 & $2.45\cdot 10^{-3}$ & $1.99\cdot 10^{-3}$\\
	\hline
Model 4 & $9.15\cdot 10^{-3}$ & $6.02\cdot 10^{-3}$\\
	\hline
\end{tabular}
\medskip
\caption{Mean MSEs of 100 looCV adaptive NW estimations compared to the oracle estimations.}\label{table_MSE}
\end{center}
\end{table}
\begin{remark}\label{t_0_practice}
Note that to take $t_0\geqslant 1$ (here $t_0 = 1$) is recommended even in numerical experiments. Indeed, for instance, the mean MSE of 10 looCV estimations for Model 1 is significantly lower with $t_0 = 1$ ($2.49\cdot 10^{-4}$) than with $t_0 = 0$ ($3.82\cdot 10^{-3}$).
\end{remark}
%


%
\section{Concluding remarks}\label{section_conclusion}
In this paper, first, a risk bound on our continuous-time Nadaraya-Watson estimator of $b$ has been established. This bound is satisfactory because it leads to a rate of same order than in the classic nonparametric regression framework (see Comte \cite{COMTE17}, Chapter 4), and of same order than in Della Maestra and Hoffmann \cite{DMH21} for their estimator of the drift function in McKean-Vlasov models. Then, a risk bound on a discrete-time approximate estimator of $b$ has been established too. The bound is satisfactory when $b$ and $\sigma$ are bounded, but a bit degraded when $b$ is unbounded. To improve this bound will be the subject of future investigations.
\\
\\
In a second part, two bandwidth selection methods are provided. The first one is an extension of the PCO method to the 2bNW estimator of $b$ in the spirit of Comte and Marie \cite{CM21}. An oracle inequality is established but under the condition $(Nh^3)^{-1}\leqslant 1$ (instead of $(Nh)^{-1}\leqslant 1$) on the bandwidths collection. Unfortunately, it seems difficult to bypass this condition because of some constants involved in Bernstein's inequality and in the concentration inequality for U-statistics of Gin\'e and Nickl \cite{GN15} (see Subsection \ref{section_proof_PCO}), but as explained at Remark \ref{uncommon_condition_bandwidths} this condition is not so bad. The second bandwidth selection method is an extension of the looCV procedure for the discrete-time approximate estimator written has a convex combination. As in the nonparametric regression framework, this method is numerically satisfactory but it seems difficult to establish a theoretical risk bound on the associated adaptive estimator.
\\
\\
Finally, the estimation of $b$ has been only investigated in the case of one-dimensional diffusion processes because of its simplicity, but by following the same ideas than Halconruy and Marie used in the nonparametric regression framework in \cite{HM20}, the major part of the results of the present paper should be extendable to multidimensional diffusion processes.
%


%

%
\appendix
%


%
\section{Proofs}\label{section_proofs}
%


%
\subsection{Proof of Corollary \ref{Nikolskii_f}}
First of all, since $p_t(x_0,x) > 0$ for every $(t,x)\in (0,T]\times\mathbb R$,
\begin{displaymath}
f(x) =
\frac{1}{T - t_0}\int_{t_0}^{T}p_t(x_0,x)dt > 0
\end{displaymath}
for every $x\in\mathbb R$. Consider $\ell\in\{0,\dots,\beta - 1\}$ and $\theta\in\mathbb R_+$. Thanks to the bound on $(t,x)\mapsto\partial_{x}^{\ell + 1}p_t(x_0,x)$ given in Assumption \ref{assumption_bound_derivative_density_x},
\begin{eqnarray*}
 & &\|f^{(\ell)}(\cdot +\theta) - f^{(\ell)}\|_{2}^{2} =
 \int_{-\infty}^{\infty}[f^{(\ell)}(x + x_0 +\theta) - f^{(\ell)}(x + x_0)]^2dx\\
 & &
 \hspace{1cm}\leqslant
 \frac{1}{T - t_0}\int_{t_0}^{T}\int_{-\infty}^{\infty}
 (\partial_{2}^{\ell}p_t(x_0,x + x_0 +\theta) -\partial_{2}^{\ell}p_t(x_0,x + x_0))^2dxdt\\
 & &
 \hspace{1cm}\leqslant
 \frac{\theta^2}{T - t_0}\int_{t_0}^{T}\int_{-\infty}^{\infty}
 \sup_{z\in [x,x +\theta]}|\partial_{2}^{\ell + 1}p_t(x_0,z + x_0)|^2dxdt\\
 & &
 \hspace{1cm}\leqslant
 \mathfrak c_{\ref{assumption_bound_derivative_density_x},2}(\ell + 1)^2
 \frac{\theta^2}{T - t_0}\int_{t_0}^{T}
 \frac{1}{t^{2q_2(\ell + 1)}}\int_{-\infty}^{\infty}\sup_{z\in [x,x +\theta]}
 \exp\left(-2\mathfrak m_{\ref{assumption_bound_derivative_density_x},2}(\ell + 1)\frac{z^2}{t}\right)dxdt\\
 & &
 \hspace{1cm} =
 \mathfrak c_{\ref{assumption_bound_derivative_density_x},2}(\ell + 1)^2
 \frac{\theta^2}{T - t_0}\int_{t_0}^{T}
 \frac{1}{t^{2q_2(\ell + 1)}}\times\\
 & &
 \hspace{2cm}
 \left[
 \int_{-\infty}^{-\theta}\exp\left(-2\mathfrak m_{\ref{assumption_bound_derivative_density_x},2}(\ell + 1)\frac{(x +\theta)^2}{t}\right)dx +\theta +
 \int_{0}^{\infty}\exp\left(-2\mathfrak m_{\ref{assumption_bound_derivative_density_x},2}(\ell + 1)\frac{x^2}{t}\right)dx\right]dt\\
 & &
 \hspace{1cm}\leqslant
 \frac{1}{t_{0}^{2q_2(\ell + 1)}}\left[
 \mathfrak c_1\theta^2 +
 \theta^3\max_{k\in\{0,\dots,\beta - 1\}}
 \mathfrak c_{\ref{assumption_bound_derivative_density_x},2}(k + 1)^2\right]
\end{eqnarray*}
with
\begin{displaymath}
\mathfrak c_1 =
2\max_{k\in\{0,\dots,\beta - 1\}}\left\{\mathfrak c_{\ref{assumption_bound_derivative_density_x},2}(k + 1)^2
\int_{0}^{\infty}\exp\left(-2\mathfrak m_{\ref{assumption_bound_derivative_density_x},2}(k + 1)\frac{x^2}{T}\right)dx\right\},
\end{displaymath}
and the same way,
\begin{eqnarray*}
 & & \|f^{(\ell)}(\cdot -\theta) - f^{(\ell)}\|_{2}^{2}
 \leqslant
 \frac{\theta^2}{T - t_0}\int_{t_0}^{T}\int_{-\infty}^{\infty}
 \sup_{z\in [x -\theta,x]}|\partial_{2}^{\ell + 1}p_t(x_0,z + x_0)|^2dxdt\\
 & &
 \hspace{1cm}\leqslant
 \mathfrak c_{\ref{assumption_bound_derivative_density_x},2}(\ell + 1)^2
 \frac{\theta^2}{T - t_0}\int_{t_0}^{T}
 \frac{1}{t^{2q_2(\ell + 1)}}\times\\
 & &
 \hspace{2cm}
 \left[
 \int_{-\infty}^{0}\exp\left(-2\mathfrak m_{\ref{assumption_bound_derivative_density_x},2}(\ell + 1)\frac{x^2}{t}\right)dx +
 \theta +
 \int_{\theta}^{\infty}\exp\left(-2\mathfrak m_{\ref{assumption_bound_derivative_density_x},2}(\ell + 1)\frac{(x -\theta)^2}{t}\right)dx\right]dt\\
 & &
 \hspace{1cm}\leqslant
 \frac{1}{t_{0}^{2q_2(\ell + 1)}}\left[
 \mathfrak c_1\theta^2 +
 \theta^3\max_{k\in\{0,\dots,\beta - 1\}}
 \mathfrak c_{\ref{assumption_bound_derivative_density_x},2}(k + 1)^2\right].
\end{eqnarray*}
This concludes the proof.
%


%
\subsection{Proof of Proposition \ref{risk_bound_denominator}}
First of all, the bias of $\widehat f_{N,h}(x)$ is denoted by $\mathfrak b(x)$ and its variance by $\mathfrak v(x)$. Moreover, let us recall the bias-variance decomposition of the $\mathbb L^2$-risk of $\widehat f_{N,h}$:
\begin{eqnarray*}
\mathbb E(\|\widehat f_{N,h} - f\|_{2}^{2}) =
\int_{-\infty}^{\infty}\mathfrak b(x)^2dx +
\int_{-\infty}^{\infty}\mathfrak v(x)dx.
\end{eqnarray*}
On the one hand, let us find a suitable bound on the integrated variance of $\widehat f_{N,h}$. Since $X^1,\dots,X^N$ are i.i.d. copies of $X$, and thanks to Jensen's inequality,
\begin{eqnarray*}
 \mathfrak v(x) & = &
 \mathrm{var}\left(\frac{1}{N(T - t_0)}\sum_{i = 1}^{N}\int_{t_0}^{T}K_h(X_{t}^{i} - x)dt\right)\\
 & = &
 \frac{1}{N(T - t_0)^2}\mathrm{var}\left(\int_{t_0}^{T}K_h(X_t - x)dt\right)
 \leqslant
 \frac{1}{N}\mathbb E\left[\left(\int_{t_0}^{T}K_h(X_t - x)\frac{dt}{T - t_0}\right)^2\right]\\
 & \leqslant &
 \frac{1}{N(T - t_0)}\int_{t_0}^{T}\mathbb E(K_h(X_t - x)^2)dt
 =\frac{1}{N}\int_{-\infty}^{\infty}K_h(z - x)^2f(z)dz.
\end{eqnarray*}
Thus, since $K$ is symmetric,
\begin{eqnarray*}
 \int_{-\infty}^{\infty}\mathfrak v(x)dx
 & \leqslant &
 \frac{1}{N}
 \int_{-\infty}^{\infty}f(z)\int_{-\infty}^{\infty}K_h(z - x)^2dxdz\\
 & = &
 \frac{1}{Nh}\left(\int_{-\infty}^{\infty}f(z)dz\right)\left(\int_{-\infty}^{\infty}K(x)^2dx\right)
 =\frac{\|K\|_{2}^{2}}{Nh}.
\end{eqnarray*}
On the other hand, let us find a suitable bound on the integrated squared-bias of $\widehat f_{N,h}(x)$. Since $X^1,\dots,X^N$ are i.i.d. copies of $X$,
\begin{eqnarray*}
 \mathfrak b(x) & = &
 \frac{1}{T - t_0}\int_{t_0}^{T}\mathbb E(K_h(X_t - x))dt - f(x)\\
 & = &
 \frac{1}{h}
 \int_{-\infty}^{\infty}
 K\left(\frac{z - x}{h}\right)f(z)dz - f(x)\\
 & = &
 \int_{-\infty}^{\infty}
 K(z)(f(hz + x) - f(x))dz.
\end{eqnarray*}
First, assume that $\beta = 1$. By Assumption \ref{assumption_K_2}, the generalized Minkowski inequality and Corollary \ref{Nikolskii_f},
\begin{eqnarray*}
 \int_{-\infty}^{\infty}\mathfrak b(x)^2dx
 & \leqslant &
 \int_{-\infty}^{\infty}\left(
 \int_{-\infty}^{\infty}K(z)(f(hz + x) - f(x))dz\right)^2dx\\
 & \leqslant &
 \left[\int_{-\infty}^{\infty}K(z)\left(
 \int_{-\infty}^{\infty}(f(hz + x) - f(x))^2dx\right)^{1/2}dz\right]^2
 \leqslant\mathfrak c_1(t_0)h^2
\end{eqnarray*}
with
\begin{displaymath}
\mathfrak c_1(t_0) =
\frac{\mathfrak c_{\ref{Nikolskii_f}}}{t_{0}^{2q_2(1)}}
\left(\int_{-\infty}^{\infty}|z|(1 + |z|^{1/2})|K(z)|dz\right)^2.
\end{displaymath}
Now, assume that $\beta\geqslant 2$. By the Taylor formula with integral remainder, for every $z\in\mathbb R$,
\begin{displaymath}
f(hz + x) - f(x) =
\mathbf 1_{\beta\geqslant 3}
\sum_{\ell = 1}^{\beta - 2}\frac{(hz)^{\ell}}{\ell !}f^{(\ell)}(x) +
\frac{(hz)^{\beta - 1}}{(\beta - 2)!}
\int_{0}^{1}(1 -\tau)^{\beta - 2}f^{(\beta - 1)}(\tau hz + x)d\tau.
\end{displaymath}
Then, by Assumption \ref{assumption_K_2}, the generalized Minkowski inequality (two times) and Corollary \ref{Nikolskii_f},
\begin{eqnarray*}
 & &
 \int_{-\infty}^{\infty}\mathfrak b(x)^2dx
 =\int_{-\infty}^{\infty}\left(
 \int_{-\infty}^{\infty}K(z)(f(hz + x) - f(x))dz\right)^2dx\\
 & &
 \hspace{1cm} =
 \frac{h^{2(\beta - 1)}}{|(\beta - 2)!|^2}
 \int_{-\infty}^{\infty}\left(
 \int_{-\infty}^{\infty}z^{\beta - 1}K(z)
 \int_{0}^{1}(1 -\tau)^{\beta - 2}
 [f^{(\beta - 1)}(\tau hz + x) - f^{(\beta - 1)}(x)]d\tau dz\right)^2dx\\
 & &
 \hspace{1cm}\leqslant
 \frac{h^{2(\beta - 1)}}{|(\beta - 2)!|^2}\times\\
 & &
 \hspace{2cm}
 \left[\int_{-\infty}^{\infty}
 |z|^{\beta - 1}|K(z)|\int_{0}^{1}(1 -\tau)^{\beta - 2}
 \left(\int_{-\infty}^{\infty}
 [f^{(\beta - 1)}(\tau hz + x) - f^{(\beta - 1)}(x)]^2dx\right)^{1/2}d\tau dz\right]^2\\
 & &
 \hspace{1cm}\leqslant
 \frac{\mathfrak c_{\ref{Nikolskii_f}}h^{2\beta}}{|(\beta - 2)!|^2t_{0}^{2q_2(\beta)}}
 \left(\int_{-\infty}^{\infty}
 |z|^{\beta - 1}|K(z)|\int_{0}^{1}(1 -\tau)^{\beta - 2}
 [\tau |z| + (\tau |z|)^{3/2}]
 d\tau dz\right)^2
 \leqslant\mathfrak c_2(t_0)h^{2\beta}
\end{eqnarray*}
with
\begin{displaymath}
\mathfrak c_2(t_0) =
\frac{\mathfrak c_{\ref{Nikolskii_f}}}{|(\beta - 2)!|^2t_{0}^{2q_2(\beta)}}
\left(\int_{-\infty}^{\infty}|z|^{\beta}(1 + |z|^{1/2})|K(z)|dz\right)^2.
\end{displaymath}
This concludes the proof.
%


%
\subsection{Proof of Proposition \ref{risk_bound_numerator}}
First of all,
\begin{eqnarray*}
\mathbb E(\|\widehat{bf}_{N,h} - bf\|_{2}^{2}) = \int_{-\infty}^{\infty}\mathfrak b(x)^2dx
+\int_{-\infty}^{\infty}\mathfrak v(x)dx
\end{eqnarray*}
where $\mathfrak b(x)$ (resp. $\mathfrak v(x)$) is the bias (resp. the variance) term of $\widehat{bf}_{N,h}(x)$ for any $x\in\mathbb R$. On the one hand, let us find a suitable bound on the integrated variance of $\widehat{bf}_{N,h}$. Since $X^1,\dots,X^N$ are i.i.d. copies of $X$,
\begin{eqnarray*}
 \mathfrak v(x)
 & = &
 \textrm{var}\left(\frac{1}{N(T - t_0)}
 \sum_{i = 1}^{N}\int_{t_0}^{T}K_h(X_{t}^{i} - x)dX_{t}^{i}\right)
 \leqslant
 \frac{1}{N(T - t_0)^2}
 \mathbb E\left[\left(\int_{t_0}^{T}K_h(X_t - x)dX_t\right)^2\right]\\
 & \leqslant &
 \frac{2}{N}\mathbb E\left[
 \left(\int_{t_0}^{T}K_h(X_t - x)b(X_t)\frac{dt}{T - t_0}\right)^2 +
 \frac{1}{(T - t_0)^2}\left(\int_{t_0}^{T}K_h(X_t - x)\sigma(X_t)dW_t\right)^2\right].
\end{eqnarray*}
In the right-hand side of the previous inequality, Jensen's inequality on the first term and the isometry property for It\^o's integral on the second one give
\begin{eqnarray*}
 \mathfrak v(x) & \leqslant &
 \frac{2}{N(T - t_0)}\int_{t_0}^{T}\mathbb E[K_h(X_t - x)^2b(X_t)^2]dt +
 \frac{2}{N(T - t_0)^2}\int_{t_0}^{T}\mathbb E[K_h(X_t - x)^2\sigma(X_t)^2]dt\\
 & = &
 \frac{2}{N}\int_{-\infty}^{\infty}K_h(z - x)^2b(z)^2f(z)dz +
 \frac{2}{N(T - t_0)}\int_{-\infty}^{\infty}K_h(z - x)^2\sigma(z)^2f(z)dz.
\end{eqnarray*}
Moreover, $K$ is symmetric and $K\in\mathbb L^2(\mathbb R,dx)$ by Assumption \ref{assumption_K_1}, and $b,\sigma\in\mathbb L^2(\mathbb R,f(x)dx)$ by Remark \ref{integrability_functionals_X}. Then,
\begin{eqnarray*}
 \int_{-\infty}^{\infty}\mathfrak v(x)dx
 & \leqslant &
 \frac{2}{N}\int_{\mathbb R^2}K_h(z - x)^2b(z)^2f(z)dzdx +
 \frac{2}{N(T - t_0)}\int_{\mathbb R^2}K_h(z - x)^2\sigma(z)^2f(z)dzdx\\
 & = &
 \frac{2}{Nh}\int_{-\infty}^{\infty}b(z)^2f(z)
 \int_{-\infty}^{\infty}K(x)^2dxdz +
 \frac{2}{N(T - t_0)h}\int_{-\infty}^{\infty}\sigma(z)^2f(z)
 \int_{-\infty}^{\infty}K(x)^2dxdz\\
 & \leqslant &
 \frac{2\|K\|_{2}^{2}}{Nh}\left(\int_{-\infty}^{\infty}b(z)^2f(z)dz +
 \frac{1}{T - t_0}\int_{-\infty}^{\infty}\sigma(z)^2f(z)dz\right).
\end{eqnarray*}
On the other hand, let us find a suitable bound on the integrated squared-bias of $\widehat{bf}_{N,h}(x)$. Again, since $X^1,\dots,X^N$ are i.i.d. copies of $X$, and since It\^o's integral restricted to $\mathbb H^2$ is a martingale-valued map,
\begin{eqnarray*}
 \mathfrak b(x)
 & = &
 \mathbb E\left[\frac{1}{N(T - t_0)}
 \sum_{i = 1}^{N}\int_{t_0}^{T}K_h(X_{t}^{i} - x)dX_{t}^{i}\right] - b(x)f(x)\\
 & = &
 \frac{1}{T - t_0}
 \mathbb E\left(\int_{t_0}^{T}K_h(X_t - x)dX_t\right) - b(x)f(x)\\
 & = &
 \frac{1}{T - t_0}\left[\mathbb E\left(\int_{t_0}^{T}K_h(X_t - x)b(X_t)dt\right) +
 \mathbb E\left(\int_{t_0}^{T}K_h(X_t - x)\sigma(X_t)dW_t\right)\right] - b(x)f(x)\\
 & = &
 \frac{1}{T - t_0}\int_{t_0}^{T}\mathbb E(K_h(X_t - x)b(X_t))dt - b(x)f(x)
 =\int_{-\infty}^{\infty}K_h(z - x)b(z)f(z)dz - b(x)f(x).
\end{eqnarray*}
Then,
\begin{eqnarray*}
\mathfrak b(x)^2 = ((bf)_h - bf)(x)^2
\quad
\textrm{with}
\quad
(bf)_h = K_h\ast (bf).
\end{eqnarray*}
Therefore, since $f$ is bounded and $b$ belongs to $\mathbb L^2(\mathbb R,f(x)dx)$ by Remark \ref{integrability_functionals_X},
\begin{displaymath}
\int_{-\infty}^{\infty}\mathfrak b(x)^2dx =
\|bf - (bf)_h\|_{2}^{2}.
\end{displaymath}
This concludes the proof.
%


%
\subsection{Proof of Proposition \ref{risk_bound_NW}}
First of all,
\begin{displaymath}
\widehat b_{N,h,h'} - b =
\left[\frac{\widehat{bf}_{N,h} - bf}{\widehat f_{N,h'}} +
\left(\frac{1}{\widehat f_{N,h'}} -\frac{1}{f}\right)bf\right]
\mathbf 1_{\widehat f_{N,h'}(.) > m/2} - 
b\mathbf 1_{\widehat f_{N,h'}(.)\leqslant m/2}.
\end{displaymath}
Then,
\begin{eqnarray*}
 \|\widehat b_{N,h,h'} - b\|_{f,A,B}^{2}
 & \leqslant &
 2\left\|\left[\frac{\widehat{bf}_{N,h} - bf}{\widehat f_{N,h'}} +
 \left(\frac{1}{\widehat f_{N,h'}} -\frac{1}{f}\right)bf\right]
 \mathbf 1_{\widehat f_{N,h'}(.) > m/2}\right\|_{f,A,B}^{2} + 
 2\|b\mathbf 1_{\widehat f_{N,h'}(.)\leqslant m/2}\|_{f,A,B}^{2}.
\end{eqnarray*}
Moreover, for any $x\in [A,B]$, since $f(x) > m$, for every $\omega\in\{\widehat f_{N,h'}(.)\leqslant m/2\}$,
\begin{displaymath}
|f(x) -\widehat f_{N,h'}(x,\omega)|
\geqslant
f(x) -\widehat f_{N,h'}(x,\omega) > m -\frac{m}{2} =\frac{m}{2}.
\end{displaymath}
Thus,
\begin{eqnarray*}
 \|\widehat b_{N,h,h'} - b\|_{f,A,B}^{2}
 & \leqslant &
 \frac{8}{m^2}\|\widehat{bf}_{N,h} - bf\|_{2,f}^{2} +
 \frac{8}{m^2}\|(f -\widehat f_{N,h'})b\|_{f,A,B}^{2} +
 2\|b\mathbf 1_{|f(.) -\widehat f_{N,h'}(.)| > m/2}\|_{f,A,B}^{2}\\
 & \leqslant &
 \frac{8}{m^2}\int_{-\infty}^{\infty}(\widehat{bf}_{N,h} - bf)(x)^2f(x)dx\\
 & & \hspace{1.5cm} + 
 \frac{8}{m^2}\int_{A}^{B}(f(x) -\widehat f_{N,h'}(x))^2b(x)^2f(x)dx\\
 & & \hspace{3cm} + 
 2\int_{A}^{B}b(x)^2f(x)\mathbf 1_{|f(x) -\widehat f_{N,h'}(x)| > m/2}dx.
\end{eqnarray*}
Since $f$ has a sub-Gaussian tail by Assumption \ref{assumption_bound_derivative_density_x}, and since $b$ has at most linear growth because it is Lipschitz continuous from $\mathbb R$ into itself (see Assumption \ref{assumption_b_sigma}), $b^2f$ is bounded on $\mathbb R$. So,
\begin{eqnarray*}
 \|\widehat b_{N,h,h'} - b\|_{f,A,B}^{2}
 & \leqslant &
 \frac{8\|f\|_{\infty}}{m^2}\|\widehat{bf}_{N,h} - bf\|_{2}^{2}\\
 & &
 \quad
 +\frac{8\|b^2f\|_{\infty}}{m^2}\|\widehat f_{N,h'} - f\|_{2}^{2}
 + 2\|b^2f\|_{\infty}\int_{-\infty}^{\infty}\mathbf 1_{|f(x) -\widehat f_{N,h'}(x)| > m/2}dx.
\end{eqnarray*}
Therefore, thanks to Markov's inequality,
\begin{eqnarray*}
 \mathbb E(\|\widehat b_{N,h,h'} - b\|_{f,A,B}^{2})
 & \leqslant &
 \frac{8\|f\|_{\infty}}{m^2}\mathbb E(\|\widehat{bf}_{N,h} - bf\|_{2}^{2})\\
 & &
 \quad
 +\frac{8\|b^2f\|_{\infty}}{m^2}\mathbb E(\|\widehat f_{N,h'} - f\|_{2}^{2})
 +\frac{8\|b^2f\|_{\infty}}{m^2}\int_{-\infty}^{\infty}\mathbb E(|f(x) -\widehat f_{N,h'}(x)|^2)dx\\
 & \leqslant &
 \frac{8(\|f\|_{\infty}\vee\|b^2f\|_{\infty})}{m^2}[\mathbb E(\|\widehat{bf}_{N,h} - bf\|_{2}^{2})
 + 2\mathbb E(\|\widehat f_{N,h'} - f\|_{2}^{2})].
\end{eqnarray*}
Propositions \ref{risk_bound_numerator} and \ref{risk_bound_denominator} allow to conclude.
%


%
\subsection{Proof of Proposition \ref{risk_bound_approximate_denominator}}
First of all, note that
\begin{eqnarray*}
 \mathbb E(\|\widehat f_{n,N,h} - f\|_{2}^{2})
 & \leqslant &
 2\mathbb E(\|\widehat f_{N,h} - f\|_{2}^{2}) +
 2\mathbb E(\|\widehat f_{N,h} -\widehat f_{n,N,h}\|_{2}^{2})\\
 & \leqslant &
 2\left[
 \mathfrak c_{\ref{risk_bound_denominator}}(t_0)h^{2\beta} +\frac{1}{Nh} +
 \int_{-\infty}^{\infty}\mathbb E(\widehat f_{N,h}(x) -\widehat f_{n,N,h}(x))^2dx\right.\\
 & &
 \hspace{4cm}\left. +
 \int_{-\infty}^{\infty}{\rm var}(\widehat f_{N,h}(x) -\widehat f_{n,N,h}(x))dx
 \right]
\end{eqnarray*}
by Proposition \ref{risk_bound_denominator}, and note also that
\begin{displaymath}
\widehat f_{N,h}(x) -\widehat f_{n,N,h}(x) =
\frac{1}{N(T - t_0)}\sum_{i = 1}^{N}\sum_{j = 0}^{n - 1}
\int_{t_j}^{t_{j + 1}}(K_{h,x}(X_{t}^{i}) - K_{h,x}(X_{t_j}^{i}))dt
\end{displaymath}
with $K_{h,x}(.) := K_h(\cdot - x)$. On the one hand, for every $s,u\in [t_0,T]$ such that $s\leqslant u$, by It\^o's formula, Jensen's inequality, the isometry property for It\^o's integral and Remark \ref{integrability_functionals_X},
\begin{eqnarray*}
 \int_{-\infty}^{\infty}
 \mathbb E[(K_{h,x}(X_u) - K_{h,x}(X_s))^2]dx & = &
 \int_{-\infty}^{\infty}
 \mathbb E\left[\left(
 \int_{s}^{u}K_{h,x}'(X_t)dX_t +\frac{1}{2}\int_{s}^{u}K_{h,x}''(X_t)d\langle X\rangle_t\right)^2\right]dx\\
 & \leqslant &
 \mathfrak c_1\int_{-\infty}^{\infty}
 \left[\mathbb E\left[\left(
 \int_{s}^{u}K_{h,x}'(X_t)b(X_t)dt\right)^2\right]\right.\\
 & &
 \hspace{1.5cm} +
 \mathbb E\left[\left(\int_{s}^{u}K_{h,x}''(X_t)\sigma(X_t)^2dt\right)^2\right]\\
 & &
 \hspace{3cm}
 \left. +
 \mathbb E\left[\left(\int_{s}^{u}K_{h,x}'(X_t)\sigma(X_t)dW_t\right)^2\right]\right]dx\\
 & \leqslant &
 \mathfrak c_1\left[
 (u - s)\int_{s}^{u}\mathbb E\left(b(X_t)^2\int_{-\infty}^{\infty}K_{h,x}'(X_t)^2dx\right)dt\right.\\
 & &
 \hspace{1.5cm} +
 (u - s)\int_{s}^{u}\mathbb E\left(\sigma(X_t)^4\int_{-\infty}^{\infty}K_{h,x}''(X_t)^2dx\right)dt\\
 & &
 \hspace{3cm}
 \left. +
 \int_{s}^{u}\mathbb E\left(\sigma(X_t)^2\int_{-\infty}^{\infty}K_{h,x}'(X_t)^2dx\right)dt\right]\\
 & \leqslant &
 \mathfrak c_2\left[
 \frac{(u - s)^2}{h^3} +
 \frac{(u - s)^2}{h^5} +
 \frac{u - s}{h^3}\right]
\end{eqnarray*}
where $\mathfrak c_1$ and $\mathfrak c_2$ are two positive constants not depending on $s$, $u$, $h$, $N$, $n$ and $t_0$. Then,
\begin{eqnarray*}
 & &
 \int_{-\infty}^{\infty}{\rm var}(\widehat f_{n,N,h}(x) -\widehat f_{N,h}(x))dx\\
 & &
 \hspace{2cm}
 =\frac{1}{N(T - t_0)^2}\int_{-\infty}^{\infty}{\rm var}\left[\sum_{j = 0}^{n - 1}
 \int_{t_j}^{t_{j + 1}}(K_{h,x}(X_t) - K_{h,x}(X_{t_j}))dt\right]dx\\
 & &
 \hspace{2cm}
 \leqslant\frac{1}{N(T - t_0)}\sum_{j = 0}^{n - 1}
 \int_{t_j}^{t_{j + 1}}\int_{-\infty}^{\infty}\mathbb E[(K_{h,x}(X_t) - K_{h,x}(X_{t_j}))^2]dxdt
 \leqslant\frac{\mathfrak c_3}{Nnh^3}
\end{eqnarray*}
where the constant $\mathfrak c_3 > 0$ is not depending on $h$, $N$, $n$ and $t_0$. On the other hand, by Assumption \ref{assumption_bound_derivative_density_t},
\begin{eqnarray*}
 |\mathbb E(\widehat f_{N,h}(x) -\widehat f_{n,N,h}(x))| & \leqslant &
 \frac{1}{T - t_0}\sum_{j = 0}^{n - 1}
 \int_{t_j}^{t_{j + 1}}|\mathbb E(K_{h,x}(X_t)) -\mathbb E(K_{h,x}(X_{t_j}))|dt\\
 & \leqslant &
 \frac{1}{T - t_0}\sum_{j = 0}^{n - 1}
 \int_{t_j}^{t_{j + 1}}\int_{-\infty}^{\infty}|K_h(z - x)|\cdot|p_t(x_0,z) - p_{t_j}(x_0,z)|dzdt\\
 & \leqslant &
 \frac{1}{T - t_0}\sum_{j = 0}^{n - 1}\left[
 \int_{t_j}^{t_{j + 1}}(t - t_j)dt\right]\\
 & &
 \hspace{3cm}
 \times\left[\int_{-\infty}^{\infty}|K_h(z - x)|\sup_{u\in [t_0,T]}|\partial_up_u(x_0,z)|dz\right]\\
 & \leqslant &
 \mathfrak c_{\ref{assumption_bound_derivative_density_t},3}
 \frac{T - t_0}{nt_{0}^{q_3}}\int_{-\infty}^{\infty}|K(z)|
 \exp\left[-\mathfrak m_{\ref{assumption_bound_derivative_density_t},3}\frac{(hz + x - x_0)^2}{T}\right]dz.
\end{eqnarray*}
Then, by Jensen's inequality,
\begin{eqnarray*}
 \int_{-\infty}^{\infty}\mathbb E(\widehat f_{N,h}(x) -\widehat f_{n,N,h}(x))^2dx
 & \leqslant &
 \frac{\mathfrak c_5}{n^2t_{0}^{2q_3}}\int_{-\infty}^{\infty}|K(z)|
 \int_{-\infty}^{\infty}
 \exp\left[-2\mathfrak m_{\ref{assumption_bound_derivative_density_t},3}\frac{(hz + x - x_0)^2}{T}\right]dxdz\\
 & = &
 \frac{\mathfrak c_6}{n^2t_{0}^{2q_3}}
\end{eqnarray*}
where
\begin{displaymath}
\mathfrak c_6 =
\mathfrak c_5\|K\|_1
\int_{-\infty}^{\infty}
\exp\left[-2\mathfrak m_{\ref{assumption_bound_derivative_density_t},3}\frac{(x - x_0)^2}{T}\right]dx
\end{displaymath}
and the constant $\mathfrak c_5 > 0$ is not depending on $h$, $N$, $n$ and $t_0$. This concludes the proof.
%


%
\subsection{Proof of Proposition \ref{risk_bound_approximate_numerator}}
The proof of Proposition \ref{risk_bound_approximate_numerator} relies on the two following technical lemmas.
%


%
\begin{lemma}\label{extended_BDG_inequality}
Consider a symmetric and continuous function $\varphi_1 :\mathbb R\rightarrow\mathbb R$ such that $\overline\varphi_1 : z\mapsto z\varphi_1(z)$ belongs to $\mathbb L^2(\mathbb R,dx)$. Consider also $\varphi_2,\psi\in C^0(\mathbb R)$ having polynomial growth. Under Assumptions \ref{assumption_b_sigma} and \ref{assumption_bound_derivative_density_x}, for every $p > 0$, there exists a constant $\mathfrak c_{\ref{extended_BDG_inequality}}(p) > 0$, not depending on $\varphi_1$ and $t_0$, such that for every $s,t\in [t_0,T]$ satisfying $s < t$,
\begin{eqnarray*}
 & &
 \int_{-\infty}^{\infty}\mathbb E\left[
 \left(\int_{s}^{t}\varphi_1(x - X_u)\varphi_2(X_u)dW_u\right)^2\psi(X_t)^2\right]dx\\
 & &
 \hspace{4cm}\leqslant
 \mathfrak c_{\ref{extended_BDG_inequality}}(p)(t - s)\left[
 \|\varphi_1\|_{2}^{2} +\|\overline\varphi_1\|_{2}^{2} +
 \frac{1}{t_{0}^{1/(2p)}}
 \left(\int_{-\infty}^{\infty}\varphi_1(z)^{2p}dz\right)^{1/p}
 \right].
\end{eqnarray*}
\end{lemma}
%


%
\begin{lemma}\label{extended_Jensen_inequality}
Consider $\varphi\in C^0(\mathbb R)$. Under Assumptions \ref{assumption_b_sigma} and \ref{assumption_bound_derivative_density_x}, for every $s,t\in [t_0,T]$ such that $s < t$,
\begin{displaymath}
\int_{-\infty}^{\infty}\mathbb E(K_{h,x}(X_s)\varphi(X_s,X_t))^2dx
\leqslant
\frac{\mathfrak c_{\ref{assumption_bound_derivative_density_x},1}\|K\|_{1}^{2}}{t_{0}^{1/2}}
\mathbb E[\varphi(X_s,X_t)^2].
\end{displaymath}
\end{lemma}
\noindent
The proof of Lemma \ref{extended_BDG_inequality} (resp. Lemma \ref{extended_Jensen_inequality}) is postponed to Subsubection \ref{section_proof_lemma_extended_BDG_inequality} (resp. Subsubsection \ref{section_proof_lemma_extended_Jensen_inequality}).
\\
\\
First of all, note that
\begin{eqnarray*}
 \mathbb E(\|\widehat{bf}_{n,N,h} - bf\|_{2}^{2})
 & \leqslant &
 2\mathbb E(\|\widehat{bf}_{N,h} - bf\|_{2}^{2}) +
 2\mathbb E(\|\widehat{bf}_{N,h} -\widehat{bf}_{n,N,h}\|_{2}^{2})\\
 & \leqslant &
 2\left[
 \|(bf)_h - bf\|_{2}^{2} +\frac{\mathfrak c_{\ref{risk_bound_numerator}}(t_0)}{Nh} +
 \int_{-\infty}^{\infty}\mathbb E(\widehat{bf}_{N,h}(x) -\widehat{bf}_{n,N,h}(x))^2dx\right.\\
 & &
 \hspace{5cm}\left. +
 \int_{-\infty}^{\infty}{\rm var}(\widehat{bf}_{N,h}(x) -\widehat{bf}_{n,N,h}(x))dx
 \right]\\
 & =: &
 2\left[
 \|(bf)_h - bf\|_{2}^{2} +\frac{\mathfrak c_{\ref{risk_bound_numerator}}(t_0)}{Nh} +
 B_{n,N,h} + V_{n,N,h}\right]
\end{eqnarray*}
by Proposition \ref{risk_bound_numerator}, and note also that
\begin{displaymath}
\widehat{bf}_{N,h}(x) -\widehat{bf}_{n,N,h}(x) =
\frac{1}{N(T - t_0)}\sum_{i = 1}^{N}\sum_{j = 0}^{n - 1}
\int_{t_j}^{t_{j + 1}}(K_{h,x}(X_{t}^{i}) - K_{h,x}(X_{t_j}^{i}))dX_{t}^{i}.
\end{displaymath}
The proof is dissected in two steps. The term $V_{n,N,h}$ is controlled in the first step, and then $B_{n,N,h}$ is controlled in the second one.
\\
\\
\textbf{Step 1.} First of all, by Jensen's inequality,
\begin{eqnarray*}
 V_{n,N,h}
 & = &
 \frac{1}{N(T - t_0)^2}\int_{-\infty}^{\infty}
 {\rm var}\left[\sum_{j = 0}^{n - 1}\int_{t_j}^{t_{j + 1}}
 (K_{h,x}(X_t) - K_{h,x}(X_{t_j}))dX_t\right]dx\\
 & \leqslant &
 \frac{2}{N(T - t_0)}\int_{-\infty}^{\infty}\left[
 \sum_{j = 0}^{n - 1}\int_{t_j}^{t_{j + 1}}
 \mathbb E[(K_{h,x}(X_t) - K_{h,x}(X_{t_j}))^2b(X_t)^2]dt\right]dx + V_{n,N,h}^{\sigma}
\end{eqnarray*}
with
\begin{displaymath}
V_{n,N,h}^{\sigma} :=
\frac{2}{N(T - t_0)^2}\int_{-\infty}^{\infty}
\mathbb E\left[\left(\sum_{j = 0}^{n - 1}\int_{t_j}^{t_{j + 1}}
(K_{h,x}(X_t) - K_{h,x}(X_{t_j}))\sigma(X_t)dW_t\right)^2\right]dx.
\end{displaymath}
In order to control $V_{n,N,h}$ as in the proof of Proposition \ref{risk_bound_approximate_denominator}, a preliminary bound on $V_{n,N,h}^{\sigma}$  has to be established via the isometry property of It\^o's integral:
\begin{eqnarray*}
 V_{n,N,h}^{\sigma}
 & = &
 \frac{2}{N(T - t_0)^2}\int_{-\infty}^{\infty}\mathbb E\left[\left(
 \int_{t_0}^{T}\left(\sum_{j = 0}^{n - 1}(K_{h,x}(X_t) - K_{h,x}(X_{t_j}))
 \sigma(X_t)\mathbf 1_{[t_j,t_{j + 1}]}(t)\right)dW_t\right)^2\right]dx\\
 & = &
 \frac{2}{N(T - t_0)^2}\int_{-\infty}^{\infty}\int_{t_0}^{T}
 \mathbb E\left[\left(\sum_{j = 0}^{n - 1}(K_{h,x}(X_t) - K_{h,x}(X_{t_j}))
 \sigma(X_t)\mathbf 1_{[t_j,t_{j + 1}]}(t)\right)^2\right]dtdx\\
 & = &
 \frac{2}{N(T - t_0)}\int_{-\infty}^{\infty}\left(
 \sum_{j = 0}^{n - 1}\int_{t_j}^{t_{j + 1}}
 \mathbb E[(K_{h,x}(X_t) - K_{h,x}(X_{t_j}))^2\sigma(X_t)^2]dt\right)dx.
\end{eqnarray*}
Then, 
\begin{eqnarray*}
 V_{n,N,h}
 & \leqslant &
 \frac{2}{N(T - t_0)}\sum_{j = 0}^{n - 1}\int_{t_j}^{t_{j + 1}}\int_{-\infty}^{\infty}
 \mathbb E[(K_{h,x}(X_t) - K_{h,x}(X_{t_j}))^2b(X_t)^2]dxdt\\
 & &
 \hspace{2cm}
 +\frac{2}{N(T - t_0)}\sum_{j = 0}^{n - 1}\int_{t_j}^{t_{j + 1}}\int_{-\infty}^{\infty}
 \mathbb E[(K_{h,x}(X_t) - K_{h,x}(X_{t_j}))^2\sigma(X_t)^2]dxdt.
\end{eqnarray*}
For $\varphi = b$ or $\varphi =\sigma$, by It\^o's formula,
\begin{eqnarray*}
 & &
 \sum_{j = 0}^{n - 1}\int_{t_j}^{t_{j + 1}}\int_{-\infty}^{\infty}
 \mathbb E[(K_{h,x}(X_t) - K_{h,x}(X_{t_j}))^2\varphi(X_t)^2]dxdt\\
 & &
 \hspace{2cm}\leqslant
 \mathfrak c_1\sum_{j = 0}^{n - 1}\int_{t_j}^{t_{j + 1}}\int_{-\infty}^{\infty}\left[
 \mathbb E\left[\left(\int_{t_j}^{t}
 K_{h,x}'(X_u)b(X_u)du\right)^2\varphi(X_t)^2\right]\right.\\
 & &
 \hspace{4cm} +
 \mathbb E\left[\left(\int_{t_j}^{t}
 K_{h,x}''(X_u)\sigma(X_u)^2du\right)^2\varphi(X_t)^2\right]\\
 & &
 \hspace{6cm}
 \left. +
 \mathbb E\left[\left(\int_{t_j}^{t}
 K_{h,x}'(X_u)\sigma(X_u)dW_u\right)^2\varphi(X_t)^2\right]\right]dxdt
\end{eqnarray*}
where the constant $\mathfrak c_1 > 0$ is not depending on $\varphi$, $h$, $N$, $n$ and $t_0$. Moreover, for every $j\in\{0,\dots,n - 1\}$ and $t\in [t_j,t_{j + 1}]$, by Lemma \ref{extended_BDG_inequality} with $p = 1/(1 -\varepsilon)$,
\begin{eqnarray*}
 \int_{-\infty}^{\infty}\mathbb E\left[\left(\int_{t_j}^{t}
 K_{h,x}'(X_u)\sigma(X_u)dW_u\right)^2\varphi(X_t)^2\right]dx
 & \leqslant &
 \mathfrak c_{\ref{extended_BDG_inequality}}(p)(t - t_j)
 \int_{-\infty}^{\infty}K_h'(z)^2dz\\
 & &
 \hspace{1cm}
 +\mathfrak c_{\ref{extended_BDG_inequality}}(p)(t - t_j)\int_{-\infty}^{\infty}z^2K_h'(z)^2dz\\
 & &
 \hspace{2cm}
 +\frac{\mathfrak c_{\ref{extended_BDG_inequality}}(p)}{t_{0}^{1/(2p)}}(t - t_j)\left(\int_{-\infty}^{\infty}K_h'(z)^{2p}dz\right)^{1/p}\\
 & \leqslant &
 \mathfrak c_2(\varepsilon)(t - t_j)
 \left[1 +\frac{1}{h^3} +\frac{t_{0}^{-(1 -\varepsilon)/2}}{h^{3 +\varepsilon}}\right]
\end{eqnarray*}
where the constant $\mathfrak c_2(\varepsilon) > 0$ depends on $\varepsilon$, but not on $\varphi$, $j$, $t$, $h$, $N$, $n$ and $t_0$. Thus, by Jensen's inequality and Remark \ref{integrability_functionals_X},
\begin{eqnarray*}
 & &
 \sum_{j = 0}^{n - 1}\int_{t_j}^{t_{j + 1}}\int_{-\infty}^{\infty}
 \mathbb E[(K_{h,x}(X_t) - K_{h,x}(X_{t_j}))^2\varphi(X_t)^2]dxdt\\
 & &
 \hspace{1cm}\leqslant
 \mathfrak c_3(\varepsilon)\sum_{j = 0}^{n - 1}\int_{t_j}^{t_{j + 1}}\left[
 (t - t_j)\int_{t_j}^{t}
 \mathbb E\left(
 \varphi(X_t)^2b(X_u)^2\int_{-\infty}^{\infty}K_{h,x}'(X_u)^2dx\right)du\right.\\
 & &
 \hspace{2cm}\left. +
 (t - t_j)\int_{t_j}^{t}
 \mathbb E\left(
 \varphi(X_t)^2\sigma(X_u)^4\int_{-\infty}^{\infty}K_{h,x}''(X_u)^2dx\right)du
 + (t - t_j)\left[1 +\frac{1}{h^3} +
 \frac{t_{0}^{-(1 -\varepsilon)/2}}{h^{3 +\varepsilon}}\right]\right]dt\\
 & &
 \hspace{1cm}\leqslant
 \frac{\mathfrak c_4(\varepsilon)}{\min\{1,t_{0}^{(1 -\varepsilon)/2}\}}(T - t_0)^3
 \left(\frac{1}{n^2h^3} +\frac{1}{n^2h^5} +\frac{1}{nh^{3 +\varepsilon}}\right)
\end{eqnarray*}
where $\mathfrak c_3(\varepsilon)$ and $\mathfrak c_4(\varepsilon)$ are two positive constants depending on $\varepsilon$, but not on $\varphi$, $h$, $N$, $n$ and $t_0$. Therefore,
\begin{displaymath}
V_{n,N,h}\leqslant
\frac{\mathfrak c_5(\varepsilon)}{\min\{1,t_{0}^{(1 -\varepsilon)/2}\}}
\cdot\frac{1}{Nnh^{3 +\varepsilon}}
\end{displaymath}
where the constant $\mathfrak c_5(\varepsilon) > 0$ depends on $\varepsilon$, but not on $h$, $N$, $n$ and $t_0$.
\\
\\
\textbf{Step 2.} First of all, since It\^o's integral restricted to $\mathbb H^2$ is a martingale-valued map, since $K_{h,x}$ is a kernel, by Lemma \ref{extended_Jensen_inequality}, by Assumptions \ref{assumption_bound_derivative_density_x} and \ref{assumption_bound_derivative_density_t}, and since $b$ is Lipschitz continuous,
\begin{eqnarray*}
 & &
 \int_{-\infty}^{\infty}
 \mathbb E(\widehat{bf}_{N,h}(x) -\widehat{bf}_{n,N,h}(x))^2dx\\
 & &
 \hspace{2cm} =
 \int_{-\infty}^{\infty}
 \left(\frac{1}{T - t_0}\sum_{j = 0}^{n - 1}\int_{t_j}^{t_{j + 1}}
 \mathbb E((K_{h,x}(X_t) - K_{h,x}(X_{t_j}))b(X_t))dt\right.\\
 & &
 \hspace{4cm}\left. +
 \frac{1}{T - t_0}\sum_{j = 0}^{n - 1}\mathbb E\left[\int_{t_j}^{t_{j + 1}}
 (K_{h,x}(X_t) - K_{h,x}(X_{t_j}))\sigma(X_t)dW_t\right]\right)^2dx\\
 & &
 \hspace{2cm}\leqslant
 \frac{2}{T - t_0}\sum_{j = 0}^{n - 1}\int_{t_j}^{t_{j + 1}}\int_{-\infty}^{\infty}
 \mathbb E(K_{h,x}(X_t)b(X_t) - K_{h,x}(X_{t_j})b(X_{t_j}))^2dxdt\\
 & &
 \hspace{4cm} +
 \frac{2}{T - t_0}\sum_{j = 0}^{n - 1}\int_{t_j}^{t_{j + 1}}\int_{-\infty}^{\infty}
 \mathbb E(|K_{h,x}(X_{t_j})|\cdot|b(X_t) - b(X_{t_j})|)^2dxdt\\
 & &
 \hspace{2cm}\leqslant
 \frac{2\|K\|_1}{T - t_0}\sum_{j = 0}^{n - 1}\int_{t_j}^{t_{j + 1}}\int_{-\infty}^{\infty}
 \left(\int_{-\infty}^{\infty}|K_{h,x}(z)|dx\right)b(z)^2(p_t(x_0,z) - p_{t_j}(x_0,z))^2dzdt\\
 & &
 \hspace{4cm} +
 \frac{2\mathfrak c_{\ref{assumption_bound_derivative_density_x},1}
 \|K\|_{1}^{2}}{t_{0}^{1/2}(T - t_0)}
 \sum_{j = 0}^{n - 1}\int_{t_j}^{t_{j + 1}}
 \mathbb E[(b(X_t) - b(X_{t_j}))^2]dt\\
 & &
 \hspace{2cm}\leqslant
 \frac{4\|K\|_{1}^{2}}{T - t_0}\sum_{j = 0}^{n - 1}\left[\int_{t_j}^{t_{j + 1}}(t - t_j)^2dt\right]
 \left[\int_{-\infty}^{\infty}b(z)^2\sup_{u\in[t_0,T]}|\partial_up_u(x_0,z)|^2dz\right]\\
 & &
 \hspace{4cm} +
 \frac{2\mathfrak c_{\ref{assumption_bound_derivative_density_x},1}\|K\|_{1}^{2}}{t_{0}^{1/2}(T - t_0)}
 \|b'\|_{\infty}^{2}\sum_{j = 0}^{n - 1}\int_{t_j}^{t_{j + 1}}
 \mathbb E[(X_t - X_{t_j})^2]dt\\
 & &
 \hspace{2cm}\leqslant
 \mathfrak c_6\left(\frac{1}{t_{0}^{2q_3}n^2} +\frac{1}{t_{0}^{1/2}(T - t_0)}\sum_{j = 0}^{n - 1}\int_{t_j}^{t_{j + 1}}
 \mathbb E[(X_t - X_{t_j})^2]dt\right)
\end{eqnarray*}
where the constant $\mathfrak c_6 > 0$ is not depending on $h$, $N$, $n$ and $t_0$. Moreover, for any $j\in\{0,\dots,n - 1\}$ and $t\in [t_j,t_{j + 1}]$,
\begin{displaymath}
X_t - X_{t_j} =
\int_{t_j}^{t}b(X_u)du +\int_{t_j}^{t}\sigma(X_u)dW_u
\end{displaymath}
and then, by Jensen's inequality, the isometry property of It\^o's integral and Remark \ref{integrability_functionals_X},
\begin{eqnarray*}
 \mathbb E[(X_t - X_{t_j})^2] 
 & \leqslant &
 (t - t_j)\int_{t_j}^{t}\mathbb E(b(X_u)^2)du +
 \int_{t_j}^{t}\mathbb E(\sigma(X_u)^2)du\\
 & \leqslant &
 (t - t_j)^2\sup_{u\in[t_0,T]}\mathbb E(b(X_u)^2) +
 (t - t_j)\sup_{u\in[t_0,T]}\mathbb E(\sigma(X_u)^2)
 \leqslant
 \mathfrak c_7(t - t_j)
\end{eqnarray*}
where the constant $\mathfrak c_7 > 0$ is not depending on $j$, $t$, $h$, $N$, $n$ and $t_0$. Therefore,
\begin{displaymath}
\int_{-\infty}^{\infty}\mathbb E(\widehat{bf}_{N,h}(x) -\widehat{bf}_{n,N,h}(x))^2dx
\leqslant
\frac{\mathfrak c_8}{\min\{t_{0}^{1/2},t_{0}^{2q_3}\}}\left(\frac{1}{n^2} +\frac{1}{n}\right)
\end{displaymath}
where the constant $\mathfrak c_8 > 0$ is not depending on $n$, $N$, $h$ and $t_0$.
%


%
\begin{remark}\label{special_case_approximate_numerator_details}
Assume that $b$ and $\sigma$ are bounded. Then, in Step 1, for $\varphi = b$ or $\varphi =\sigma$,
\begin{eqnarray*}
 & &
 \sum_{j = 0}^{n - 1}\int_{t_j}^{t_{j + 1}}\int_{-\infty}^{\infty}
 \mathbb E[(K_{h,x}(X_t) - K_{h,x}(X_{t_j}))^2\varphi(X_t)^2]dxdt\\
 & &
 \hspace{1.5cm}\leqslant
 \mathfrak c_1\sum_{j = 0}^{n - 1}\int_{t_j}^{t_{j + 1}}\int_{-\infty}^{\infty}\left[
 \mathbb E\left[\left(\int_{t_j}^{t}
 K_{h,x}'(X_u)b(X_u)du\right)^2\varphi(X_t)^2\right]\right.\\
 & &
 \hspace{3cm} +
 \mathbb E\left[\left(\int_{t_j}^{t}
 K_{h,x}''(X_u)\sigma(X_u)^2du\right)^2\varphi(X_t)^2\right]\\
 & &
 \hspace{4.5cm}
 \left. +
 \mathbb E\left[\left(\int_{t_j}^{t}
 K_{h,x}'(X_u)\sigma(X_u)dW_u\right)^2\varphi(X_t)^2\right]\right]dxdt\\
 & &
 \hspace{1.5cm}\leqslant
 \mathfrak c_1\|\varphi\|_{\infty}^{2}\sum_{j = 0}^{n - 1}\int_{t_j}^{t_{j + 1}}\int_{-\infty}^{\infty}\left[
 \mathbb E\left[\left(\int_{t_j}^{t}
 K_{h,x}'(X_u)b(X_u)du\right)^2\right]\right.\\
 & &
 \hspace{3cm}\left. +
 \mathbb E\left[\left(\int_{t_j}^{t}
 K_{h,x}''(X_u)\sigma(X_u)^2du\right)^2\right] +
 \mathbb E\left[\left(\int_{t_j}^{t}
 K_{h,x}'(X_u)\sigma(X_u)dW_u\right)^2\right]\right]dxdt.
\end{eqnarray*}
So, in this special case, the bound on $V_{n,N,h}$ is established by using the exact same arguments than in the proof of Proposition \ref{risk_bound_approximate_denominator}. In particular, one can take $\varepsilon = 0$, the additional conditions $K\in\mathbb L^4(\mathbb R,dx)$ and $z\mapsto zK(z)$ belongs to $\mathbb L^2(\mathbb R,dx)$ are not required, and the bound on $V_{n,N,h}$ is of order $1/(Nnh^3)$ and doesn't depend on $t_0$. When $\varphi = b$ or $\varphi =\sigma$ is not bounded, since $\varphi(X_t)$ is not $\sigma(W_u)$-measurable for every $u\in [t_j,t)$ ($j\in\{0,\dots,n - 1\}$), the H\"older inequality has to be used to get a suitable bound on
\begin{displaymath}
\int_{-\infty}^{\infty}\mathbb E\left[\left(
\int_{t_j}^{t}K_{h,x}'(X_u)\sigma(X_u)dW_u\right)^2\varphi(X_t)^2\right]dx
\end{displaymath}
(see the proof of Lemma \ref{extended_BDG_inequality}), and for this reason the variance term in the bound of Proposition \ref{risk_bound_approximate_numerator} is of order $1/(Nnh^{3 +\varepsilon})$ instead of $1/(Nnh^3)$ as when $b$ and $\sigma$ are bounded.
\end{remark}
%


%
\subsubsection{Proof of Lemma \ref{extended_BDG_inequality}}\label{section_proof_lemma_extended_BDG_inequality}
Consider $\varphi(x,z) :=\varphi_1(x - z)\varphi_2(z)$ for every $z\in\mathbb R$, $q > 0$ such that $1/p + 1/q = 1$, and $s,t\in [0,T]$ such that $s < t$. First of all, by the isometry property of It\^o's integral, Burkholder-Davis-Gundy's inequality, H\"older's inequality, Markov's inequality, Remark \ref{integrability_functionals_X}, and the generalized Minkowski inequality,
\begin{eqnarray*}
 & &
 \mathbb E\left[
 \left(\int_{s}^{t}\varphi(x,X_u)dW_u\right)^2\psi(X_t)^2
 \right]\\
 & &
 \hspace{1cm}\leqslant
 x^2\mathbb E\left[
 \left(\int_{s}^{t}\varphi(x,X_u)dW_u\right)^2\mathbf 1_{\psi(X_t)^2\leqslant x^2}
 \right]\\
 & &
 \hspace{2cm} +
 \mathbb E\left[
 \left(\int_{s}^{t}\varphi(x,X_u)dW_u\right)^{2p}
 \right]^{1/p}
 \mathbb E(\psi(X_t)^{4q})^{1/(2q)}\mathbb P(\psi(X_t)^2 > x^2)^{1/(2q)}\\
 & &
 \hspace{1cm}\leqslant
 x^2\int_{s}^{t}\mathbb E[\varphi(x,X_u)^2]du +
 \mathfrak c_1(p)\mathbb E\left[
 \left(\int_{s}^{t}\varphi(x,X_u)^2du\right)^{p}
 \right]^{1/p}
 \mathbb E(\psi(X_t)^{4q})^{1/(2q)}\\
 & &
 \hspace{6cm}\times
 \left[
 \frac{\mathbb E(\psi(X_t)^2)^{1/(2q)}}{x^{1/q}}\mathbf 1_{[-1,1]}(x) +
 \frac{\mathbb E(\psi(X_t)^{4q})^{1/(2q)}}{x^2}\mathbf 1_{\mathbb R\backslash [-1,1]}(x)\right]\\
 & &
 \hspace{1cm}\leqslant
 x^2\int_{s}^{t}\mathbb E[\varphi(x,X_u)^2]du +
 \mathfrak c_2(p)
 \left(\int_{s}^{t}\mathbb E[\varphi(x,X_u)^{2p}]^{1/p}du\right)\left(
 \frac{1}{x^{1/q}}\mathbf 1_{[-1,1]}(x) +
 \frac{1}{x^2}\mathbf 1_{\mathbb R\backslash [-1,1]}(x)\right)
\end{eqnarray*}
where $\mathfrak c_1(p)$ and $\mathfrak c_2(p)$ are two positive constants depending on $p$, but not on $x$, $s$, $t$, $\varphi$ and $t_0$. On the one hand, since $\varphi_2$ has polynomial growth, by Remark \ref{integrability_functionals_X}, for every $u\in [s,t]$,
\begin{eqnarray*}
 \int_{-\infty}^{\infty}x^2
 \mathbb E[\varphi(x,X_u)^2]dx & = &
 \int_{-\infty}^{\infty}
 \int_{-\infty}^{\infty}x^2\varphi_1(x - z)^2\varphi_2(z)^2p_u(x_0,z)dxdz\\
 & = &
 \int_{-\infty}^{\infty}
 \int_{-\infty}^{\infty}(x + z)^2\varphi_1(x)^2\varphi_2(z)^2p_u(x_0,z)dxdz\\
 & \leqslant &
 2\left(\int_{-\infty}^{\infty}x^2\varphi_1(x)^2dx\right)
 \left(\int_{-\infty}^{\infty}\varphi_2(z)^2p_u(x_0,z)dz\right)\\
 & &
 \hspace{2cm} +
 2\left(\int_{-\infty}^{\infty}\varphi_1(x)^2dx\right)
 \left(\int_{-\infty}^{\infty}z^2\varphi_2(z)^2p_u(x_0,z)dz\right)\\
 & \leqslant &
 \mathfrak c_3(\|\overline\varphi_1\|_{2}^{2} +\|\varphi_1\|_{2}^{2})
\end{eqnarray*}
where the constant $\mathfrak c_3 > 0$ is not depending on $u$, $\varphi_1$ and $t_0$. Then,
\begin{displaymath}
\int_{-\infty}^{\infty}x^2\int_{s}^{t}\mathbb E[\varphi(x,X_u)^2]dudx
\leqslant
\mathfrak c_3(t - s)(\|\overline\varphi_1\|_{2}^{2} +\|\varphi_1\|_{2}^{2}).
\end{displaymath}
On the other hand, since $\varphi_2$ has polynomial growth, by Assumption \ref{assumption_bound_derivative_density_x}, for every $x\in\mathbb R$,
\begin{eqnarray*}
 \int_{s}^{t}\mathbb E[\varphi(x,X_u)^{2p}]^{1/p}du & = &
 \int_{s}^{t}\left(\int_{-\infty}^{\infty}\varphi_1(z)^{2p}\varphi_2(z + x)^{2p}p_u(x_0,z + x)dz\right)^{1/p}du\\
 & \leqslant &
 \frac{\mathfrak c_4(p)}{t_{0}^{1/(2p)}}(t - s)
 \left(\int_{-\infty}^{\infty}\varphi_1(z)^{2p}dz\right)^{1/p}
\end{eqnarray*}
where the constant $\mathfrak c_4(p) > 0$ depends on $p$, but not on $x$, $s$, $t$, $\varphi_1$ and $t_0$. Then,
\begin{eqnarray*}
 & &
 \int_{-\infty}^{\infty}\left(\int_{s}^{t}\mathbb E[\varphi(x,X_u)^{2p}]^{1/p}du\right)\left(
 \frac{1}{x^{1/q}}\mathbf 1_{[-1,1]}(x) +
 \frac{1}{x^2}\mathbf 1_{\mathbb R\backslash [-1,1]}(x)\right)dx\\
 & &
 \hspace{7cm}\leqslant
 \frac{\mathfrak c_5(p)}{t_{0}^{1/(2p)}}(t - s)
 \left(\int_{-\infty}^{\infty}\varphi_1(z)^{2p}dz\right)^{1/p}
\end{eqnarray*}
with
\begin{displaymath}
\mathfrak c_5(p) =\mathfrak c_4(p)\left(
\int_{-1}^{1}\frac{dx}{x^{1/q}} +\int_{\mathbb R\backslash [-1,1]}\frac{dx}{x^2}\right) <\infty.
\end{displaymath}
%


%
\subsubsection{Proof of Lemma \ref{extended_Jensen_inequality}}\label{section_proof_lemma_extended_Jensen_inequality}
Consider $s,t\in [t_0,T]$ such that $s < t$, and let $p_{s,t}$ (resp. $p_{t|s}$) be the density of $(X_s,X_t)$ (resp. the conditional density of $X_t$ with respect to $X_s$). Moreover, for the sake of readability, $p_s(x_0,.)$ is denoted by $p_s(.)$ in this proof. By Assumption \ref{assumption_bound_derivative_density_x},
\begin{eqnarray*}
 \mathbb E(K_{h,x}(X_s)\varphi(X_s,X_t))^2
 & = &
 \|K\|_{1}^{2}
 \left[\int_{-\infty}^{\infty}\frac{K_{h,x}(y)}{\|K_{h,x}\|_1}\left(\int_{-\infty}^{\infty}
 \varphi(y,z)p_{t|s}(z|y)dz\right)p_s(y)dy\right]^2\\
 & \leqslant &
 \|K\|_1
 \int_{-\infty}^{\infty}|K_{h,x}(y)|\left(\int_{-\infty}^{\infty}
 \varphi(y,z)p_{t|s}(z|y)dz\right)^2p_s(y)^2dy\\
 & \leqslant &
 \|K\|_1\sup_{s\in [t_0,T]}\left\{
 \sup_{y\in\mathbb R}p_s(y)
 \right\}
 \int_{-\infty}^{\infty}|K_{h,x}(y)|\left(\int_{-\infty}^{\infty}\varphi(y,z)^2p_{t|s}(z|y)dz\right)p_s(y)dy\\
 & \leqslant &
 \frac{\mathfrak c_{\ref{assumption_bound_derivative_density_x},1}\|K\|_1}{t_{0}^{1/2}}
 \int_{-\infty}^{\infty}|K_{h,x}(y)|\int_{-\infty}^{\infty}\varphi(y,z)^2p_{s,t}(y,z)dzdy.
\end{eqnarray*}
Therefore,
\begin{eqnarray*}
 \int_{-\infty}^{\infty}\mathbb E(K_{h,x}(X_s)\varphi(X_s,X_t))^2dx & \leqslant &
 \frac{\mathfrak c_{\ref{assumption_bound_derivative_density_x},1}\|K\|_1}{t_{0}^{1/2}}
 \int_{-\infty}^{\infty}\left(
 \int_{-\infty}^{\infty}|K_{h,x}(y)|dx\right)\left(
 \int_{-\infty}^{\infty}\varphi(y,z)^2p_{s,t}(y,z)dz\right)dy\\
 & = &
 \frac{\mathfrak c_{\ref{assumption_bound_derivative_density_x},1}\|K\|_{1}^{2}}{t_{0}^{1/2}}
 \int_{-\infty}^{\infty}\int_{-\infty}^{\infty}
 \varphi(y,z)^2p_{s,t}(y,z)dzdy\\
 & = &
 \frac{\mathfrak c_{\ref{assumption_bound_derivative_density_x},1}
 \|K\|_{1}^{2}}{t_{0}^{1/2}}
 \mathbb E[\varphi(X_s,X_t)^2].
\end{eqnarray*}
%


%
\subsection{Proof of Theorem \ref{risk_bound_PCO_estimator_bf}}\label{section_proof_PCO}
Throughout this subsection, $\mathcal K$ is a primitive function of the kernel $K$. The proof of Theorem \ref{risk_bound_PCO_estimator_bf}.(1) relies on the following technical lemmas proved at the end of this subsection. The proof of Theorem \ref{risk_bound_PCO_estimator_bf}.(2) is left to the reader because it is similar but simpler than the proof of Theorem \ref{risk_bound_PCO_estimator_bf}.(1) detailed in this subsection.
%


%
\begin{lemma}\label{alternative_expression_numerator}
Under Assumptions \ref{assumption_b_sigma}, \ref{assumption_bound_derivative_density_x}, \ref{assumption_K_1} and \ref{assumption_K_sigma_PCO},
\begin{displaymath}
\frac{1}{T - t_0}\int_{t_0}^{T}K_h(X_t - x)dX_t =\Phi_h(X,x)
\textrm{ $;$ }
\forall x\in\mathbb R\textrm{, }
\forall h > 0,
\end{displaymath}
where $(x,h,\varphi)\mapsto\Phi_h(\varphi,x)$ is the map from $\mathbb R\times (0,\infty)\times C^0([t_0,T];\mathbb R)$ into $\mathbb R$ defined by
\begin{displaymath}
\Phi_h(\varphi,x) :=
\frac{1}{T - t_0}\left[\mathcal K\left(\frac{\varphi(T) - x}{h}\right) -
\mathcal K\left(\frac{\varphi(t_0) - x}{h}\right) -
\frac{1}{2h^2}\int_{t_0}^{T}K'\left(\frac{\varphi(t) - x}{h}\right)\sigma(\varphi(t))^2dt
\right]
\end{displaymath}
for every $x\in\mathbb R$, $h > 0$ and $\varphi\in C^0([t_0,T];\mathbb R)$. Moreover,
\begin{enumerate}
 \item For every $x\in\mathbb R$, $h > 0$ and $\varphi\in C^0([t_0,T];\mathbb R)$,
 \begin{displaymath}
 |\Phi_h(\varphi,x)|\leqslant
 \frac{2\|\mathcal K\|_{\infty}}{T - t_0} +\frac{\|\sigma\|_{\infty}^{2}\|K'\|_{\infty}}{2h^2}.
 \end{displaymath}
 \item For every $h > 0$ and $\varphi\in C^0([t_0,T];\mathbb R)$,
 \begin{displaymath}
 \|\Phi_h(\varphi,.)\|_{2,\delta}^{2}
 \leqslant
 \frac{6\|\mathcal K\|_{\infty}^{2}}{(T - t_0)^{2}} +
 \frac{\|\delta\|_{\infty}\|\sigma\|_{\infty}^{4}\|K'\|_{2}^{2}}{h^3}.
 \end{displaymath}
 \item There exists a deterministic constant $\mathfrak c_{\ref{alternative_expression_numerator},1} > 0$ such that, for every $h,h' > 0$,
 \begin{displaymath}
 \mathbb E(\langle\Phi_h(X^1,.),\Phi_{h'}(X^2,.)\rangle_{2,\delta}^{2})
 \leqslant
 \mathfrak c_{\ref{alternative_expression_numerator},1}\mathfrak m(h')
 \end{displaymath}
 with
 \begin{displaymath}
 \mathfrak m(h') =\mathbb E(\|\Phi_{h'}(X,.)\|_{2,\delta}^{2}).
 \end{displaymath}
 \item There exists a deterministic constant $\mathfrak c_{\ref{alternative_expression_numerator},2} > 0$ such that, for every $h > 0$ and $\varphi\in\mathbb L^2(\mathbb R,dx)$,
 \begin{displaymath}
 \mathbb E(\langle\Phi_h(X,.),\varphi\rangle_{2,\delta}^{2})
 \leqslant\mathfrak c_{\ref{alternative_expression_numerator},2}
 \|\varphi\|_{2,\delta}^{2}.
 \end{displaymath}
 \item There exists a deterministic constant $\mathfrak c_{\ref{alternative_expression_numerator},3} > 0$ such that, for every $h,h'\in\mathcal H_N$,
 \begin{displaymath}
 |\langle\Phi_h(X,.),(bf)_{h'}\rangle_{2,\delta}|
 \leqslant\mathfrak c_{\ref{alternative_expression_numerator},3}
 \quad
 {\rm a.s.}
 \end{displaymath}
\end{enumerate}
\end{lemma}
%


%
\begin{lemma}\label{bound_U_statistics}
Consider
\begin{equation}\label{U_statistic_definition}
U_{h,h'}(N) :=\sum_{i\not= j}
\langle\Phi_h(X^i,.) - (bf)_h,
\Phi_{h'}(X^j,.) - (bf)_{h'}\rangle_{2,\delta}
\textrm{ $;$ }
\forall h,h'\in\mathcal H_N.
\end{equation}
Under Assumptions \ref{assumption_b_sigma}, \ref{assumption_bound_derivative_density_x}, \ref{assumption_K_1} and \ref{assumption_K_sigma_PCO}, there exists a deterministic constant $\mathfrak c_{\ref{bound_U_statistics}} > 0$, not depending on $N$, such that for every $\theta\in (0,1)$ and $\lambda > 0$, with probability larger than $1 - 5.4|\mathcal H_N|e^{-\lambda}$,
\begin{displaymath}
\sup_{h\in\mathcal H_N}
\left\{\frac{|U_{h,h_0}(N)|}{N^2}
-\frac{\theta\mathfrak m(h)}{N}\right\}
\leqslant
\frac{\mathfrak c_{\ref{bound_U_statistics}}(1 +\lambda)^3}{\theta N}
\end{displaymath}
and
\begin{displaymath}
\sup_{h\in\mathcal H_N}
\left\{\frac{|U_{h,h}(N)|}{N^2}
-\frac{\theta\mathfrak m(h)}{N}\right\}
\leqslant
\frac{\mathfrak c_{\ref{bound_U_statistics}}(1 +\lambda)^3}{\theta N}.
\end{displaymath}
\end{lemma}
%


%
\begin{lemma}\label{bound_trace_term}
Consider
\begin{displaymath}
V_h(N) :=\frac{1}{N}\sum_{i = 1}^{N}\|\Phi_h(X^i,.) - (bf)_h\|_{2,\delta}^{2}
\textrm{ $;$ }
\forall h\in\mathcal H_N.
\end{displaymath}
Under Assumptions \ref{assumption_b_sigma}, \ref{assumption_bound_derivative_density_x}, \ref{assumption_K_1} and \ref{assumption_K_sigma_PCO}, there exists a deterministic constant $\mathfrak c_{\ref{bound_trace_term}} > 0$, not depending on $N$, such that for every $\theta\in (0,1)$ and $\lambda > 0$, with probability larger than $1 - 2|\mathcal H_N|e^{-\lambda}$,
\begin{displaymath}
\sup_{h\in\mathcal H_N}\left\{
\frac{1}{N}|V_h(N) -\mathfrak m(h)| -\frac{\theta\mathfrak m(h)}{N}\right\}
\leqslant
\frac{\mathfrak c_{\ref{bound_trace_term}}(1 +\lambda)}{\theta N}.
\end{displaymath}
\end{lemma}
%


%
\begin{lemma}\label{bound_crossed_term}
Consider
\begin{equation}\label{W_statistic_definition}
W_{h,h'}(N) :=
\langle\widehat{bf}_{N,h} - (bf)_h,(bf)_{h'} - bf\rangle_{2,\delta}
\textrm{ $;$ }
\forall h,h'\in\mathcal H_N.
\end{equation}
Under Assumptions \ref{assumption_b_sigma}, \ref{assumption_bound_derivative_density_x}, \ref{assumption_K_1} and \ref{assumption_K_sigma_PCO}, there exists a deterministic constant $\mathfrak c_{\ref{bound_crossed_term}} > 0$, not depending on $N$, such that for every $\theta\in (0,1)$ and $\lambda > 0$, with probability larger than $1 - 2|\mathcal H_N|e^{-\lambda}$,
\begin{eqnarray*}
 \sup_{h\in\mathcal H_N}\{
 |W_{h,h_0}(N)| -\theta\|(bf)_{h_0} - bf\|_{2,\delta}^{2}\}
 & \leqslant &
 \frac{\mathfrak c_{\ref{bound_crossed_term}}(1 +\lambda)^2}{\theta N},\\
 \sup_{h\in\mathcal H_N}\{
 |W_{h_0,h}(N)| -\theta\|(bf)_h - bf\|_{2,\delta}^{2}\}
 & \leqslant &
 \frac{\mathfrak c_{\ref{bound_crossed_term}}(1 +\lambda)^2}{\theta N}\textrm{ and}\\
 \sup_{h\in\mathcal H_N}\{
 |W_{h,h}(N)| -\theta\|(bf)_h - bf\|_{2,\delta}^{2}\}
 & \leqslant &
 \frac{\mathfrak c_{\ref{bound_crossed_term}}(1 +\lambda)^2}{\theta N}.
\end{eqnarray*}
\end{lemma}
%


%
\subsubsection{Steps of the proof}
The proof of Theorem \ref{risk_bound_PCO_estimator_bf}.(1) is dissected in four steps.
\\
\\
\textbf{Step 1.} This first step provides a suitable decomposition of $\|\widehat{bf}_{N,\widehat h} - bf\|_{2,\delta}^{2}$. First,
\begin{eqnarray*}
 \|\widehat{bf}_{N,\widehat h} - bf\|_{2,\delta}^{2}
 & = &
 \|\widehat{bf}_{N,\widehat h} -\widehat{bf}_{N,h_0}\|_{2,\delta}^{2} +
 \|\widehat{bf}_{N,h_0} - bf\|_{2,\delta}^{2}\\
 & &
 \hspace{3cm}
 -2\langle\widehat{bf}_{N,h_0} -\widehat{bf}_{N,\widehat h},
 \widehat{bf}_{N,h_0} - bf\rangle_{2,\delta}.
\end{eqnarray*}
Then, by (\ref{penalty_proposal_bf}) and the definition of ${\rm pen}(.)$ (see (\ref{penalty_PCO})), for any $h\in\mathcal H_N$, 
\begin{eqnarray}
 \|\widehat{bf}_{N,\widehat h} - bf\|_{2,\delta}^{2}
 & \leqslant &
 \|\widehat{bf}_{N,h} -\widehat{bf}_{N,h_0}\|_{2,\delta}^{2} +
 \textrm{pen}(h) -\textrm{pen}(\widehat h)
 \nonumber\\
 & &
 \hspace{2cm}
 +\|\widehat{bf}_{N,h_0} - bf\|_{2,\delta}^{2} -
 2\langle\widehat{bf}_{N,h_0} -\widehat{bf}_{N,\widehat h},
 \widehat{bf}_{N,h_0} - bf\rangle_{2,\delta}
 \nonumber\\
 & \leqslant &
 \|\widehat{bf}_{N,h} - bf\|_{2,\delta}^{2} +
 \textrm{pen}(h) -\textrm{pen}(\widehat h)
 \nonumber\\
 & &
 \hspace{2cm}
 +\|\widehat{bf}_{N,h_0} - bf\|_{2,\delta}^{2} -
 2\langle\widehat{bf}_{N,h} -\widehat{bf}_{N,\widehat h},
 \widehat{bf}_{N,h_0} - bf\rangle_{2,\delta}
 \nonumber\\
 \label{risk_bound_PCO_estimator_bf_1}
 & = &
 \|\widehat{bf}_{N,h} - bf\|_{2,\delta}^{2} -\psi_N(h) +\psi_N(\widehat h)
\end{eqnarray}
where 
\begin{displaymath}
\psi_N(h) :=
2\langle\widehat{bf}_{N,h} - bf,\widehat{bf}_{N,h_0} -
bf\rangle_{2,\delta} -\textrm{pen}(h).
\end{displaymath}
Let's complete the decomposition of $\|\widehat{bf}_{N,\widehat h} - bf\|_{2,\delta}^{2}$ by writing
\begin{displaymath}
\psi_N(h) = 2(\psi_{1,N}(h) + \psi_{2,N}(h) + \psi_{3,N}(h)),
\end{displaymath}
where
\begin{eqnarray*}
 \psi_{1,N}(h) & := &
 \frac{1}{(T - t_0)^2N^2}
 \sum_{i = 1}^{N}\left\langle\int_{t_0}^{T}K_h(X_{s}^{i} -\cdot)dX_{s}^{i},
 \int_{t_0}^{T}K_{h_0}(X_{s}^{i} -\cdot)dX_{s}^{i}\right\rangle_{2,\delta}
 +\frac{U_{h,h_0}(N)}{N^2} -\frac{1}{2}{\rm pen}(h)\\
 & &
 \hspace{2cm} =\frac{U_{h,h_0}(N)}{N^2},\\
 \psi_{2,N}(h) & := & -\frac{1}{N^2}\left(
 \sum_{i = 1}^{N}
 \left\langle\frac{1}{T - t_0}
 \int_{t_0}^{T}K_{h_0}(X_{s}^{i} -\cdot)dX_{s}^{i},(bf)_h\right\rangle_{2,\delta} +\right.\\
 & &
 \hspace{2cm}
 \left. +\sum_{i = 1}^{N}
 \left\langle\frac{1}{T - t_0}
 \int_{t_0}^{T}K_h(X_{s}^{i} -\cdot)dX_{s}^{i},(bf)_{h_0}\right\rangle_{2,\delta}\right)
 +\frac{1}{N}\langle (bf)_{h_0},(bf)_h\rangle_{2,\delta}\textrm{ and}\\
 \psi_{3,N}(h)
 & := &
 W_{h,h_0}(N) + W_{h_0,h}(N) +\langle (bf)_h - bf,(bf)_{h_0} - bf\rangle_{2,\delta}.
\end{eqnarray*}	
\textbf{Step 2.} This step deals with bounds on $\mathbb E(\psi_{j,N}(h))$ and $\mathbb E(\psi_{j,N}(\widehat h))$ for $j = 1,2,3$.
\begin{itemize}
 \item By Lemma \ref{bound_U_statistics}, for any $\lambda > 0$ and $\theta\in (0,1)$, with probability larger than $1 - 5.4|\mathcal H_N|e^{-\lambda}$,
 \begin{displaymath}
 |\psi_{1,N}(h)|
 \leqslant \frac{\theta\mathfrak m(h)}{N} +
 \frac{\mathfrak c_{\ref{bound_U_statistics}}(1 +\lambda)^3}{\theta N}
 \quad {\rm and}\quad
 |\psi_{1,N}(\widehat h)|
 \leqslant\frac{\theta\mathfrak m(\widehat h)}{N} +
 \frac{\mathfrak c_{\ref{bound_U_statistics}}(1 +\lambda)^3}{\theta N}.
 \end{displaymath}
 \item On the one hand, for any $h,h'\in\mathcal H_N$, consider 
 \begin{displaymath}
 \Psi_{2,N}(h,h') :=
 \frac{1}{N}\sum_{i = 1}^{N}
 \langle\Phi_h(X^i,.),(bf)_{h'}\rangle_{2,\delta}.
 \end{displaymath}
 By Lemma \ref{alternative_expression_numerator},
 \begin{displaymath}
 |\Psi_{2,N}(h,h')|
 \leqslant
 \frac{1}{N}\sum_{i = 1}^{N}
 \left|\int_{-\infty}^{\infty}
 \Phi_h(X^i,x)(bf)_{h'}(x)\delta(x)dx\right|\\
 \leqslant
 \mathfrak c_{\ref{alternative_expression_numerator},3}
 \quad {\rm a.s.}
 \end{displaymath}
 On the other hand,
 \begin{displaymath}
 |\langle (bf)_h,(bf)_{h_0}\rangle_{2,\delta}| 
 \leqslant
 \|\delta\|_{\infty}\|K_h\ast (bf)\|_{\infty}\|K_{h_0}\ast (bf)\|_1
 \leqslant
 \|\delta\|_{\infty}\|K\|_{1}^{2}\|bf\|_{\infty}\|bf\|_1.
 \end{displaymath}	
 Then, there exists a deterministic constant $\mathfrak c_1 > 0$, not depending on $N$ and $h$, such that
 \begin{displaymath}
 |\psi_{2,N}(h)|\leqslant \frac{\mathfrak c_1}{N}
 \quad {\rm and}\quad
 |\psi_{2,N}(\widehat h)|\leqslant
 \sup_{h'\in\mathcal H_N}
 |\psi_{2,N}(h')|\leqslant
 \frac{\mathfrak c_1}{N}
 \quad {\rm a.s.}
 \end{displaymath}
\item By Lemma \ref{bound_crossed_term} and Cauchy-Schwarz's inequality, with probability larger that $1 - |\mathcal H_N|e^{-\lambda}$,
 \begin{eqnarray*}
  |\psi_{3,N}(h)|
  & \leqslant &
  \frac{\theta}{4}(\|(bf)_h - bf\|_{2,\delta}^{2} +
  \|(bf)_{h_0} - bf\|_{2,\delta}^{2}) +
  \frac{8\mathfrak c_{\ref{bound_crossed_term}}(1 +\lambda)^2}{\theta N}\\
  & &
  \hspace{2cm}
  + 2\times\frac{1}{2^{1/2}}\left(\frac{\theta}{2}\right)^{1/2}
  \|(bf)_h - bf\|_{2,\delta}\times
  \frac{1}{2^{1/2}}\left(\frac{2}{\theta}\right)^{1/2}\|(bf)_{h_0} - bf\|_{2,\delta}\\
  & \leqslant &
  \frac{\theta}{2}\|(bf)_h - bf\|_{2,\delta}^{2} +
  \left(\frac{\theta}{4} +\frac{1}{\theta}\right)
  \|(bf)_{h_0} - bf\|_{2,\delta}^{2} +
  \frac{8\mathfrak c_{\ref{bound_crossed_term}}(1 +\lambda)^2}{\theta N}
 \end{eqnarray*}	
 and
 \begin{displaymath}
 |\psi_{3,N}(\widehat h)|\leqslant
 \frac{\theta}{2}\|(bf)_{\widehat h} - bf\|_{2,\delta}^{2} +
 \left(\frac{\theta}{4} +\frac{1}{\theta}\right)
 \|(bf)_{h_0} - bf\|_{2,\delta}^{2} +
 \frac{8\mathfrak c_{\ref{bound_crossed_term}}(1 +\lambda)^2}{\theta N}.
 \end{displaymath}
\end{itemize}
\textbf{Step 3.} Let us establish that there exist two deterministic constants $\mathfrak c_2,\overline{\mathfrak c}_2 > 0$, not depending on $N$ and $\theta$, such that with probability larger than $1 -\overline{\mathfrak c}_2|\mathcal H_N|e^{-\lambda}$,
\begin{displaymath}
\sup_{h\in\mathcal H_N}\left\{
\|\widehat{bf}_{N,h} - bf\|_{2,\delta}^{2} -
(1 +\theta)\left(\|(bf)_h - bf\|_{2,\delta}^{2} +
\frac{\mathfrak m(h)}{N}\right)\right\}
\leqslant
\frac{\mathfrak c_2(1 +\lambda)^3}{\theta N}
\end{displaymath}
and
\begin{displaymath}
\sup_{h\in\mathcal H_N}\left\{
\|(bf)_h - bf\|_{2,\delta}^{2} +\frac{\mathfrak m(h)}{N} -
\frac{1}{1 -\theta}\|\widehat{bf}_{N,h} - bf\|_{2,\delta}^{2}\right\}
\leqslant
\frac{\mathfrak c_2(1 +\lambda)^3}{\theta(1 -\theta)N}.
\end{displaymath}
On the one hand, note that
\begin{displaymath}
\|\widehat{bf}_{N,h} - bf\|_{2,\delta}^{2} -
(1 +\theta)\left(\|(bf)_h - bf\|_{2,\delta}^{2} +\frac{\mathfrak m(h)}{N}\right)
\end{displaymath}
can be written
\begin{displaymath}
\|\widehat{bf}_{N,h} - (bf)_h\|_{2,\delta}^{2} -\frac{(1 +\theta)\mathfrak m(h)}{N} +
2W_h(N) -\theta\|(bf)_h - bf\|_{2,\delta}^{2},
\end{displaymath}
where $W_h(N) := W_{h,h}(N)$ (see (\ref{W_statistic_definition})). Moreover, for any $h\in\mathcal H_N$,
\begin{equation}\label{risk_bound_PCO_estimator_bf_2}
\|\widehat{bf}_{N,h} - (bf)_h\|_{2,\delta}^{2} =
\frac{U_h(N)}{N^2} +\frac{V_h(N)}{N}
\end{equation}
with $U_h(N) = U_{h,h}(N)$ (see (\ref{U_statistic_definition})). So, with probability larger than $1 -\overline{\mathfrak c}_2|\mathcal H_N|e^{-\lambda}$,
\begin{displaymath}
\sup_{h\in\mathcal H_N}\left\{
\left|\|\widehat{bf}_{N,h} - (bf)_h\|_{2,\delta}^{2} -
\frac{\mathfrak m(h)}{N}\right|
-\frac{\theta\mathfrak m(h)}{N}\right\}
\leqslant
\frac{2(\mathfrak c_{\ref{bound_U_statistics}} +\mathfrak c_{\ref{bound_trace_term}})(1 +\lambda)^3}{\theta N}
\end{displaymath}
by Lemmas \ref{bound_U_statistics} and \ref{bound_trace_term}, and then
\begin{displaymath}
\sup_{h\in\mathcal H_N}\left\{
\|\widehat{bf}_{N,h} - bf\|_{2,\delta}^{2} - (1 +\theta)\left(\|(bf)_h - bf\|_{2,\delta}^{2} +\frac{\mathfrak m(h)}{N}\right)\right\}
\leqslant\frac{\mathfrak c_2(1 +\lambda)^3}{\theta N}
\end{displaymath}
by Lemma \ref{bound_crossed_term}. On the other hand, for any $h\in\mathcal H_N$,
\begin{displaymath}
\|(bf)_h - bf\|_{2,\delta}^{2} =
\|\widehat{bf}_{N,h} - bf\|_{2,\delta}^{2}
-\|\widehat{bf}_{N,h} - (bf)_h\|_{2,\delta}^{2} - W_h(N).
\end{displaymath}
Then,
\begin{displaymath}
(1 -\theta)\left(\|(bf)_h - bf\|_{2,\delta}^{2} +\frac{\mathfrak m(h)}{N}\right)
-\|\widehat{bf}_{N,h} - bf\|_{2,\delta}^{2}
\leqslant
|W_h(N)| -\theta\|(bf)_h - bf\|_{2,\delta}^{2} +
\Lambda_h(N) -\frac{\theta\mathfrak m(h)}{N}
\end{displaymath}
where
\begin{displaymath}
\Lambda_h(N) :=\left|
\|\widehat{bf}_{N,h} - (bf)_h\|_{2,\delta}^{2} -
\frac{\mathfrak m(h)}{N}\right|.
\end{displaymath}
By Equality (\ref{risk_bound_PCO_estimator_bf_2}),
\begin{displaymath}
\Lambda_h(N) =
\left|\frac{U_h(N)}{N^2} +\frac{V_h(N)}{N} -\frac{\mathfrak m(h)}{N}\right|.
\end{displaymath}
By Lemmas \ref{bound_trace_term} and \ref{bound_U_statistics}, there exist two deterministic constants $\mathfrak c_3,\overline{\mathfrak c}_3 > 0$, not depending $N$ and $\theta$, such that with probability larger than $1 -\overline{\mathfrak c}_3|\mathcal H_N|e^{-\lambda}$,
\begin{displaymath}
\sup_{h\in\mathcal H_N}\left\{
\Lambda_h(N) -\theta\frac{\mathfrak m(h)}{N}\right\}
\leqslant
\frac{\mathfrak c_3(1 +\lambda)^3}{\theta N}.
\end{displaymath}
By Lemma \ref{bound_crossed_term}, with probability larger than $1 - 2|\mathcal H_N|e^{-\lambda}$,
\begin{displaymath}
\sup_{h\in\mathcal H_N}\{|W_h(N)| -\theta\|(bf)_h - bf\|_{2,\delta}^{2}
\}\leqslant
\frac{\mathfrak c_{\ref{bound_crossed_term}}(1 +\lambda)^2}{\theta N}.
\end{displaymath}
Therefore, with probability larger than $1 -\overline{\mathfrak c}_2|\mathcal H_N|e^{-\lambda}$,
\begin{displaymath}
\sup_{h\in\mathcal H_N}\left\{
\|(bf)_h - bf\|_{2,\delta}^{2} +\frac{\mathfrak m(h)}{N} -
\frac{1}{1 -\theta}\|\widehat{bf}_{N,h} - bf\|_{2,\delta}^{2}\right\}
\leqslant
\frac{\mathfrak c_2(1 +\lambda)^3}{\theta(1 -\theta)N}.
\end{displaymath}
\textbf{Step 4.} By step 2, there exist two deterministic constants $\mathfrak c_4,\overline{\mathfrak c}_4 > 0$, not depending on $N$, $\theta$, $h$ and $h_0$, such that with probability larger than $1 -\overline{\mathfrak c}_4|\mathcal H_N|e^{-\lambda}$,
\begin{displaymath}
|\psi_N(h)|
\leqslant\theta\left(\|(bf)_h - bf\|_{2,\delta}^{2} +\frac{\mathfrak m(h)}{N}\right) +
\left(\frac{\theta}{2} +\frac{2}{\theta}\right)\|(bf)_{h_0} - bf\|_{2,\delta}^{2} +
\frac{\mathfrak c_4(1 +\lambda)^3}{\theta N}
\end{displaymath} 
and
\begin{displaymath}
|\psi_N(\widehat h)|
\leqslant\theta\left(\|(bf)_{\widehat h} - bf\|_{2,\delta}^{2} +
\frac{\mathfrak m(\widehat h)}{N}\right) +
\left(\frac{\theta}{2} +\frac{2}{\theta}\right)\|(bf)_{h_0} - bf\|_{2,\delta}^{2} +
\frac{\mathfrak c_4(1 +\lambda)^3}{\theta N}.
\end{displaymath} 
Then, by step 3, there exist two deterministic constants $\mathfrak c_5,\overline{\mathfrak c}_5 > 0$, not depending on $N$, $\theta$, $h$ and $h_0$, such that with probability larger than $1 -\overline{\mathfrak c}_5|\mathcal H_N|e^{-\lambda}$,
\begin{displaymath}
|\psi_N(h)|
\leqslant
\frac{\theta}{1 -\theta}
\|\widehat{bf}_{N,h} - bf\|_{2,\delta}^{2} +
\left(\frac{\theta}{2} +\frac{2}{\theta}\right)\|(bf)_{h_0} - bf\|_{2,\delta}^{2} +
\mathfrak c_5\left(\frac{1}{\theta} +\frac{1}{1 -\theta}\right)
\frac{(1 +\lambda)^3}{N}
\end{displaymath}	
and
\begin{displaymath}
|\psi_N(\widehat h)|
\leqslant
\frac{\theta}{1 -\theta}
\|\widehat{bf}_{N,\widehat h} - bf\|_{2,\delta}^{2} +
\left(\frac{\theta}{2} +\frac{2}{\theta}\right)\|(bf)_{h_0} - bf\|_{2,\delta}^{2} +
\mathfrak c_5\left(\frac{1}{\theta} +\frac{1}{1 -\theta}\right)
\frac{(1 +\lambda)^3}{N}.
\end{displaymath}	
By the decomposition \eqref{risk_bound_PCO_estimator_bf_1}, there exist two deterministic constants $\mathfrak c_6,\overline{\mathfrak c}_6 > 0$, not depending on $N$, $\theta$, $h$ and $h_0$, such that with probability larger than $1 -\overline{\mathfrak c}_6|\mathcal H_N|e^{-\lambda}$,
\begin{eqnarray*}
 \|\widehat{bf}_{N,\widehat h} - bf\|_{2,\delta}^{2}
 & \leqslant &
 \|\widehat{bf}_{N,h} - bf\|_{2,\delta}^{2} + |\psi_N(h)| + |\psi_N(\widehat h)|\\
 & \leqslant &
 \left(1 +\frac{\theta}{1 -\theta}\right)
 \|\widehat{bf}_{N,h} - bf\|_{2,\delta}^{2} +
 \frac{\theta}{1 -\theta}\|\widehat{bf}_{N,\widehat h} - bf\|_{2,\delta}^{2}\\
 & &
 \hspace{2cm} +
 \frac{\mathfrak c_6}{\theta}\|(bf)_{h_0} - bf\|_{2,\delta}^{2} +
 \frac{\mathfrak c_6}{\theta(1 -\theta)}\cdot\frac{(1 +\lambda)^3}{N}.
\end{eqnarray*}	
This concludes the proof.
%


%
\subsubsection{Proof of Lemma \ref{alternative_expression_numerator}}
First of all, for any $x\in\mathbb R$ and $h > 0$, by It\^o's formula,
\begin{displaymath}
\mathcal K\left(\frac{X_T - x}{h}\right) =
\mathcal K\left(\frac{X_{t_0} - x}{h}\right) +
\int_{t_0}^{T}K_h(X_t - x)dX_t +
\frac{1}{2h^2}\int_{t_0}^{T}K'\left(\frac{X_t - x}{h}\right)d\langle X\rangle_t.
\end{displaymath}
So,
\begin{eqnarray*}
 \int_{t_0}^{T}K_h(X_t - x)dX_t & = &
 \mathcal K\left(\frac{X_T - x}{h}\right) -
 \mathcal K\left(\frac{X_{t_0} - x}{h}\right)\\
 & & \hspace{2cm} -
 \frac{1}{2h^2}\int_{t_0}^{T}K'\left(\frac{X_t - x}{h}\right)\sigma(X_t)^2dt
 = (T - t_0)\Phi_h(X,x).
\end{eqnarray*}
Lemma \ref{alternative_expression_numerator}.(1) is a straightforward consequence of the previous equality and Lemma \ref{alternative_expression_numerator}.(2) is easy to establish: for every $h > 0$ and $\varphi\in C^0([t_0,T];\mathbb R)$,
\begin{eqnarray*}
 (T - t_0)^2\|\Phi_h(\varphi,.)\|_{2,\delta}^{2}
 & \leqslant &
 2\int_{-\infty}^{\infty}\mathcal K\left(\frac{\varphi(T) - x}{h}\right)^2\delta(x)dx +
 4\int_{-\infty}^{\infty}\mathcal K\left(\frac{\varphi(t_0) - x}{h}\right)^2\delta(x)dx\\
 & &
 \hspace{3cm} +
 \frac{1}{h^4}\int_{-\infty}^{\infty}
 \left[\int_{t_0}^{T}K'\left(\frac{\varphi(t) - x}{h}\right)
 \sigma(\varphi(t))^2dt\right]^2\delta(x)dx\\
 & \leqslant &
 6\|\mathcal K\|_{\infty}^{2} +
 \frac{T - t_0}{h^4}\int_{t_0}^{T}\sigma(\varphi(t))^4
 \int_{-\infty}^{\infty}K'\left(\frac{\varphi(t) - x}{h}\right)^2\delta(x)dxdt\\
 & \leqslant &
 6\|\mathcal K\|_{\infty}^{2} +
 \frac{(T - t_0)^2\|\delta\|_{\infty}\|\sigma\|_{\infty}^{4}\|K'\|_{2}^{2}}{h^3}.
\end{eqnarray*}
Let us prove Lemma \ref{alternative_expression_numerator}.(3,4,5). First, for any $h,h' > 0$,
\begin{eqnarray*}
 & &
 \mathbb E(\langle\Phi_h(X^1,.),\Phi_{h'}(X^2,.)\rangle_{2,\delta}^{2})\\
 & &
 \hspace{2cm} =
 \frac{1}{(T - t_0)^4}\mathbb E\left[\left(\int_{-\infty}^{\infty}
 \left(\int_{t_0}^{T}K_h(X_{t}^{1} - x)dX_{t}^{1}\right)
 \left(\int_{t_0}^{T}K_{h'}(X_{t}^{2} - x)dX_{t}^{2}\right)\delta(x)dx\right)^2
 \right]\\
 & &
 \hspace{2cm}\leqslant
 \frac{2}{(T - t_0)^4}(\mathbb E(A_{h,h'}^{2})
 +\mathbb E(B_{h,h'}^{2}))
\end{eqnarray*}
with
\begin{eqnarray*}
 A_{h,h'} & := &
 \int_{-\infty}^{\infty}
 \left(\int_{t_0}^{T}K_h(X_{t}^{1} - x)\sigma(X_{t}^{1})dW_{t}^{1}\right)
 \left(\int_{t_0}^{T}K_{h'}(X_{t}^{2} - x)dX_{t}^{2}\right)\delta(x)dx
 \textrm{ and}\\
 B_{h,h'} & := &
 \int_{-\infty}^{\infty}
 \left(\int_{t_0}^{T}K_h(X_{t}^{1} - x)b(X_{t}^{1})dt\right)
 \left(\int_{t_0}^{T}K_{h'}(X_{t}^{2} - x)dX_{t}^{2}\right)
 \delta(x)dx.
\end{eqnarray*}
{\bf Bound on $\mathbb E(A_{h,h'}^{2})$.} Since $(X^1,W^1)$ and $X^2$ are independent,
\begin{eqnarray*}
 \mathbb E(A_{h,h'}^{2}) & = &
 \int_{-\infty}^{\infty}\int_{-\infty}^{\infty}
 \mathbb E\left[\left(\int_{t_0}^{T}K_h(X_{t}^{1} - x)\sigma(X_{t}^{1})dW_{t}^{1}\right)
 \left(\int_{t_0}^{T}K_h(X_{t}^{1} - y)\sigma(X_{t}^{1})dW_{t}^{1}\right)\right]\\
 & &
 \hspace{2cm}\times
 \mathbb E\left[\left(\int_{t_0}^{T}K_{h'}(X_{t}^{2} - x)dX_{t}^{2}\right)
 \left(\int_{t_0}^{T}K_{h'}(X_{t}^{2} - y)dX_{t}^{2}\right)\right]\delta(x)\delta(y)dxdy.
\end{eqnarray*}
On the one hand, for every $x,y\in\mathbb R$, by the isometry property of It\^o's integral and the definition of $f$,
\begin{eqnarray*}
 & &
 \mathbb E\left[\left(\int_{t_0}^{T}K_h(X_{t}^{1} - x)\sigma(X_{t}^{1})dW_{t}^{1}\right)
 \left(\int_{t_0}^{T}K_h(X_{t}^{1} - y)\sigma(X_{t}^{1})dW_{t}^{1}\right)\right]\\
 & &
 \hspace{1cm}
 =\int_{t_0}^{T}\mathbb E(K_h(X_{t}^{1} - x)K_h(X_{t}^{1} - y)
 \sigma(X_{t}^{1})^2)dt
 = (T - t_0)\int_{-\infty}^{\infty}K_h(z - x)K_h(z - y)\sigma(z)^2f(z)dz.
\end{eqnarray*}
Then,
\begin{eqnarray*}
 \mathbb E(A_{h,h'}^{2}) & = &
 (T - t_0)
 \int_{-\infty}^{\infty}\int_{-\infty}^{\infty}
 \int_{-\infty}^{\infty}K_h(z - x)K_h(z - y)\sigma(z)^2f(z)\\
 & &
 \hspace{2cm}\times
 \mathbb E\left[
 \left(\int_{t_0}^{T}K_{h'}(X_{t}^{2} - x)dX_{t}^{2}\right)
 \left(\int_{t_0}^{T}K_{h'}(X_{t}^{2} - y)dX_{t}^{2}\right)\right]\delta(x)\delta(y)dxdydz\\
 & = &
 (T - t_0)\int_{-\infty}^{\infty}\sigma(z)^2f(z)
 \mathbb E\left[\left(\int_{-\infty}^{\infty}
 K_h(z - x)\delta(x)
 \int_{t_0}^{T}K_{h'}(X_{t}^{2} - x)dX_{t}^{2}dx\right)^2\right]dz.
\end{eqnarray*}
On the other hand, for every $z\in\mathbb R$, $x\mapsto |K_h(z - x)|/\|K\|_1$ is a density function. Then, by Jensen's inequality,
\begin{eqnarray*}
 \mathbb E(A_{h,h'}^{2}) & \leqslant &
 (T - t_0)\|K\|_1
 \int_{-\infty}^{\infty}\sigma(z)^2f(z)\int_{-\infty}^{\infty}
 |K_h(z - x)|\delta(x)^2\mathbb E\left[\left(
 \int_{t_0}^{T}K_{h'}(X_{t}^{2} - x)dX_{t}^{2}\right)^2\right]dxdz\\
 & \leqslant &
 (T - t_0)\|\sigma^2f\|_{\infty}
 \|K\|_{1}^{2}\|\delta\|_{\infty}\int_{-\infty}^{\infty}
 \delta(x)\mathbb E\left[\left(
 \int_{t_0}^{T}K_{h'}(X_{t}^{2} - x)dX_{t}^{2}\right)^2\right]dx\\
 & &
 \hspace{8cm}\leqslant
 (T - t_0)^3\|\sigma^2f\|_{\infty}
 \|K\|_{1}^{2}\|\delta\|_{\infty}\mathfrak m(h').
\end{eqnarray*}
{\bf Bound on $\mathbb E(B_{h,h'}^{2})$.} Since $x\mapsto |K_h(X_t(\omega) - x)|/\|K\|_1$ is a density function for every $(t,\omega)\in [t_0,T]\times\Omega$, by Jensen's inequality,
\begin{eqnarray*}
 \mathbb E(B_{h,h'}^{2}) & = &
 \mathbb E\left[\left(\int_{t_0}^{T}\int_{-\infty}^{\infty}K_h(X_{t}^{1} - x)b(X_{t}^{1})\delta(x)
 \int_{t_0}^{T}K_{h'}(X_{s}^{2} - x)dX_{s}^{2}dxdt\right)^2\right]\\
 & \leqslant &
 (T - t_0)\|K\|_1
 \int_{t_0}^{T}\int_{-\infty}^{\infty}\mathbb E(|K_h(X_{t}^{1} - x)|b(X_{t}^{1})^2)
 \delta(x)^2\mathbb E\left[
 \left(\int_{t_0}^{T}K_{h'}(X_{s}^{2} - x)dX_{s}^{2}\right)^2\right]dxdt\\
 & = &
 (T - t_0)^2\|K\|_1
 \int_{-\infty}^{\infty}\left(\int_{-\infty}^{\infty}
 |K_h(z - x)|b(z)^2f(z)dz\right)
 \delta(x)^2\mathbb E\left[
 \left(\int_{t_0}^{T}K_{h'}(X_{s}^{2} - x)dX_{s}^{2}\right)^2\right]dx\\
 & \leqslant &
 (T - t_0)^2\|b^2f\|_{\infty}
 \|K\|_{1}^{2}\|\delta\|_{\infty}\int_{-\infty}^{\infty}
 \delta(x)\mathbb E\left[
 \left(\int_{t_0}^{T}K_{h'}(X_{s}^{2} - x)dX_{s}^{2}\right)^2\right]dx\\
 & &
 \hspace{8cm}\leqslant
 (T - t_0)^4\|b^2f\|_{\infty}\|K\|_{1}^{2}\|\delta\|_{\infty}\mathfrak m(h').
\end{eqnarray*}
Now, for any $h > 0$ and $\varphi\in\mathbb L^2(\mathbb R,dx)$,
\begin{eqnarray*}
 \mathbb E(\langle\Phi_h(X,.),\varphi\rangle_{2,\delta}^{2}) & = &
 \frac{1}{(T - t_0)^2}
 \mathbb E\left[\left(\int_{-\infty}^{\infty}
 \varphi(x)\delta(x)\int_{t_0}^{T}K_h(X_t - x)dX_tdx\right)^2
 \right]\\
 & \leqslant &
 \frac{2}{(T - t_0)^2}(\mathbb E(C_{h}^{2}) +\mathbb E(D_{h}^{2}))
\end{eqnarray*}
with
\begin{eqnarray*}
 C_h & := &
 \int_{-\infty}^{\infty}
 \varphi(x)\delta(x)\int_{t_0}^{T}K_h(X_t - x)\sigma(X_t)dW_tdx
 \textrm{ and}\\
 D_h & := &
 \int_{-\infty}^{\infty}
 \varphi(x)\delta(x)\int_{t_0}^{T}K_h(X_t - x)b(X_t)dtdx.
\end{eqnarray*}
{\bf Bound on $\mathbb E(C_{h}^{2})$.} By the isometry property of It\^o's integral and the definition of $f$,
\begin{eqnarray*}
 \mathbb E(C_{h}^{2}) & = &
 \int_{-\infty}^{\infty}\int_{-\infty}^{\infty}
 \varphi(x)\varphi(y)\delta(x)\delta(y)\int_{t_0}^{T}\mathbb E(K_h(X_t - x)K_h(X_t - y)
 \sigma(X_t)^2)dtdxdy\\
 & = &
 (T - t_0)\int_{-\infty}^{\infty}\int_{-\infty}^{\infty}
 \varphi(x)\varphi(y)\delta(x)\delta(y)
 \int_{-\infty}^{\infty}K_h(z - x)K_h(z - y)\sigma(z)^2f(z)dzdxdy\\
 & = &
 (T - t_0)\int_{-\infty}^{\infty}(K_h\ast (\varphi\delta))(z)^2\sigma(z)^2f(z)dz
 \leqslant
 (T - t_0)\|\sigma^2f\|_{\infty}\|K\|_{1}^{2}\|\delta\|_{\infty}\|\varphi\|_{2,\delta}^{2}.
\end{eqnarray*}
{\bf Bound on $\mathbb E(D_{h}^{2})$.} By the definition of $f$,
\begin{eqnarray*}
 \mathbb E(D_{h}^{2}) & = &
 \mathbb E\left[\left(\int_{t_0}^{T}
 b(X_t)\int_{-\infty}^{\infty}
 K_h(X_t - x)\varphi(x)\delta(x)dxdt\right)^2\right]
 \leqslant
 (T - t_0)\int_{t_0}^{T}\mathbb E(b(X_t)^2(K_h\ast (\varphi\delta))(X_t)^2)dt\\
 & \leqslant &
 (T - t_0)^2\int_{-\infty}^{\infty}(K_h\ast (\varphi\delta))(z)^2b(z)^2f(z)dz
 \leqslant (T - t_0)^2\|b^2f\|_{\infty}
 \|K\|_{1}^{2}\|\delta\|_{\infty}\|\varphi\|_{2,\delta}^{2}.
\end{eqnarray*}
Finally, since $X$ is a semi-martingale and since the map $(t,\omega,x)\mapsto K_h(X_t(\omega) - x)(bf)_{h'}(x)\delta(x)$ is measurable and bounded for any $h,h'\in\mathcal H_N$, by the stochastic Fubini theorem and It\^o's formula,
\begin{eqnarray*}
 (T - t_0)\langle\Phi_h(X,.),(bf)_{h'}\rangle_{2,\delta}
 & = &
 (T - t_0)\int_{-\infty}^{\infty}
 \Phi_h(X,x)(bf)_{h'}(x)\delta(x)dx\\
 & = &
 \int_{t_0}^{T}
 \int_{-\infty}^{\infty}
 K_h(X_t - x)(bf)_{h'}(x)\delta(x)dxdX_t
 \quad {\rm a.s.}\\
 & = &
 \int_{t_0}^{T}[K_h\ast ((bf)_{h'}\delta)](X_t)dX_t\\
 & = &
 \Psi_{h,h'}(X_T) -\Psi_{h,h'}(X_{t_0}) -
 \frac{1}{2}\int_{t_0}^{T}\psi_{h,h'}'(X_t)\sigma(X_t)^2dt
\end{eqnarray*}
where
\begin{displaymath}
\psi_{h,h'} := K_h\ast ((bf)_{h'}\delta),
\quad {\rm and}\quad
\Psi_{h,h'} :=\mathcal K(./h)\ast ((bf)_{h'}\delta)
\end{displaymath}
is a primitive function of $\psi_{h,h'}$. On the one hand,
\begin{eqnarray*}
 \psi_{h,h'}' & = &
 K_h\ast ((bf)_{h'}\delta)'\\
 & = &
 K_h\ast ((bf)_{h'}\delta') +
 K_h\ast ((K_{h'}\ast (bf)')\delta)\\
 & = &
 K_h\ast ((bf)_{h'}\delta') +
 K_h\ast ((K_{h'}\ast (bf'))\delta) +
 K_h\ast ((K_{h'}\ast (b'f))\delta).
\end{eqnarray*}
Then, since $bf$, $bf'$ and $b'f$ are bounded under Assumption \ref{assumption_bound_derivative_density_x},
\begin{eqnarray*}
 \|\psi_{h,h'}'\|_{\infty}
 & \leqslant &
 \|K_h\|_1\|K_{h'}\ast (bf)\|_{\infty}\|\delta'\|_{\infty} +
 \|K_h\|_1\|K_{h'}\ast (bf')\|_{\infty}\|\delta\|_{\infty} +
 \|K_h\|_1\|K_{h'}\ast (b'f))\|_{\infty}\|\delta\|_{\infty}\\
 & \leqslant &
 \|K\|_{1}^{2}\|bf\|_{\infty}\|\delta'\|_{\infty} +
 \|K\|_{1}^{2}\|bf'\|_{\infty}\|\delta\|_{\infty} +
 \|K\|_{1}^{2}\|b'f\|_{\infty}\|\delta\|_{\infty} <\infty.
\end{eqnarray*}
On the other hand,
\begin{displaymath}
\|\Psi_{h,h'}\|_{\infty}
\leqslant
\|\mathcal K(./h)\|_{\infty}\|(K_{h'}\ast (bf))\delta\|_1
\leqslant
\|\mathcal K\|_{\infty}\|\delta\|_{\infty}
\|K\|_1\|bf\|_1 <\infty.
\end{displaymath}
This concludes the proof because
\begin{displaymath}
(T - t_0)|\langle\Phi_h(X,.),(bf)_{h'}\rangle_{2,\delta}|
\leqslant
2\|\Psi_{h,h'}\|_{\infty} +
\frac{T - t_0}{2}\|\sigma\|_{\infty}^{2}\|\psi_{h,h'}'\|_{\infty}
\quad {\rm a.s.}
\end{displaymath}
%


%
\subsubsection{Proof of Lemma \ref{bound_U_statistics}}
For any $h,h'\in\mathcal H_N$,
\begin{displaymath}
U_{h,h'}(N) =
\sum_{i \not= j}
g_{h,h'}(X^i,X^j)
\end{displaymath}
with, for every $\varphi_1,\varphi_2\in E = C^0([0,T];\mathbb R)$,
\begin{displaymath}
g_{h,h'}(\varphi_1,\varphi_2) :=
\langle\Phi_h(\varphi_1,.) - (bf)_h,\Phi_{h'}(\varphi_2,.) - (bf)_{h'}\rangle_{2,\delta}.
\end{displaymath}
On the one hand, since $\mathbb E(g_{h,h'}(\varphi,X)) = 0$ for every $\varphi\in E$, by Gin\'e and Nickl \cite{GN15}, Theorem 3.4.8, there exists a universal constant $\mathfrak m\geqslant 1$ such that for any $\lambda > 0$, with probability larger than $1 - 5.4e^{-\lambda}$,
\begin{displaymath}
\frac{|U_{h,h'}(N)|}{N^2}
\leqslant\frac{\mathfrak m}{N^2}(\mathfrak c_{h,h'}(N)\lambda^{1/2} +\mathfrak d_{h,h'}(N)\lambda +\mathfrak b_{h,h'}(N)\lambda^{3/2} +\mathfrak a_{h,h'}(N)\lambda^2)
\end{displaymath}
where the constants $\mathfrak a_{h,h'}(N)$, $\mathfrak b_{h,h'}(N)$, $\mathfrak c_{h,h'}(N)$ and $\mathfrak d_{h,h'}(N)$ are defined and controlled later. First, note that
\begin{equation}\label{bound_U_statistics_1}
U_{h,h'}(N) =
\sum_{i\not= j}
(\overline g_{h,h'}(X^i,X^j)
-\widetilde g_{h,h'}(X^i)
-\widetilde g_{h',h}(X^j)
+\mathbb E(\overline g_{h,h'}(X^i,X^j)))
\end{equation}
where, for every $\eta,\eta'\in\mathcal H_N$ and $\varphi_1,\varphi_2,\psi\in E$,
\begin{displaymath}
\overline g_{\eta,\eta'}(\varphi_1,\varphi_2) :=
\langle\Phi_{\eta}(\varphi_1,.),\Phi_{\eta'}(\varphi_2,.)\rangle_{2,\delta}
\quad {\rm and}\quad
\widetilde g_{\eta,\eta'}(\psi) :=
\langle\Phi_{\eta}(\psi,.),(bf)_{\eta'}\rangle_{2,\delta} =
\mathbb E(\overline g_{\eta,\eta'}(\psi,X)).
\end{displaymath}
Let us now control $\mathfrak a_{h,h'}(N)$, $\mathfrak b_{h,h'}(N)$, $\mathfrak c_{h,h'}(N)$ and $\mathfrak d_{h,h'}(N)$:
\begin{itemize}
 \item\textbf{The constant $\mathfrak a_{h,h'}(N)$.} Consider
 \begin{displaymath}
 \mathfrak a_{h,h'}(N) :=
 \sup_{\varphi_1,\varphi_2\in E}
 |g_{h,h'}(\varphi_1,\varphi_2)|.
 \end{displaymath}
 By (\ref{bound_U_statistics_1}), Cauchy-Schwarz's inequality and Lemma \ref{alternative_expression_numerator},
 \begin{eqnarray*}
  \mathfrak a_{h,h'}(N) & \leqslant &
  4\sup_{\varphi_1,\varphi_2\in E}
  |\langle\Phi_h(\varphi_1,.),\Phi_{h'}(\varphi_2,.)\rangle_{2,\delta}|
  \leqslant
  4\left(\sup_{\varphi_1\in E}\|\Phi_h(\varphi_1,.)\|_{2,\delta}\right)
  \left(\sup_{\varphi_2\in E}\|\Phi_{h'}(\varphi_2,.)\|_{2,\delta}\right)\\
  & \leqslant &
  4\left[\frac{6\|\mathcal K\|_{\infty}^{2}}{(T - t_0)^2} +
  \frac{\|\delta\|_{\infty}\|\sigma\|_{\infty}^{4}
  \|K'\|_{2}^{2}}{h^3}\right]^{1/2}
  \left[\frac{6\|\mathcal K\|_{\infty}^{2}}{(T - t_0)^2} +
  \frac{\|\delta\|_{\infty}\|\sigma\|_{\infty}^{4}
  \|K'\|_{2}^{2}}{(h')^3}\right]^{1/2}
  \leqslant\frac{\mathfrak c_1}{h_{0}^{3}}
 \end{eqnarray*}
 with
 \begin{displaymath}
 \mathfrak c_1 =
 4\left[\frac{6\|\mathcal K\|_{\infty}^{2}}{(T - t_0)^2} +
 \|\delta\|_{\infty}\|\sigma\|_{\infty}^{4}\|K'\|_{2}^{2}\right].
 \end{displaymath}
 So, since $(Nh_{0}^{3})^{-1}\leqslant 1$,
 \begin{displaymath}
 \frac{\mathfrak a_{h,h'}(N)\lambda^2}{N^2}
 \leqslant
 \frac{\mathfrak c_1\lambda^2}{N^2h_{0}^{3}}
 \leqslant
 \frac{\mathfrak c_1\lambda^2}{N}.
 \end{displaymath}
 \item\textbf{The constant $\mathfrak b_{h,h'}(N)$.} Consider
 \begin{displaymath}
 \mathfrak b_{h,h'}(N)^2 :=
 N\sup_{\varphi\in E}
 \mathbb E(g_{h,h'}(\varphi,X)^2).
 \end{displaymath}
 By (\ref{bound_U_statistics_1}), Cauchy-Schwarz's inequality and Lemma \ref{alternative_expression_numerator},
 \begin{eqnarray*}
  \mathfrak b_{h,h'}(N)^2 & \leqslant &
  16N\sup_{\varphi\in E}
  \mathbb E(\langle\Phi_h(\varphi,.),\Phi_{h'}(X,.)\rangle_{2,\delta}^{2})\\
  & \leqslant &
  16N\mathbb E(
  \|\Phi_{h'}(X,.)\|_{2,\delta}^{2})
  \sup_{\varphi\in E}
  \|\Phi_h(\varphi,.)\|_{2,\delta}^{2}
  \leqslant
  \frac{\mathfrak c_2\mathfrak m(h')N}{h^3}
  \quad {\rm with}\quad
  \mathfrak c_2 = 4\mathfrak c_1.
 \end{eqnarray*}
 So, for any $\theta\in (0,1)$, since $(Nh_{0}^{3})^{-1}\leqslant 1$,
 \begin{eqnarray*}
  \frac{\mathfrak b_{h,h'}(N)\lambda^{3/2}}{N^2}
  & \leqslant &
  2\left(\frac{\theta}{3\mathfrak m}\right)^{1/2}
  \frac{\mathfrak m(h')^{1/2}}{Nh^{3/2}}\times
  \left(\frac{3\mathfrak m}{\theta}\right)^{1/2}
  \frac{\mathfrak c_{2}^{1/2}\lambda^{3/2}}{N^{1/2}}\\
  & \leqslant &
  \frac{\theta\mathfrak m(h')}{3\mathfrak mN^2h^3} +
  \frac{3\mathfrak c_2\mathfrak m\lambda^3}{\theta N}
  \leqslant
  \frac{\theta\mathfrak m(h')}{3\mathfrak mN} +
  \frac{3\mathfrak c_2\mathfrak m\lambda^3}{\theta N}.
 \end{eqnarray*}
 \item\textbf{The constant $\mathfrak c_{h,h'}(N)$.} Consider
 \begin{displaymath}
 \mathfrak c_{h,h'}(N)^2 :=
 N^2\mathbb E(g_{h,h'}(X^1,X^2)^2).
 \end{displaymath}
 By (\ref{bound_U_statistics_1}) and Lemma \ref{alternative_expression_numerator},
 \begin{eqnarray*}
  \mathfrak c_{h,h'}(N)^2 & \leqslant &
  16N^2\mathbb E(\langle\Phi_h(X^1,.),\Phi_{h'}(X^2,.)\rangle_{2,\delta}^{2})\\
  & \leqslant &
  \mathfrak c_3\mathfrak m(h')N^2
  \quad {\rm with}\quad
  \mathfrak c_3 = 16\mathfrak c_{\ref{alternative_expression_numerator},1}.
 \end{eqnarray*}
 So, as previously,
 \begin{displaymath}
 \frac{\mathfrak c_{h,h'}(N)\lambda^{1/2}}{N^2}
 \leqslant
 \frac{\theta\mathfrak m(h')}{3\mathfrak mN} +
 \frac{3\mathfrak c_3\mathfrak m\lambda}{\theta N}.
 \end{displaymath}
 \item\textbf{The constant $\mathfrak d_{h,h'}(N)$.} Consider
 \begin{displaymath}
 \mathfrak d_{h,h'}(N) :=
 \sup_{(a,b)\in\mathcal A}
 \mathbb E\left[\sum_{i < j}a_i(X^i)b_j(X^j)g_{h,h'}(X^i,X^j)\right],
 \end{displaymath}
 where
 \begin{displaymath}
 \mathcal A :=
 \left\{(a,b) :
 \sum_{i = 1}^{N - 1}\mathbb E(a_i(X^i)^2)\leqslant 1
 \textrm{ and }
 \sum_{j = 2}^{N}\mathbb E(b_j(X^j)^2)\leqslant 1\right\}.
 \end{displaymath}
 By (\ref{bound_U_statistics_1}), Cauchy-Schwarz's inequality, Jensen's inequality and Lemma \ref{alternative_expression_numerator},
 \begin{eqnarray*}
  \mathfrak d_{h,h'}(N) & \leqslant &
  4\sup_{(a,b)\in\mathcal A}\mathbb E\left(
  \sum_{i = 1}^{N - 1}
  \sum_{j = i + 1}^{N}
  |a_i(X^i)b_j(X^j)
  \overline g_{h,h'}(X^i,X^j)|\right)\\
  & \leqslant &
  4N\mathbb E(\langle\Phi_h(X^1,.),\Phi_{h'}(X^2,.)\rangle_{2,\delta}^{2})^{1/2}
  \leqslant
  \mathfrak c_{3}^{1/2}\mathfrak m(h')^{1/2}N.
 \end{eqnarray*}
 So, as previously,
 \begin{displaymath}
 \frac{\mathfrak d_{h,h'}(N)\lambda}{N^2}
 \leqslant
 \frac{\theta\mathfrak m(h')}{3\mathfrak mN} +
 \frac{3\mathfrak c_3\mathfrak m\lambda^2}{\theta N}.
 \end{displaymath}
\end{itemize}
Therefore, there exists a deterministic constant $\mathfrak c_4 > 0$, not depending on $N$, $h$ and $h'$, such that with probability larger than $1 - 5.4e^{-\lambda}$,
\begin{displaymath}
\frac{|U_{h,h'}(N)|}{N^2}
\leqslant
\frac{\theta\mathfrak m(h')}{N} +
\frac{\mathfrak c_4(1 +\lambda)^3}{\theta N}.
\end{displaymath}
In conclusion, with probability larger than $1 - 5.4|\mathcal H_N|e^{-\lambda}$,
\begin{displaymath}
\sup_{h\in\mathcal H_N}
\left\{\frac{|U_{h,h_0}(N)|}{N^2}
-\frac{\theta\mathfrak m(h)}{N}
\right\}
\leqslant\frac{\mathfrak c_4(1 +\lambda)^3}{\theta N}
\end{displaymath}
and
\begin{displaymath}
\sup_{h\in\mathcal H_N}
\left\{\frac{|U_{h,h}(N)|}{N^2}
-\frac{\theta\mathfrak m(h)}{N}
\right\}
\leqslant\frac{\mathfrak c_4(1 +\lambda)^3}{\theta N}.
\end{displaymath}
%


%
\subsubsection{Proof of Lemma \ref{bound_trace_term}}
First, the two following results are used several times in the sequel:
\begin{eqnarray}
 \|(bf)_h\|_{2,\delta}^{2} & \leqslant &
 \|\delta\|_{\infty}\int_{-\infty}^{\infty}\left(\int_{-\infty}^{\infty}
 K_h(y - x)b(y)f(y)dy\right)^2dx
 \nonumber\\
 \label{bound_trace_term_1}
 & \leqslant &
 \|\delta\|_{\infty}\int_{-\infty}^{\infty}b(y)^2f(y)\int_{-\infty}^{\infty}K_h(y - x)^2dxdy
 \leqslant
 \frac{\|\delta\|_{\infty}\|K\|_{2}^{2}\|b^2f\|_1}{h}
\end{eqnarray}
and
\begin{eqnarray}
 \mathbb E(V_h(N)) & = &
 \mathbb E(\|\Phi_h(X,.) - (bf)_h\|_{2,\delta}^{2})
 \nonumber\\
 & = &
 \mathbb E(\|\Phi_h(X,.)\|_{2,\delta}^{2}) +
 \|(bf)_h\|_{2,\delta}^{2} - 2\int_{-\infty}^{\infty}(bf)_h(x)\mathbb E(\Phi_h(X,x))\delta(x)dx
 \nonumber\\
 \label{bound_trace_term_2}
 & = &
 \mathbb E(\|\Phi_h(X,.)\|_{2,\delta}^{2}) -
 \|(bf)_h\|_{2,\delta}^{2}.
\end{eqnarray}
Consider
\begin{displaymath}
v_h(N) := V_h(N) -\mathbb E(V_h(N)) =
\frac{1}{N}\sum_{i = 1}^{N}
(g_h(X^i) -\mathbb E(g_h(X^i)))
\end{displaymath}
with
\begin{displaymath}
g_h(\varphi) :=
\|\Phi_h(\varphi,.) - (bf)_h\|_{2,\delta}^{2}
\textrm{ $;$ }
\forall\varphi\in E.
\end{displaymath}
By Bernstein's inequality, for any $\lambda > 0$, with probability larger than $1 - 2e^{-\lambda}$,
\begin{displaymath}
|v_h(N)|
\leqslant
\sqrt{\frac{2\mathfrak v_h\lambda}{N}}
+\frac{\mathfrak c_h\lambda}{N}
\end{displaymath}
where
\begin{displaymath}
\mathfrak c_h =\frac{\|g_h\|_{\infty}}{3}
\quad {\rm and}\quad
\mathfrak v_h =
\mathbb E(g_h(X)^2).
\end{displaymath}
Moreover, by Inequality (\ref{bound_trace_term_1}) and Lemma \ref{alternative_expression_numerator},
\begin{eqnarray*}
 \mathfrak c_h & = &
 \frac{1}{3}\sup_{\varphi\in E}
 \|\Phi_h(\varphi,.) - (bf)_h\|_{2,\delta}^{2}\leqslant
 \frac{2}{3}\left(\sup_{\varphi\in E}\|\Phi_h(\varphi,.)\|_{2,\delta}^{2} +\|(bf)_h\|_{2,\delta}^{2}
 \right)\\
 & \leqslant &
 \frac{\mathfrak c_1}{h^3}
 \quad {\rm with}\quad
 \mathfrak c_1 =
 \frac{2}{3}\left[
 \frac{6\|\mathcal K\|_{\infty}^{2}}{(T - t_0)^{2}} +
 \|\delta\|_{\infty}\|\sigma\|_{\infty}^{4}\|K'\|_{2}^{2} +
 \|\delta\|_{\infty}\|K\|_{2}^{2}\|b^2f\|_1\right]
\end{eqnarray*}
and, by Inequality (\ref{bound_trace_term_1}), Equality (\ref{bound_trace_term_2}) and Lemma \ref{alternative_expression_numerator},
\begin{eqnarray*}
 \mathfrak v_h & \leqslant &
 \|g_h\|_{\infty}\mathbb E(V_h(N))
 \leqslant
 \frac{3\mathfrak c_1}{h^3}
 (\mathbb E(\|\Phi_h(X,.)\|_{2,\delta}^{2}) -
 \|(bf)_h\|_{2,\delta}^{2})\\
 & \leqslant &
 \frac{\mathfrak c_2\mathfrak m(h)}{h^3}
 \quad {\rm with}\quad
 \mathfrak c_2 = 3\mathfrak c_1.
\end{eqnarray*}
Then, for any $\theta\in (0,1)$, since $(Nh_{0}^{3})^{-1}\leqslant 1$, with probability larger than $1 - 2e^{-\lambda}$,
\begin{eqnarray*}
 |v_h(N)|
 & \leqslant &
 2\sqrt{\frac{\mathfrak c_2\mathfrak m(h)\lambda}{Nh^3}}
 +\frac{\mathfrak c_1\lambda}{Nh^3}\\
 & \leqslant &
 \theta\mathfrak m(h) +
 \frac{(\mathfrak c_1 +\mathfrak c_2)\lambda}{\theta Nh^3}
 \leqslant
 \theta\mathfrak m(h) +
 \frac{(\mathfrak c_1 +\mathfrak c_2)\lambda}{\theta}.
\end{eqnarray*}
So, with probability larger than $1 - 2|\mathcal H_N|e^{-\lambda}$,
\begin{displaymath}
\sup_{h\in\mathcal H_N}
\left\{\frac{|v_h(N)|}{N}
-\frac{\theta\mathfrak m(h)}{N}\right\}
\leqslant
\frac{(\mathfrak c_1 +\mathfrak c_2)\lambda}{\theta N}.
\end{displaymath}
Therefore, by Equality (\ref{bound_trace_term_2}), with probability larger that $1 - 2|\mathcal H_N|e^{-\lambda}$,
\begin{eqnarray*}
 & &
 \sup_{h\in\mathcal H_N}\left\{
 \frac{1}{N}|V_h(N) -\mathbb E(\|\Phi_h(X,.)\|_{2,\delta}^{2})| -\frac{\theta\mathfrak m(h)}{N}\right\}\\
 & &
 \hspace{1.5cm}
 \leqslant
 \sup_{h\in\mathcal H_N}
 \left\{\frac{|v_h(N)|}{N}
 -\frac{\theta\mathfrak m(h)}{N}\right\} +
 \frac{1}{N}\|K_h\ast (bf)\|_{2,\delta}^{2}
 \leqslant
 \frac{(\mathfrak c_1 +\mathfrak c_2 +\|\delta\|_{\infty}
 \|K\|_{1}^{2}\|bf\|_{2}^{2})(1 +\lambda)}{\theta N}.
\end{eqnarray*}
%


%
\subsubsection{Proof of Lemma \ref{bound_crossed_term}}
For any $h,h'\in\mathcal H_N$,
\begin{displaymath}
W_{h,h'}(N) =
\frac{1}{N}\sum_{i = 1}^{N}
(g_{h,h'}(X^i) -\mathbb E(g_{h,h'}(X^i)))
\end{displaymath}
with, for every $\varphi\in E$,
\begin{displaymath}
g_{h,h'}(\varphi) :=
\langle\Phi_h(\varphi,.),(bf)_{h'} - bf\rangle_{2,\delta}.
\end{displaymath}
By Bernstein's inequality, for any $\lambda > 0$, with probability larger than $1 - 2e^{-\lambda}$,
\begin{displaymath}
|W_{h,h'}(N)|\leqslant
\sqrt{\frac{2\mathfrak v_{h,h'}\lambda}{N}}
+\frac{\mathfrak c_{h,h'}\lambda}{N}
\end{displaymath}
where
\begin{displaymath}
\mathfrak c_{h,h'} =\frac{\|g_{h,h'}\|_{\infty}}{3}
\quad {\rm and}\quad
\mathfrak v_{h,h'} =
\mathbb E(g_{h,h'}(X)^2).
\end{displaymath}
Moreover, by Lemma \ref{alternative_expression_numerator},
\begin{eqnarray*}
 \mathfrak c_{h,h'} & = &
 \frac{1}{3}\sup_{\varphi\in E}
 |\langle\Phi_h(\varphi,.),(bf)_{h'} - bf\rangle_{2,\delta}|
 \leqslant
 \frac{1}{3}\|(bf)_{h'} - bf\|_{2,\delta}
 \sup_{\varphi\in E}\|\Phi_h(\varphi,.)\|_{2,\delta}\\
 & \leqslant &
 \frac{\mathfrak c_1}{h^{3/2}}
 \|(bf)_{h'} - bf\|_{2,\delta}
 \quad {\rm with}\quad
 \mathfrak c_1 =
 \frac{1}{3}\left[\frac{6\|\mathcal K\|_{\infty}^{2}}{(T - t_0)^{2}} +
 \|\delta\|_{\infty}\|\sigma\|_{\infty}^{4}\|K'\|_{2}^{2}\right]^{1/2}
\end{eqnarray*}
and
\begin{displaymath}
\mathfrak v_{h,h'}\leqslant
\mathbb E(\langle\Phi_h(X,.),(bf)_{h'} - bf\rangle_{2,\delta}^{2})\\
\leqslant\mathfrak c_{\ref{alternative_expression_numerator},2}
\|(bf)_{h'} - bf\|_{2,\delta}^{2}.
\end{displaymath}
Then, for any $\theta\in (0,1)$, with probability larger than $1 - 2e^{-\lambda}$,
\begin{eqnarray*}
 |W_{h,h'}(N)|
 & \leqslant &
 2\sqrt{\frac{\mathfrak c_{\ref{alternative_expression_numerator},2}\lambda}{N}
 \|(bf)_{h'} - bf\|_{2,\delta}^{2}}
 +\frac{\mathfrak c_1\lambda}{Nh^{3/2}}\|(bf)_{h'} - bf\|_{2,\delta}\\
 & \leqslant &
 \theta\|(bf)_{h'} - bf\|_{2,\delta}^{2} +
 \frac{2\mathfrak c_{\ref{alternative_expression_numerator},2}\lambda}{\theta N} +
 \frac{2\mathfrak c_{1}^{2}\lambda^2}{\theta N^2h^3}
 \leqslant
 \theta\|(bf)_{h'} - bf\|_{2,\delta}^{2} +
 \frac{2(\mathfrak c_{\ref{alternative_expression_numerator},2} +
 \mathfrak c_{1}^{2})(1 +\lambda)^2}{\theta N}.
\end{eqnarray*}
So, with probability larger than $1 - 2|\mathcal H_N|e^{-\lambda}$,
\begin{eqnarray*}
 \sup_{h\in\mathcal H_N}\{
 |W_{h,h_0}(N)| -\theta\|(bf)_{h_0} - bf\|_{2,\delta}^{2}\}
 & \leqslant &
 \frac{\mathfrak c_{\ref{bound_crossed_term}}(1 +\lambda)^2}{\theta N},\\
 \sup_{h\in\mathcal H_N}\{
 |W_{h_0,h}(N)| -\theta\|(bf)_h - bf\|_{2,\delta}^{2}\}
 & \leqslant &
 \frac{\mathfrak c_{\ref{bound_crossed_term}}(1 +\lambda)^2}{\theta N}\textrm{ and}\\
 \sup_{h\in\mathcal H_N}\{
 |W_{h,h}(N)| -\theta\|(bf)_h - bf\|_{2,\delta}^{2}\}
 & \leqslant &
 \frac{\mathfrak c_{\ref{bound_crossed_term}}(1 +\lambda)^2}{\theta N}.
\end{eqnarray*}
%


%
\subsection{Proof of Corollary \ref{risk_bound_NW_PCO}}
On the one hand, as in the proof of Proposition \ref{risk_bound_NW} and since $\delta(x) > m$ for every $x\in [A,B]$,
\begin{eqnarray*}
 \mathbb E(\|\widehat b_{N,\widehat h,\widehat h'} - b\|_{f,A,B}^{2})
 & \leqslant &
 \frac{\mathfrak c_{\ref{risk_bound_NW}}}{m^2}[
 \mathbb E(\|\widehat{bf}_{N,\widehat h} - bf\|_{2,A,B}^{2})
 + 2\mathbb E(\|\widehat f_{N,\widehat h'} - f\|_{2}^{2})]\\
 & \leqslant &
 \frac{2\mathfrak c_{\ref{risk_bound_NW}}}{m^3}[
 \mathbb E(\|\widehat{bf}_{N,\widehat h} - bf\|_{2,\delta}^{2})
 +\mathbb E(\|\widehat f_{N,\widehat h'} - f\|_{2}^{2})].
\end{eqnarray*}
On the other hand, by Theorem \ref{risk_bound_PCO_estimator_bf} and {\it union bounds},
\begin{eqnarray*}
 \mathbb E(\|\widehat{bf}_{N,\widehat h} - bf\|_{2,\delta}^{2})
 & \leqslant &
 (1 +\vartheta)\min_{h\in\mathcal H_N}
 \mathbb E(\|\widehat{bf}_{N,h} - bf\|_{2,\delta}^{2}) +
 \frac{\mathfrak c_{\ref{risk_bound_PCO_estimator_bf},2}}{\vartheta}
 \left(\|(bf)_{h_0} - bf\|_{2,\delta}^{2} +\frac{1}{N}\right)\\
 & \leqslant &
 (1\vee\|\delta\|_{\infty})\left[
 (1 +\vartheta)\min_{h\in\mathcal H_N}
 \mathbb E(\|\widehat{bf}_{N,h} - bf\|_{2}^{2}) +
 \frac{\mathfrak c_{\ref{risk_bound_PCO_estimator_bf},2}}{\vartheta}
 \left(\|(bf)_{h_0} - bf\|_{2}^{2} +\frac{1}{N}\right)\right]
\end{eqnarray*}
and
\begin{displaymath}
\mathbb E(\|\widehat f_{N,\widehat h'} - f\|_{2}^{2})
\leqslant
(1 +\vartheta)\min_{h'\in\mathcal H_N}
\mathbb E(\|\widehat f_{N,h'} - f\|_{2}^{2}) +
\frac{\overline{\mathfrak c}_{\ref{risk_bound_PCO_estimator_bf},2}}{\vartheta}
\left(\|f_{h_0'} - f\|_{2}^{2} +\frac{1}{N}\right).
\end{displaymath}
Therefore,
\begin{eqnarray*}
 \mathbb E(\|\widehat b_{N,\widehat h,\widehat h'} - b\|_{f,A,B}^{2})
 & \leqslant &
 \frac{2\mathfrak c_{\ref{risk_bound_NW}}(1\vee\|\delta\|_{\infty})}{m^3}\left[
 (1 +\vartheta)\min_{(h,h')\in\mathcal H_N\times\mathcal H_N'}
 \{\mathbb E(\|\widehat{bf}_{N,h} - bf\|_{2}^{2}) +
 \mathbb E(\|\widehat f_{N,h'} - f\|_{2}^{2})\}\right.\\
 & &
 \hspace{4cm}\left.
 +\frac{\mathfrak c_{\ref{risk_bound_NW_PCO}}}{\vartheta}\left(
 \|(bf)_{h_0} - bf\|_{2}^{2} +\|f_{h_0'} - f\|_{2}^{2} +\frac{1}{N}
 \right)\right].
\end{eqnarray*}
{\bf Acknowledgments.} Many thanks to Fabienne Comte and Valentine Genon-Catalot for their constructive comments about this work.
\end{document}